\documentclass[10pt, reqno]{amsart}

\usepackage{epigraph}

\usepackage[pagewise]{lineno}
\usepackage{tikz}
\usepackage{curves}
\usepackage[widespace]{fourier}
\usepackage{xcolor}
\usepackage{hyperref}
\usepackage[utf8]{inputenc}
\hypersetup{colorlinks=true,linkcolor=blue, citecolor=blue}

\usepackage{amsmath,amsthm,amscd,amsfonts,amssymb,euscript}
\usepackage{graphics}
\usepackage{tikz-cd}
\usepackage{pdfpages}
\usepackage{mathtools}

\usepackage{frcursive}
\usepackage[all]{xypic}
\usepackage{latexsym}
\usepackage{color}
\theoremstyle{plain}

\font\manual=manfnt
\def\dbend{{\manual\char127}} 

\def\danger{\begin{trivlist}\begin{footnotesize}\item[]\noindent%
\begingroup\hangindent=3pc\hangafter=-2
\def\par{\endgraf\endgroup}%
\hboX to0pt{\hskip-\hangindent\dbend\hfill}\ignorespaces}
\def\enddanger{\par\end{footnotesize}\end{trivlist}}

\def\ddanger{\begin{trivlist}\begin{footnotesize}\item[]\noindent%
\begingroup\hangindent=3pc\hangafter=-2
\def\par{\endgraf\endgroup}%
\hboX to0pt{\hskip-\hangindent\dbend\kern2pt\dbend\hfill}\ignorespaces}
\def\endddanger{\par\end{footnotesize}\end{trivlist}}

\DeclareFontFamily{OT1}{rsfs}{}
\DeclareFontShape{OT1}{rsfs}{n}{it}{<-> rsfs10}{}
\DeclareMathAlphabet{\mathscr}{OT1}{rsfs}{n}{it}

\DeclareFontFamily{U}{stixbbit}{}
\DeclareFontShape{U}{stixbbit}{m}{it}{<-> stix-mathbbit}{}

\newcommand{\Hom}{{\mathrm {Hom}}}

\newcommand{\eP}{{\mathscr P}}
\newcommand{\eE}{{\mathscr E}}
\newcommand{\eG}{{\mathscr G}}

\newcommand{\eV}{{\mathscr V}}
\newcommand{\eW}{{\mathscr W}}

\newcommand{\eT}{{\mathscr T}}

\newcommand{\eR}{{\mathscr R}}

\newcommand{\qh}{\tiny\text{\cursive H}}

\newcommand{\tC}{\tilde{C}}
\newcommand{\tZ}{\tilde{Z}}
\newcommand{\tW}{\tilde{W}}
\newcommand{\tU}{\tilde{U}}
\newcommand{\tN}{\tilde{N}}
\newcommand{\tsl}{\text{SL}}
\newcommand{\tgl}{\text{GL}}
\newcommand{\cug}{\text{\cursive g}}

\makeatletter
\let\@wraptoccontribs\wraptoccontribs
\makeatother

\begin{document}
\input{amssym.def}

\newtheorem*{theorem*}{Theorem}

\newtheorem{guess}{\sc Theorem}[section]
\newcommand{\bth}{\begin{guess}$\!\!\!${\bf }~}
\newcommand{\eeth}{\end{guess}}

\newtheorem{propo}[guess]{\sc Proposition}

\newcommand{\bprop}{\begin{propo}$\!\!\!${\bf }~}
\newcommand{\eprop}{\end{propo}}

\newtheorem{lema}[guess]{\sc Lemma}
\newcommand{\blem}{\begin{lema}$\!\!\!${\bf }~}
\newcommand{\elem}{\end{lema}}

\newtheorem{defe}[guess]{\sc Definition}
\newcommand{\bdefe}{\begin{defe}$\!\!\!${\it }~}
\newcommand{\edefe}{\end{defe}}

\newtheorem{coro}[guess]{\sc Corollary}
\newcommand{\bcor}{\begin{coro}$\!\!\!${\bf }~}
\newcommand{\ecor}{\end{coro}}

\newtheorem{rema}[guess]{\it Remark}
\newcommand{\brem}{\begin{rema}$\!\!\!${\it }~\rm}
\newcommand{\erem}{\end{rema}}

\theoremstyle{remark}
\newtheorem{assump}[guess]{\sc Assumption}

\newcommand{\spec}{{\rm Spec}\,}
\newtheorem{notation}{Notation}[section]
\newcommand{\bnot}{\begin{notation}$\!\!\!${\bf }~~\rm}
\newcommand{\enot}{\end{notation}}

\newcommand{\bpr}{\begin{proof}}
\newcommand{\epr}{\end{proof}}
\numberwithin{equation}{guess}

\newcommand{\beqa}{\begin{eqnarray}}
\newcommand{\eeqa}{\end{eqnarray}}
\theoremstyle{definition}
\newtheorem{example}[guess]{\it Example}
\newtheorem{say}[guess]{\it }
\newcommand{\bsem}{\begin{say}$\!\!\!${\it }~~\rm}
\newcommand{\esem}{\end{say}}
\newtheorem{observe}[subsubsection]{Observation}
\newsymbol \bulletledarrowleft 1309

\newcommand{\ha}{\sf h}
\newcommand{\g}{\sf g}
\newcommand{\ta}{\mathfrak T}
\newcommand{\s}{\sf s}

\newcommand{\wt}{\widetilde}
\newcommand{\Lr}{\Longrightarrow}
\newcommand{\Aut}{\mboX{{\rm Aut}$\,$}}
\newcommand{\ul}{\underline}
\newcommand{\ol}{\overline}
\newcommand{\lr}{\longrightarrow}
\newcommand{\sh}{{\sf h}}
\newcommand{\smu}{{\sf\mu_{_d}}}
\newcommand{\be}{{\mathbb E}}
\newcommand{\ba}{{\mathbb A}}
\newcommand{\bd}{{\mathbb D}}
\newcommand{\bc}{{\mathbb C}}
\newcommand{\bp}{{\mathbb P}}
\newcommand{\bz}{{\mathbb Z}}
\newcommand{\bq}{{\mathbb Q}}
\newcommand{\bn}{{\mathbb N}}
\newcommand{\bm}{{\em}}
\newcommand{\bg}{{\mathbb G}}
\newcommand{\bbf}{{\mathbb F}}
\newcommand{\br}{{\mathbb R}}
\newcommand{{\bh}}{{\mathbb H}}
\newcommand{\bo}{{\overline \omega}'}
\newcommand{\po}{{\omega}'}
\newcommand{\dps}{\tt {DP}}

\newcommand{\ct}{{\mathcal T}}
\newcommand{\cc}{{\mathcal C}}
\newcommand{\cd}{{\mathcal D}}
\newcommand{\cl}{{\mathcal L}}
\newcommand{\cv}{{\mathcal V}}
\newcommand{\cf}{{\mathcal F}}
\newcommand{\cb}{{\Lambda}}
\newcommand{\ch}{{\mathcal H}}
\newcommand{\cw}{{\mathcal W}}
\newcommand{\bba}{{\hat\mathbb A}}

\newcommand{\mfc}{{\sf C}}
\newcommand{\ce}{{\mathcal E}}
\newcommand{\co}{{\mathcal O}}

\newcommand{\cm}{{\mathcal M}}

\newcommand{\cs}{{\mathcal S}}
\newcommand{\cg}{{\mathcal G}}
\newcommand{\ca}{{\mathcal A}}
\newcommand{\hra}{\hookrightarrow}
\newcommand{\mfu}{{\sf U}}

\newtheorem{ack}{\it Acknowledgments}       
\renewcommand{\theack}{} 
\title{Torsors on semistable curves and degenerations}
\author{V. Balaji} 
\address{{Chennai Mathematical Institute SIPCOT IT Park, Siruseri-603103, India,
balaji@cmi.ac.in}}
\date{Revised version,22 January, 2021}
\dedicatory{Dedicated to the memory of C.S. Seshadri.}
\begin{abstract} In this paper we answer two long-standing questions on the classification of $G$-torsors on curves for an almost simple, simply connected algebraic group $G$ over the field of complex numbers. The first is the construction of a flat degeneration of the moduli of  $G$-torsors on smooth projective curves when the smooth curve degenerates to an irreducible nodal curve  and the second one is to give an intrinsic definition of (semi)stability for a $G$-torsor on an {\em irreducible projective nodal curve}. A generalization of the classical Bruhat-Tits group schemes to two-dimensional regular local rings and an application of the geometric formulation of the McKay correspondence provide the key tools.
\end{abstract}
\maketitle
\footnotesize
\tableofcontents
\normalsize

\section{Introduction}

Let $G$ be an {\em almost simple and simply
connected} algebraic group  or the linear group $\text{GL}(n)$, over the field $k = \mathbb C$ of complex numbers. Let $A = \spec~ k\llbracket t \rrbracket$ and $K = \spec k(\!(t)\!)$. Let  $o \in A$ be the closed point and let $C_{_{A}} \to A$ be a proper, flat family with generic fibre $C_{_K}$ a smooth projective curve of genus $\text{g} \geq 1$ and closed fibre $C_{_o}$ an irreducible nodal curve $(C,c)$ with a single node $c \in C$. Let $\text{Bun}_{_G}(C_{_K})$ denote the stack of $G$-bundles on the curve $C_{_K}$. These stacks do not satisfy the valuative criterion for properness  and one needs to impose suitable semistability conditions to get a  separated Artin stack with a coarse space which is the proper moduli space of "slope" semistable principal $G$-bundles. These moduli spaces were constructed by A. Ramanathan in 1975 \cite{ram2}. The first examples of these spaces is in the case when $G = \text{GL}(n)$ which give the Mumford-Seshadri moduli spaces of (semi)stable vector bundles (where we need to fix the degree of the bundles).

The primary purpose of the present work is two-fold ({\em see \ref{ideainoutline} below for an outline of the basic idea in the work}). In the first part we 
construct a flat degeneration of the stack $\text{Bun}_{_G}(C_{_K})$.  This question has remained open due to the lack of a suitable analogue for the notion of torsion-free sheaves on curves in the realm of $G$-bundles (see \cite[Page 489]{faltings} and \cite[page 347]{faltings2}). 

The approach in this paper is to replace the node by a bubbling. This replaces the nodal curve $C$ by certain semistable curves $C^{^{(d)}}$ whose stable model is $C$.  We first construct a degeneration of the group $G$ to a parahoric group scheme over the semistable limits and then we define the limiting objects as certain torsors for these group schemes (see \S3 for details). Thus, torsion-free sheaves on $C$ are replaced by triples $(C^{^{(d)}}, \qh, \ce)$. This basic idea in the case of $\tgl(n)$ goes back to the work of Gieseker \cite{gies} (for $GL(2)$) and those of Nagaraj-Seshadri and Ivan Kausz for $\tgl(n), n > 2$. The novelty of the present approach is in relating this to torsors under parahoric degenerations of $G$.

The key idea for this generalization is the extension of the notion of Bruhat-Tits group schemes to the setting of regular $2$-dimensional local rings (see \S3). This local construction is based on a close study of the geometric McKay correspondence  as in \cite{gonverd}. Loosely put, this allows us to set up a kind of Fourier-Mukai transform for group schemes and torsors in two steps. In the first step, one begins with an equivariant group scheme with fibre isomorphic to $G$ on the affine plane. By a Fourier-Mukai-like operation, we construct an affine ``parahoric" group scheme on the minimal resolution of singularities of the  analytic normal surface $k\llbracket x,y,t\rrbracket/(xy - t^{^{d}})$. The second step is then to set-up a Fourier-Mukai   between the category of certain equivariant $G$-torsors on the affine plane and the category of torsors for the parahoric group scheme on the minimal resolution of singularities of the  analytic normal surface $k\llbracket x,y,t\rrbracket/(xy - t^{^{d}})$. To build the stack, the next hard step is to build a  parahoric group scheme on certain standard models for degenerations of smooth curves using the ones already constructed on the minimal resolutions for the surfaces $k\llbracket x,y,t\rrbracket/(xy - t^{^{d}})$.  Here we rely heavily on the work of Jun Li \cite{Li} on expanded degenerations and using these we build local models for Gieseker bundles (see \S4). Once this is done, one can build the stack of Gieseker-torsors for the parahoric group schemes and prove the relevant properties using standard techniques \eqref{flatavecnc} \eqref{algebraicstack}. The new idea is now to view these as logarithmic schemes. Bundles with parabolic structures on the generic points of the normal crossing divisors appear naturally and give shape to the objects.  

The resulting moduli stack over $A := \spec k\llbracket t \rrbracket$ has a closed fibre which is a s.n.c.divisor with $\ell + 1$ smooth components, $\ell = \text{rank}(G)$, indexed by the extended Dynkin diagram  of $G$. This degeneration is therefore a semistable degeneration in the sense of Mumford \cite{kkms} (see also \cite{balapandey}). Towards the very end of  \cite{kkms} Mumford constructs a relative compactification $\bar{G}_{_A}$ of $G \times A$ where $\bar{G}_{_K} = G_{_K}$. The closed fibre $\bar{G}_{_o}$ is a union of complete varieties meeting with simple normal crossing singularities. These are indexed by the vertices of the affine Dynkin diagram. This is an algebro-geometric model of the Bruhat-Tits building in the relative case. He remarks at the close of the book that his compactification can be viewed as a kind of ``N\'eron model with corners" of the semisimple group scheme over the local ring.   The degeneration in the present work can be viewed as a precise analogue for $\text{Bun}_{_G}(C_{_A})$.  

The stack, local analytically, gives a resolution of singularities  for analogues of certain matrix type singularities. These  occur on the stack of torsion-free sheaves and their links to the theory of local models and PEL Shimura varieties were already seen by Faltings \cite{faltings1}. This we hope would make our stacks  wider in their appeal. We describe in outline the structure of the closed fibre of the stack and elaborate it in the case when $G = \tsl(2)$ \eqref {components} \eqref{components1}.

In Part II and III, we work towards the coarse space. To get a  separated and proper stack in the limit  we require a definition of semistability of certain torsors under the parahoric group scheme on semistable curves. Here, even the basic case of a principal $G$-bundle on an irreducible nodal curve itself was not well understood.

A notion of $\mu$-semistability   which appropriately generalizes Ramanathan's definition, has been open  and presents serious difficulties (see \cite{faltings}, \cite{fm3}). In the third of a series of papers (\cite{fm2}, \cite{fm3}, \cite{fmw}) on principal $G$-bundles on elliptic curves and singular curves,  R. Friedman and J. Morgan write that ``{\em  there are many remaining open questions. One of the deepest is the problem of finding an intrinsic definition of semistability for G-bundles on a singular curve, and of a generalized
form of S-equivalence, which would be broad enough to include those bundles coming from
the parabolic construction}". Moreover, for the construction of coarse spaces, this notion should have a GIT interpretation.

For this study, we concentrate on a single nodal curve. The first point is to recognize that any definition of semistability of a principal $G$-bundle on an irreducible nodal curve would require one to already expand the notion of a $G$-bundle to our torsors for the parahoric group scheme. This phenomenon shows up even for vector bundles where the test objects for semistability could be torsion-free sheaves. As a step towards achieving this, we express $C$ as a coarse space of a twisted curve $\cc_{_d}$ in the sense of \cite{ollson} and express the semistability condition on torsion-free sheaves on $C$ in terms of torsors on 
$\cc_{_d}$. To get a good notion of "degree" of line bundles on these Deligne-Mumford stacks we set up a "Fourier-Mukai" like correspondence between torsors on $\cc_{_d}$ and certain objects, which we term {\em laced} torsors, on the normalization $\tC$ of $C$. These are parahoric torsors on $\tC$ with parahoric structures at the two points above the node, along with a "descent datum". The notion of a (semi)stability of $G$-torsors on $\cc_{_d}$, which is "equivalent" to the (semi)stability of torsion-free sheaves when $G = \text{GL}(n)$, is then achieved \eqref{newtfss}. We term this notion {\tt tf}-(semi)stability. The task then onwards is to show that this notion is a GIT notion \eqref{etassandschmittss} which defines  a good moduli problem. Here we draw on the work of A.Schmitt \cite{schmitt2} to express  our moduli problem in his terms and solve it using GIT methods. 

Finally, we relate the Gieseker torsors on semistable curves $C^{^{(d)}}$ to laced torsors on $\tC$ by restricting the torsor to $\tC$  and then use this to define {\tt tf}-(semi)stability of Gieseker torsors. Note that at the back of this notion is the fact that there is a "morphism" from the stack of Giesker-torsors to a "virtual" space analogous to the space of torsion-free sheaves. This needs to be carefully placed on a rigorous footing and Schmitt's construction plays the role for this. 

Using the {\tt tf}-(semi)stability and a relative polarization for this morphism, we invoke a classical principle due to Seshadri to finally get a more refined notion, that of {\tt L}-(semi)stability for Gieseker torsors. The GIT approach then shows how these notions give the construction of a coarse space for the open substack of semistable Gieseker torsors.  We summarize the main results of the paper in the following:
\begin{theorem*} \big(see \eqref{algebraicstack},\eqref{flatavecnc}, \eqref{the big theorem}\big)
\begin{enumerate}
\item The stack ${\text{Gies}}_{_{G}}(C_{_A})$ of  Gieseker torsors \eqref{giesekertorsor}  is an   algebraic stack locally of finite type, which is regular and flat over $A$. Over $K$ we have an identification ${\text{Gies}}_{_{G}}(C_{_K}) = \text{Bun}_{_G}(C_{_K})$ with the stack of  $G$-torsors on the smooth projective curve $C_{_K}$. Further the closed fibre ${\text{Gies}}_{_{G}}(C_{_o}) \subset {\text{Gies}}_{_{G}}(C_{_A})$ is a divisor with normal crossings with $\ell +1$ smooth components indexed by the extended Dynkin diagram. 
\item  The open substack ${\text{Gies}}_{_{G}}(C_{_A})^{^{{\tt L}-ss}}$ of {\tt L}-(semi)stable Gieseker torsors \eqref{elss}  has a coarse space which parametrizes $S$-equivalence classes of Gieseker torsors and which provides a {\em proper flat degeneration} of the moduli scheme of $\mu$-(semi)stable $G$-torsors on $C_{_K}$.
\end{enumerate}
\end{theorem*}
The layout of the paper is as follows. In \S3 we make the basic construction using the McKay correspondence. In \S4, we construct the group objects on Jun Li's standard models. In \S5 we discuss admissibility for the $\text{GL}(n)$ case and in \S6 we discuss the general $G$ case. In \S7 we prove the stack-theoretic properties and in \S8 we look at some examples and describe the closed fibre of the degeneration. From \S9 till \S12 we work with a single nodal curve and define the semistability of Gieseker torsors on semistable curves. In \S13 we complete the construction of the degeneration of the moduli space of $G$-bundles.
\subsubsection{The basic idea in outline}\label{ideainoutline}
We work with one parameter degenerations and with the family of curves $C_{_A}$. In the discussion below, by the {\em limiting fibre} we mean objects over the closed fibre $C_{_o} = (C,c)$. To tackle the degeneration problem for vector bundles, in the early eighties, Gieseker and Seshadri approached the problem in two different ways. Seshadri in \cite{sesasterisque} took the  approach of degenerating (slope semistable) vector bundles to torsion-free sheaves on the nodal limit $(C,c)$, an approach  which was initiated in a paper by Mayer-Mumford and Oda-Seshadri in the degeneration of the Picard variety. Seshadri's limiting moduli space was slope semistable torsion-free sheaves of fixed rank and degree and his strategy was GIT. Gieseker (\cite{gies}) approached the problem again by GIT but his strategy had its seed in his approach to the construction of moduli of bundles on curves and surfaces. This was by studying smooth curves embedded in Grassmannians using vector bundles generated by sections. The Hilbert scheme of such curves had a natural action of a suitable linear group and GIT on the Hilbert scheme gave rise to "semistable" curves embedded in the Grassmannian; these were limits of smooth curves. This was Gieseker's GIT construction of the coarse space for Deligne-Mumford compactification of $\bar{\mathcal M}_{_g}$.
The tautological vector bundle on the Grassmannian when restricted to the limiting curves  gave the "semistable" objects in the problem. The limiting moduli over the nodal curve was then classified as a "list" of semistable curves $C^{^{(d)}}$ (with a fixed stable model $(C,c)$) together with a class of vector bundles on them. These bundles on the chain of $\bp^{^1}$'s in $C^{^{(d)}}$, were from a fixed list of vector bundles which Gieseker called "standard"; they were bundles whose direct summands on each $\bp^{^1}$ had only $\co$ or $\co(1)$. Gieseker's approach had an added feature, viz, the "stack" of objects he obtained (in the modern language) was regular over $k$ and the limiting fibre was reduced with normal crossing singularities.

Gieseker's approach faced a serious block in going to higher rank bundles, since identifying the semistable limits by GIT became very unwieldy when the rank of the vector bundle exceeded $2$. Nagaraj-Seshadri \cite{ns2} and Kausz \cite{kausz}  combined the two approaches, namely Seshadri's and Gieseker's, to solve the higher rank degeneration problem with s.n.c property. En route, they obtain the {\em standard list} by two bits of data (this is my interpretation): (1) a local data, i.e. the local types of the torsion-free sheaves on $(C,c)$, which was encoded in the number of summands of the maximal ideal and (2) a data of realizing a torsion-free sheaf on $(C,c)$ as  limit of vector bundles on $C_{_K}$. The idea of {\em bubbling} was then to {\em blow-up} the torsion-free sheaves on the surface gotten by base changing $C_{_A}$ by $t \mapsto t^{^d}$. One obtains new one-parameter families of curves where the original family of vector bundles now had vector bundles as limits. The valuative criterion for the functor needs to account for such base change. The two bits of data gave the Gieseker-type list of semistable curves $C^{^{(d)}}$, together with "standard" vector bundles on them which came with a configuration of $\co$'s and $\co(1)$'s on the chain of rational curves. 

The problem for $G$-bundles and their degeneration was that,  on irreducible nodal curves there was no satisfactory solution \`a la Seshadri. More precisely, there was no {\em torsion-free} analogue except in the classical case of the symplectic and orthogonal cases both of which were exploited by Faltings \cite{faltings}; in either case, there is a basic representation and the problem gets resolved as one on torsion-free sheaves equipped with degenerate forms.  My approach is to work around this lacuna and get to the bubbling directly, i.e. by circumventing the torsion-free route. The basic principle was to identify the two bits of data for "possible limits" of principal $G$-bundles. The list was essentially {\em local data} on the $\bp^{^1}$-chains which were then glued to $G$-bundles on the normalization $\tC$ of $C$ at two marked end points of the chain. 

In the 1980's, Gonzalez-Springberg and Verdier \cite{gonverd} and Artin-Verdier \cite{artinverd} studied reflexive sheaves on normal surface singularities in the context of a geometric  McKay corresondence. These objects had been studied in depth in the paper by Lipman \cite{lipman}. From my standpoint, the paper \cite{gonverd} gives an alternate approach to reaching the bubbling data and vector bundles on chains. This was done by simultaneously considering the minimal resolution of the local normal singularities and viewing them also as quotient singularities by actions of finite Kleinian groups on affine planes. The minimal resolution of singularities ${\sf N}^{^{(d)}}$ (see Figure \eqref{minres}),of the normal singularity was realized as a "minimal platificateur" (see  \eqref{platif} and \cite[Corollaire 7, page 448]{gonverd}). The bundles were obtained by a  "Fourier-Mukai" from equivariant bundles on the affine plane to bundles on the minimal resolutions; the scheme  $D^{^{(d)}}$ in \eqref{keydiag1} acts as a correspondence  for a pull-back and an invariant push-forward from $D$ to ${\sf N}^{^{(d)}}$, indeed ${\sf N}^{^{(d)}}$ can be identified with a certain "Hilbert scheme" classifying equivariant zero-cycles on $D$, and $D^{^{(d)}}$ the universal space (work of Ito-Nakamura).

In my paper with Seshadri \cite{base}, we had studied the principal bundle analogue of parabolic vector bundles. The parahoric group schemes  were realized via what we termed "invariant direct images" from equivariant affine group schemes on ramified covers, or more precisely "orbifold stacks". Invariant direct images of group schemes were simply taking Weil restrictions of scalars under Galois coverings and then taking invariants by the Galois group; this process works well in characteristic zero and also in the "tame" cases. 

How does all this come together in the degeneration question? The idea is to first get to the "basic list" by a Fourier-Mukai like construction of affine group schemes on minimal resolution of singularities of normal surface singularities of type $A_{_d}$; these singularities were simply $ {\ba}^{^2}/\smu$. The basic list is essentially local analytic in its content when viewed on the regular surface, being data along the rational curve-chains. The new group schemes, which we term {\sf 2}BT group schemes, comes by the following process. 

We begin by considering equivariant affine group schemes under the action of $\smu$ on $D$ (which was the analytic disc at the origin $0$ in $\ba^{^2}$). The $d$ are allowed to vary. The action is essentially  given by the data of conjugacy class of representations $\rho:\smu \to G$ which we call "type $\tau$" following an old terminology due to Weil-Seshadri. First take the trivial $G$-bundle on $D$ with a twisted action by $\smu$, i.e. $\ct_{_D} := D \times ^{^{\rho}} G$ and then take the "adjoint group scheme" $\ct_{_D}(G)$ \eqref{basictorsor}, where $G$ acts on itself by inner conjugation. This gives the basic equivariant group schemes for each local type $\tau$. We then perform a "Fourier-Mukai" to these group schemes to get smooth affine group schemes  on regular analytic surfaces $N^{^{(d)}}$ \eqref{nsgrpscheme1.5}.  These group schemes now get "parahoric structures" at the generic points of the rational chain. The { minimal platificateur} property mentioned above implies that the map $D^{^{(d)}} \to {\sf N}^{^{(d)}}$ is finite flat, and this is essential here. Each local type $\tau$ gives an affine {\sf 2}BT group scheme $\qh^{^G}_{_{\tau,{\sf N}^{^{(d)}}}}$  on the regular surface ${\sf N}^{^{(d)}}$ \eqref{nsgrpscheme2} and we arrive at the {\em standard list} of group schemes $\Big\{{\qh^{^G}_{_{\tau,{\sf N}^{^{(d)}}}}}\Big\}_{_{\tau}}$  on the regular surface ${\sf N}^{^{(d)}}$ (see \eqref{typeandlength} for the nomenclature). Bruhat-Tits theory has been studied extensively on discrete valuation rings, but there is, as of now, no "affine building" approach for higher dimensional regular local rings. Our objects give a large class of examples of such group schemes and from the philosophy of Bruhat-Tits, knowing the group schemes gives a hold on the possible "parahoric" subgroups.

This basic list then gives global group algebraic spaces on regular surfaces $S^{^{(d)}}$ which are proper over $A$ \eqref{properminmodel}; these group objects are obtained by a "gluing" \eqref{algspacelikeonsurfaces} the {\sf 2}BT group schemes to constant group schemes. These group algebraic spaces give degenerations of the constant group scheme $G$ on the generic fibre $C_{_K}$ to non-reductive limits on $C^{^{(d)}} \subset S^{^{(d)}}$, for varying $d$. The next step was to replace "torsion-free" sheaves on $(C,c)$ by torsors for the affine {\em non-reductive} group schemes on semistable curves $C^{^{(d)}}$ with fixed stable model being the curve $(C,c)$.  The torsors are also obtained by the Fourier-Mukai operation which were used to construct the group schemes, i.e. begin with equivariant torsors on the disc $D$ for the group scheme $\ce_{_D}(G)$ and realize them as torsors for the {\sf 2}BT group schemes $\qh^{^G}_{_{\tau,{\sf N}^{^{(d)}}}}$  on the regular surface ${\sf N}^{^{(d)}}$. The pairs $\Big\{\big[{\qh^{^G}_{_{\tau,{\sf N}^{^{(d)}}}}}, E\big]\Big\}_{_{\tau}}$  on the regular surface ${\sf N}^{^{(d)}}$ give the "admissible" list of objects. These can be globalized by gluing.

A basic off-shoot which emerges even in the case of $G = \text{GL}(n)$ is a certain "Tannakian" principle. {\em A priori} the admissible list of vector bundles on the rational chains, which allows only $\co$ and $\co(1)$ as summands, is clearly not closed under "tensor" operations on vector bundles. However, there is an underlying "parabolic structure" on these bundles when they are viewed as restrictions of bundles on the regular surface ${\sf N}^{^{(d)}}$. The rational chain is a normal crossing divisor (revealing a logarithmic structure) and the admissible data becomes a "parabolic data". This gives rise to a "parabolic tensor structure" which explains the phenomenon. This observation plays a central role in eventually constructing the stack.  

How does one do GIT for these objects? The approach is similar. Firstly, the local picture \eqref{keydiag1} at the level of surfaces, when restricted to curves, gives  global objects. More precisely, the disc $D$ can be replaced by a proper Deligne-Mumford stack called  twisted curves \eqref{twistedstuff} (\cite{ollson}), and the minimal resolution gets replaced by the normalization $\tC$ of $C$ \eqref{stackynormalization}.  One then sets up a Fourier-Mukai machinery between pairs consisting of [group scheme, torsor] on twisted curves  and pairs on $\tC$ which we call [balanced group schemes,{\tt laced} torsors]. In Nagaraj-Seshadri, semistability of Gieseker bundles on $C^{^{(d)}}$ had two ingredients, semistability of the torsion-free sheaf on $(C,c)$ (obtained by taking direct images via $C^{^{(d)}} \to (C,c)$) and a vertical component along the fibre of Gieseker vector bundles over a fixed torsion-free sheaf. We follows this approach. The torsion-free component is missing in the $G$-bundle setting and we replace it with $G$-torsors on twisted curves $\cc_{_d} \to C$. A heuristic "semistability" \eqref{newtfss} which abstractly captures the semistability of the "underlying" torsion-free object is then defined on the twisted curve $\cc_{_d}$. To make these heuristic objects concrete, we take the Fourier-Mukai path to get to the normalization $\tC$ of $C$. Here we can define numerical invariants such as {\em parabolic degrees} which allow us to concretize the slope semistability. Finally, ideas from GIT provide the precise {\tt L}-(semi)stability for Gieseker torsors, the ${\tt L}$ standing for an ample line bundle on a suitable "Quot-scheme"-like space which is an "atlas" for the stack.  To be precise, we define {\tt L}-(semi)stability for Gieseker torsors, which are pairs $(\qh^{^G}_{_{\tau,{C}^{^{(d)}}}}, \ce)$ consisting of [group scheme, torsor]  on semistable curves $C^{^{(d)}}$ together with a technical  "admissibility" condition. We restrict $(\qh^{^G}_{_{\tau,{C}^{^{(d)}}}}, \ce)$ to the normalization $\tC \subset C^{^{(d)}}$ which give rise to {\tt laced} torsors on $\tC$. The {\tt L}-(semi)stability of these {\tt laced} torsors is then used to define the notion for the pair $(\qh^{^G}_{_{\tau,{C}^{^{(d)}}}}, \ce)$.

\subsubsection{Related works} In the early nineties, in several papers, Bhosle introduced the notion of  "generalised parabolic bundles'' as a very useful tool to study the moduli space of torsion-free sheaves by working with objects on the normalization of the singular curve, but this had a intrinsic problem, that it was not amenable to the question of degeneration (however see \cite{bhosle1} and \cite{schmitt2}). Teixidor in several paper considered the moduli of bundles on singular curves (see \cite{teix}). T. Abe \cite{abe} solved the Gieseker construction for $\tsl(2)$, and Schmitt \cite{schmitt} constructed the universal Gieseker moduli over $\bar{\mathcal M}_{_g}$. M. Thaddeus had also considered the $\tsl(2)$ case in his thesis. The paper by Kiem and Li \cite{kiem-li} studied more explicit geometry of the Gieseker spaces towards applications. There is also a preprint by P. Solis \cite{solis} which should be of some interest.

In \cite{ns1}, Nagaraj and Seshadri had made some conjectures towards the problem for the case of $\text{\text{SL}(n)}$ in terms of the "determinant'' morphism on the moduli space of torsion-free sheaves. These conjectures were answered fully by Sun in \cite{sun0} and \cite{sun}. In 2000-2003, Friedman and Morgan  wrote several important papers (one with Witten)  (\cite{fm2}, \cite{fm3}, \cite{fmw}) on $G$-torsors on elliptic curves and singular curves. In 2004-2005, A. Schmitt, in a series of papers (\cite{schmitt},\cite{schmitt1}, \cite{schmitt2}), brought back the focus on the question of moduli space of $G$-torsors on singular curves and introduced some new ideas on "decorated bundles'' and their slope (semi)stability.

\ack {\sl I firstly thank my late teacher C.S. Seshadri for his faith in the entire work. His faith supported me in this long and arduous pursuit. I thank  B.Conrad, J.Martens, Johan de Jong, M.Thaddeus and R. Fringuelli for several helpful discussions and Miles Reid who remarked that my constructions are Fourier-Mukai-like in spirit. I thank D.S. Nagaraj, J. Heinloth, M. Brion, C. Simpson and Sourav Das for their comments and questions on an earlier version. I finally thank the referee for the conscientious and  meticulous reading of the manuscript and the numerous suggestions. These have gone a long way in improving the exposition}.

\subsubsection{Notations and Conventions}
Throughout this paper, unless otherwise stated, we have the
following notations and assumptions:
{\renewcommand{\labelenumi}{{\rm (\alph{enumi})}}
\begin{enumerate}

\item We work over an algebraically closed field $k$ of characteristic zero
and without loss of generality we can take $k$ to be the field of
complex numbers ${\bc}$.

\item {\em $(C,c)$ will be an irreducible projective nodal curve over $k$ with node $c \in C$ and $\nu:\tC \to C$ the normalization.}

\item Let $G$ be an {\em almost simple, simply
connected affine algebraic group} defined over $k$  of  $\text{rank}(G) = \ell = \text{ {\tt dim}}(T) $, where $T \,\subset\, G$ is a fixed maximal torus; let
$X(T) = \Hom(T,\, \bg_{_m})$ be the group of characters of $T$ and  $Y(T)\,=\, \Hom(\bg_m,\, T)$ be the group of all
 one--parameter subgroups of $T$. Fix a Borel subgroup $B$ containing $T$, and a set $\Delta$ of simple roots $\{\alpha_{_1}, \ldots, \alpha_{_\ell} \}.$
Let $\eR\,=\,R(T,\,G)$ denote the root system of $G$. Thus for every $r \,\in\, \eR$, there is the root homomorphism $u_r \,: \, \bg_a\,\to \, G$. The {\it
standard affine apartment} $\ca_{_T}$ is the affine space under $Y(T)
\otimes_\bz \br$. and {\it we shall identify $\ca_{_T}$ with $Y(T) \otimes_\bz \br$}
(see \cite[\S~2]{base}).  
\item All group schemes considered in this paper are {\sf affine}.

\item Let $A := \spec k\llbracket t \rrbracket$ and  $K := \spec k(\!(t)\!)$ and $o \in A$ the closed point. Let $\pi:C_{_{A}} \to {A}$ be such that $\pi_{_K}:C_{_K} \to K$ is a smooth projective curve of genus $g \geq 2$ and $C_{_o} \simeq (C,c)$. {\em We assume that $C_{_{A}}$ is regular over $k$}. Let $U_{_{c}} \subset C_{_{A}}$ be an analytic neighbourhood of the node $c \in C_{_{A}}$.

\item\label{somenotas} Let $d > 0$ be a positive integer and let ${\sf\mu_{_d}} = \langle \gamma \rangle$ be the cyclic group of order $d $. The group ${\sf\mu_{_d}}$ is considered as a subgroup of $\text{SL}(2, k)$ generated by $g = \left(\begin{array}{cc} \zeta & 0 \\0 & \zeta^{^{-1}} \\ \end{array}\right)$, where $\zeta = e^{^{2i\pi/d}}$ is a primitive $d^{^{th}}$-root of unity. 
\item If $\rho:\smu \to G$ is a representation, $\tau$ will stand for its type \eqref{4'} and represent the conjugacy class of $\rho$.

\item The $S^{^{(d)}}$ \eqref{properminmodel}, are smooth surfaces with a projective morphism to $A$. These are minimal desingularizations of the surfaces with normal singularities with  local equation $x.y = t^{^d}$ obtained by base change from $C_{_A}$. The exceptional fibre $C^{^{(d)}}$, is the semistable curve with $d-1$-chain of rational curves glued to the normalization $\tC$. ${\sf N}^{^{(d)}}$ is the local analytic neighbourhood of the exceptional divisor in $S^{^{(d)}}$.
\item $W[d]$ \eqref{junlistdmodel} are Jun Li's standard models and $Z[d]$ the local standard models \eqref{thekeyconstruct}.
\item The map ${\nu}:(\tilde{C},{\bf c}) \to (C,c)$ denotes the normalization of $C$ where ${\nu}^{-1}(c) = \bf c$ and $\bf c$ stands for the pair of points $\{c_{_1},c_{_2}\}$.
\item $\cc_{_d}$ is a twisted curve \eqref{twistedstuff}  in the sense of \cite{ollson}. 
\item $E_{_\wp}$ are {\tt laced} torsors on the normalization $\tC$ \eqref{bttorsorlacedtorsor} \eqref{2-lacedbundle}.
\item $\qh^{^G}_{_{\tau,{\sf N}^{^{(d)}}}}$ are the {\tt 2}BT-group schemes on  the regular surface ${\sf N}^{^{(d)}}$ \eqref{nsgrpscheme2}.

\end{enumerate}

\begin{center}
\underline{\sc Part I}
\end{center}
\vspace{-3mm}

\section{Preliminaries}

We recall the obvious identification (\cite[2.2.8]{base})
\beqa\label{initialmumbo} 
Hom({\sf\mu_{_d}},T) \simeq {\frac{Y(T)}{d. Y(T)}}.\eeqa

Let $\rho: {\sf\mu_{_d}} \to G$ be a representation. Since ${\sf\mu_{_d}}$ is cyclic, we can suppose
that the representation $\rho$ of ${\sf\mu_{_d}}$ in $G$ factors through 
$T$ (by a suitable conjugation). The cocharacter $\alpha^{^*}_{_\rho}$ associated to $\rho$ by (\ref{initialmumbo}) gives a tuple of integers $\{a_{_1}, a_{_2}, \ldots, a_{_\ell} \}$ determined uniquely modulo $d$ and in terms of the canonical cocharacters $\{\alpha^{^*}_{_j}~\in~Y(T), j = 1, \ldots, \ell\} $ dual to the simple roots $\alpha_{_j}$ we have: 
\beqa\label{4'}
\alpha^{^*}_{_\rho} = \sum_{_{j = 1}}^{^\ell} {a_{_j}} \alpha^{^*}_{_j}.
\eeqa
We will call the tuple $\tau:= (a_{_1}, a_{_2}, \ldots, a_{_\ell})$ the  {\em type} of the representation $\rho$  and denote the association in (\ref{initialmumbo}) by:
\beqa\label{localtypeandweights}
\rho \mapsto \theta_{_{\tau}},
\eeqa
where $\theta_{_\tau} \in Y(T)/d.Y(T)$. We view $\theta_{_\tau}$ as a point in the affine apartment $\ca_{_T}$.

\subsubsection{The geometric setting and assumptions}

\bnot\label{normalization} Let ${\nu}:(\tilde{C},{\bf c}) \to (C,c)$ denote the normalization of $C$  and let ${\nu}^{-1}(c) = \{c_{_1},c_{_2}\}$. Let $E(1) \cup E(2)$ be the normalization of the analytic neighbourhood of $c \in C$. \enot
\bdefe\label{Chains} A scheme $E^{^{(m)}}$ is called a chain of rational curves if $$E^{^{(m)}} = \bigcup_{_{i=1}}^{^{m-1}} E_{_i},$$ with $E_{_i} \simeq \bp^1$, and if $i \neq j$, 
\beqa
E_{_i} \cap E_{_j}  =  \begin{cases}
\text{singleton}&   \text {if $\mid i-j\mid = 1$}\\
\emptyset& \text {otherwise}
\end{cases}
\eeqa
\edefe
\bdefe\label{Gieseker curve} Let $C^{^{(d)}}$ denote the reducible nodal curve with components being the normalization $\tilde{C}$ of $C$ and a chain $E^{^{(d)}}$ of projective lines  of length $d-1$ attached to $\tilde{C}$ at $c_{_1}$ and $c_{_2}$. Equivalently, it is a  semistable curve which has $C$ as its stable model.  If $p:C^{^{(d)}} \to C$ denotes the canonical morphism, the inverse image $p^{-1}(c)$ is the chain $E^{^{(d)}}$. \edefe
We have the diagram:
\beqa\label{normalplusmodif}
\xymatrix{
(\tC, {\bf c}) \ar[dr]_{{\nu}} \ar@{^{(}->}[rr]& &(C^{^{(d)}}, E^{^{(d)}})\ar[dl]^{p} \\
& (C,c)  &
}
\eeqa
Let $C_{_{d,A}} := C_{_A} \times_{_{\spec k\llbracket t \rrbracket}} \spec k\llbracket t \rrbracket$ via the map $t \mapsto t^{^d}$. Here and elsewhere $A = \spec~k\llbracket t \rrbracket$ and $o \in A$ the closed point.
Let $N_{_d} = \spec \frac{A\llbracket x,y \rrbracket}{(x.y - t^{^{d}})}$ be the analytic neighbourhood of $c$ in $C_{_{d,A}}$ which lies above the analytic neighbourhood $U_{_c} \subset C_{_{A}}$. We recall (\cite[Page 191]{ns2}) that $N_{_d}$ is a {\em normal} surface with an isolated singularity at $c$ of type $\text{A}_{_d}$. By the generality of $\text{A}_{_d}$-type singularities, one can realize $N_{_d}$ as a quotient $\sigma:D \to N_{_d}$ of  $D:= \spec~\frac{A\llbracket u,v \rrbracket}{(u.v - t)}$ by the cyclic group ${\sf\mu_{_d}}$, where ${\sf\mu_{_d}}$ acts on $D$ as 
follows:
\beqa\label{gammaaction}
\gamma.(u,v) = (\zeta.u, \zeta^{^{d-1}}.v)
\eeqa
and $x = u^{^{d}}, y = v^{^{d}}$.
We consider the following basic diagram for all $d > 0$ (see \cite{gonverd}):
\beqa\label{keydiag1}
\xymatrix{
{D^{^{(d)}}} \ar[r]^{f} \ar[d]_{q} &
 {\sf N}^{^{(d)}}  \ar[d]_{p_{_d}} \\
0 \in D \ar[r]^{\sigma} &  ~~{N_{_d}} \owns c \\
}
\eeqa
where  $p_{_d}:{\sf N}^{^{(d)}} \to N_{_d}$ is the {\em minimal resolution of singularities} of $N_{_d}$ obtained by successively blowing up the singularity, with the exceptional divisor $E^{^{(d)}} = p_{_d}^{^{-1}}(c)$, and 
\beqa
D^{^{(d)}} := \huge(D \times _{_{N_{_d}}} {\sf N}^{^{(d)}} \huge)_{_{red}}.
\eeqa
The {\em closed fibre} $F^{^{(d)}} (= {\sf N}^{^{(d)}}_{_o})$ of the canonical morphism ${\sf N}^{^{(d)}} \to A$ looks like:
\beqa\label{closedfibre}
F^{^{(d)}} = E^{^{(d)}} \cup E(1) \cup E(2)
\eeqa
Thus, $F^{^{(d)}} \subset {\sf N}^{^{(d)}}$ is a normal crossing divisor with $d+1$ components. 
\begin{figure}\label{chain}
\begin{center}
\includegraphics[scale=0.70]{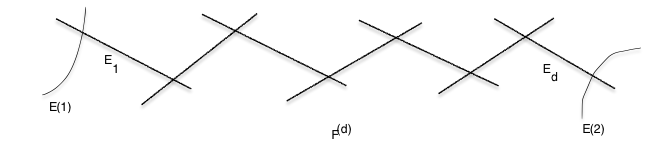}
\caption{{The fibre $F^{^{(d)}}$}}
\label{chain}
\end{center}
\end{figure} 
By \cite[Proposition 2.4]{gonverd}  the morphism  
\beqa\label{platif}
f:D^{^{(d)}} \to {\sf N}^{^{(d)}}
\eeqa
is {\em finite and flat, the minimal platificateur in the sense of Grothendieck} \cite[Cor.7, page 448]{gonverd}}. Since ${\sf N}^{^{(d)}}$ is smooth, this implies that $f$ is ramified at the generic point of each of the $d-1$ rational components of the exceptional divisor $E^{^{(d)}} = p_{_d}^{^{-1}}(c) \subset F^{^{(d)}}$. Let 
\beqa\label{properminmodel}
p_{_d}:S^{^{(d)}} \to C_{_{d,A}}
\eeqa
be {\em the} minimal smooth model for $C_{_K}$. Then we see that $ {\sf N}^{^{(d)}} \to S^{^{(d)}}$ gives an analytic neighbourhood of the exceptional fibre $p_{_d}^{^{-1}}(c)$.
\begin{figure}\label{minres}
\begin{center}
\includegraphics[scale=0.70]{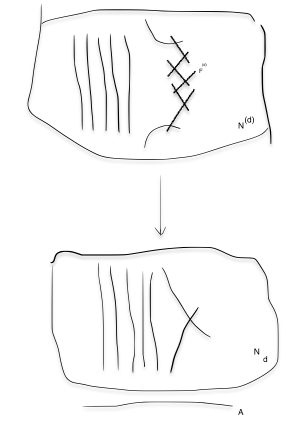}
\caption{{Minimal resolution of $N_{_d}$}}
\label{minres}
\end{center}
\end{figure} 
\brem\label{The \'etale picture} (The \'etale picture) The surface $C_{_{A}}$ over $A = \spec~k\llbracket t \rrbracket$ is assumed to be regular and hence the analytic local ring at the node $c \in C$ is $A\llbracket x,y \rrbracket/(xy -t)$. It is well known that $\spec~A\llbracket x,y \rrbracket/(xy -t)$ is the analytic local ring for a versal deformation of the simple node. By  (\cite[Proposition 2.8, page 184]{frietag}) which is somewhat delicate, one can in fact obtain an \'etale  neighbourhood $U(c)$ of $c$ in $C_{_{A}}$ which is isomorphic to an \'etale neighbourhood of the origin $0$ in $\spec~A[x,y]/(xy -t)$. By a base change by the map $k\llbracket t \rrbracket \to k\llbracket t \rrbracket$ given by $t \mapsto t^{^d}$, we see that there is an \'etale neighbourhood $N_{_{_{_{\'et,d}}}}$ of $c$ in $C_{_{d,A}}$ which is isomorphic to an \'etale neighbourhood of the origin $0$ in $\spec~A[x,y]/(xy-t^{^d})$. In other words, in the \'etale topology, we can express the neighbourhood $N_{_{_{_{\'et,d}}}}$ of $c$ in $C_{_{d,A}}$ as a quotient $\ba^{^{2}}/\smu$ for the affine space $D_{_{\'et}} := \ba^{^2} = \spec~\bc[u,v]$ for the action \eqref{gammaaction}. This gives the following \'etale picture corresponding to \eqref{keydiag1}:
\beqa\label{keydiagetale}
\xymatrix{
{D_{_{\'et}}^{^{(d)}}} \ar[r]^{f_{_{\'et}}} \ar[d]_{q} &
 {\sf N}_{_{\'et}}^{^{(d)}} \ar[d]^{p_{_d}} \\
D_{_{\'et}} \ar[r]^{\sigma_{_{\'et}}} &  N_{_{_{_{\'et,d}}}}
}
\eeqa
where $ {\sf N}_{_{\'et}}^{^{(d)}}$ is the minimal desingularization of the normal surface $N_{_{_{_{\'et,d}}}}$ and   
\beqa
{D_{_{\'et}}^{^{(d)}}} := \huge(D_{_{\'et}} \times _{_{N_{_{_{_{\'et,d}}}}}} {\sf N}_{_{\'et}}^{^{(d)}} \huge)_{_{red}}.
\eeqa
Observe that $f_{_{\'et}}:D_{_{\'et}}^{^{(d)}} \to {\sf N}_{_{\'et}}^{^{(d)}}$ is finite and flat since the map $f$ in \eqref{keydiag1} is so. If $p_{_d}$ is as in \eqref{properminmodel},  then we see that ${\sf N}_{_{\'et}}^{^{(d)}} \to S^{^{(d)}}$ gives an \'etale neighbourhood of the exceptional fibre $p_{_d}^{^{-1}}(c)$.\erem

\brem\label{balancedaction} ({\em Balanced action}) The action of ${\sf\mu_{_d}}$  is {\em balanced} in the sense of \cite{a-v} (see also \cite[2.5]{twisted}), i.e., the action of a generator $\zeta$ on the tangent spaces to each branch are inverses to each other.  For the corresponding dual action in the neighbourhood $D'_{_0}$ (the component with local coordinate $v$),  the action is by $\zeta^{^{-1}}$ (see \eqref{chinought}). If we begin with a representation $\rho:\smu \to G$ of local type $\tau$ at a point in a branch, then corresponding local type for the dual action at the point in the second branch is denoted by $\bar{\tau}$. See \eqref{balancedparabolicbundle} for the expression of "dual weights" when $G = \text{GL}(n)$.
 \erem
\subsubsection{Outline of proof strategy}
We give the broad steps of the proof.\\
\noindent
\fbox{%
\parbox{\textwidth}{

1.We work with $F^{^{(d)}}$ \eqref{closedfibre} inside the analytic surface ${\sf N}^{^{(d)}}$, i.e. the basic models are built on smooth analytic surfaces. }%
}
\noindent
\fbox{%
\parbox{\textwidth}{

 2. The geometric McKay correspondence is then used to define the local models for the group schemes in \S3.}%
}

\noindent
\noindent
\fbox{%
\parbox{\textwidth}{
 3. Unlike the vector bundle case, admissibility will be defined on the analytic surface ${\sf N}^{^{(d)}}$ and then extended to more general "modifications" defined in \eqref{modification}.}%
}
\fbox{%
\parbox{\textwidth}{
\noindent
 4. These group schemes on local analytic models get realized as invariant push-forwards from Kawamata covers \eqref{schemelikeonsurfaces}. The logarithmic structure on the surfaces becomes significant.}%
}
\fbox{%
\parbox{\textwidth}{
\noindent
 5. Use the local analytic model to define group schemes on Jun Li's local standard models $Z[d]$ \eqref{thekeyconstruct}. Realize these also from Kawamata covers of $Z[d]$ \eqref{sncinzd}, \eqref{schemelikeonsurfaces1}.}%
}
\fbox{%
\parbox{\textwidth}{
\noindent
 6. Globalize and define group algebraic spaces on Jun Li's standard models $W[d]$ via Kawamata covers of $W[d]$ \eqref{kawamatafortheversal}.}%
}
\fbox{%
\parbox{\textwidth}{
\noindent
 7. Local models for torsors on ${\sf N}^{^{(d)}}$ are made by pulling back equivariant torsors and then taking invariant push-forwards \eqref{admissibletorsor}, \eqref{admingies}.}%
}
\fbox{%
\parbox{\textwidth}{
\noindent
 8. We globalize and define group algebraic spaces on smooth  surfaces $S^{^{(d)}}$ (projective over $A$) by gluing and invariant push-forwards. Kawamata covers of $S^{^{(d)}}$ can be defined using the ramification data on $N^{^{(d)}}$.}%
}
\fbox{%
\parbox{\textwidth}{
\noindent
 9. Torsors for these group algebraic spaces on  $S^{^{(d)}}$ are obtained by invariant push-forwards of equivariant torsors on these Kawamata covers.}%
}
\fbox{%
\parbox{\textwidth}{
\noindent
 10. Admissible pairs on $S^{^{(d)}}$ \eqref{admissibletorsor} are torsors obtained by invariant push forwards of equivariant torsors using specific local data.}%
}
\fbox{%
\parbox{\textwidth}{
\noindent
 11. Define {\em parabolically associated} \eqref{admonstdmodels} vector bundles to admissible pairs. This process gets done by taking associated equivariant vector bundles on covers and then taking invariant push forwards. These parabolically associated vector bundles are the quasi-admissible vector bundles on the standard models $W[d]$. }%
}
\fbox{%
\parbox{\textwidth}{
\noindent
 12. Define admissible pairs on arbitrary modifications using Jun Li's expanded degenerations and local effectivity \eqref{junli1} and \eqref{admississame} and show that this definition is intrinsic.}%
}

\normalsize

\section{McKay correspondence and {\tt 2}Bruhat-Tits group schemes}\label{mckaybtetc} 
The aim of this section is to construct  certain smooth affine group schemes  on the ${\sf N}^{^{(d)}}$ which are  {\em generically split}, i.e. a product over the open subset ${\sf N}^{^{(d)}} - E^{^{(d)}}$ \eqref{keydiag1} with fibre $G$ and which degenerate {\em parahorically}. More precisely, these group schemes are $2$-dimensional generalizations of the classical Bruhat-Tits group schemes associated to parahoric subgroups of $G(K)$ (see \eqref{typeandlength}).  I make these constructions using the geometric McKay correspondence  of Gonzalez-Springberg and Verdier (\cite{gonverd}).

Let $\eT_{_{D}}$ be  a $({\sf\mu_{_d}},G)$-torsor  on $D$ \eqref{keydiag1} (see \eqref{basictorsor} below). Assume that it is given by a homomorphism $\rho:{\sf\mu_{_d}} \to G$ (in fact, it is easy to see that this is always the case on $D$).  This gives a homomorphism $\rho:{\sf\mu_{_d}} \to T$ into the maximal torus $T$ of $G$ of  {\em type} $\tau = (a_1, a_2, \ldots, a_{\ell})$ in the sense of \eqref{4'}. In other words, we have a ${\sf\mu_{_d}}$-action on $D \times G$,  given by 
\beqa\label{localaction0}
\gamma.(u,v,g) = (\zeta.u, \zeta^{^{-1}}.v, \rho(\gamma).g), 
\eeqa
and: 
\beqa\label{basictorsor}
\eT_{_{D}} \simeq D \times^{^{\rho}} G.
\eeqa
We observe that since the action of ${\sf\mu_{_d}}$ is {\em balanced} (\ref{balancedaction}) at the two marked points $z_{_1}$ (resp. $z_{_2}$) above the origin $0 \in D$, the local type of the action on a $G$-torsor at these points are $\tau$ (resp. $\bar\tau$). {\em Throughout the paper, we fix this $G$-torsor $\eT_{_{D}}$ on $D$ of local type $\tau$}.

Consider the  adjoint group scheme $\eT_{_{D}}(G):= \eT_{_{D}} \times ^{^{G,Ad}} G$ on $D$, where $G$ acts on itself by inner conjugation.  We define the equivariant group scheme:
\beqa\label{egtau}
E(G, \tau):= q^{*}(\eT_{_{D}}(G))
\eeqa
on $D^{^{(d)}}$ of local type $\tau$ in the sense that it comes with a ${\sf\mu_{_d}}$-action via a representation $\rho:{\sf\mu_{_d}} \to G$. 
Since the morphism $f:D^{^{(d)}} \to {\sf N}^{^{(d)}}$ is also {\em finite and flat}  we can take the Weil restriction of scalars:
\beqa\label{localnsgrpscheme1}
f_{_*}(E(G, \tau)_{_{D^{^{(d)}}}}):= {\mathscr R}_{_{{{D^{^{(d)}}}}\big/{\sf N}^{^{(d)}}}} (E(G, \tau)_{_{D^{^{(d)}}}}).
\eeqa
and since $E(G, \tau)_{_{D^{^{(d)}}}} \to D^{^{(d)}}$ is a smooth (affine) group scheme,  the basic properties of Weil restriction of scalars (\cite[Lemma 2.2]{edixhoven}) show that $f_{_*}(E(G, \tau)$ is a smooth group scheme on ${\sf N}^{^{(d)}}$, coming with a ${\sf\mu_{_d}}$-action. By taking  invariants under the action of ${\sf\mu_{_d}}$ and noting that we are over characteristic zero, by \cite[Prop 3.4]{edixhoven}, we obtain the smooth (affine) group scheme obtained by what we shall term {\bf invariant direct image}:
\beqa\label{nsgrpscheme1.5}
\qh^{^G}_{_{\tau,{\sf N}^{^{(d)}}}} = ({\tt Inv}\circ f_{_*}) \Big((E(G, \tau)_{_{D^{^{(d)}}}})\Big) :=  \Big[{{\mathscr R}_{_{{{D^{^{(d)}}}}\big/{\sf N}^{^{(d)}}}} (E(G, \tau)_{_{D^{^{(d)}}}}}\Big]^{^{{\sf\mu_{_d}}}}.
\eeqa 
on ${\sf N}^{^{(d)}}$ (see \cite[Definition 4.1.3]{base}). 
\bdefe\label{nsgrpscheme2} Let $G$ be an almost simple, simply connected algebraic group of rank $\ell$ and let $\theta_{_\tau} \in \ca_{_T}$ be a weight in the affine apartment \eqref{localtypeandweights} arising from a representation $\rho$, and $d$ a positive integer. The  {\tt 2BT}-group scheme of type $\tau$ with generic fibre $G$ of singularity type $\text{A}_{_d}$  associated to $\theta_{_\tau}$ is defined to be the affine group scheme $\qh^{^G}_{_{\tau,{\sf N}^{^{(d)}}}}$ \eqref{nsgrpscheme1.5} on the regular surface ${\sf N}^{^{(d)}}$ (see \eqref{typeandlength} for the nomenclature). This process defines  a distinguished collection of {\tt 2BT}-group scheme $\Big\{{\qh^{^G}_{_{\tau,{\sf N}^{^{(d)}}}}}\Big\}_{_{\tau}}$ indexed by the type $\tau$. \edefe
As an instance  of what we would be doing subsequently over more general modifications \eqref{modification},  we show that  we can obtain such a group algebraic space by a different more general geometric construction (this is possible since we are over characteristic $0$). This is via the Kawamata covering lemma \eqref{kawa}. Note that $D^{^{(d)}}$ in \eqref{keydiag1} is only Cohen-Macaulay in general.
\bth\label{schemelikeonsurfaces} The group scheme  $\qh^{^G}_{_{\tau,{\sf N}^{^{(d)}}}}$ can be realized as an invariant direct image \eqref{nsgrpscheme1.5} of a group scheme with fibres $G$, from a global smooth ramified covering $\tilde{\sf N}(\tau)$ of  the smooth surface ${\sf N}^{^{(d)}}$.
\eeth
\bpr  The closed fibre $F^{^{(d)}} \subset {\sf N}^{^{(d)}}$ \eqref{closedfibre} is a reduced divisor with normal crossing singularities. Further, the covering $f:D^{^{(d)}} \to  {\sf N}^{^{(d)}}$ \eqref{keydiag1} of the analytic neighbourhood of the closed fibre is ramified at the generic points of the chain of rational curves $E^{^{(d)}} \subset F^{^{(d)}}$ with a choice of ramification indices. We fix this ramification data. 

Recall the Kawamata covering lemma \eqref{kawa} (see \cite[Theorem 17]{kawamata},\cite[Lemma 2.5, page 56]{vieweg}). We have a Galois covering $h: {\sf V}_{1} \to {\sf N}^{^{(d)}}$ with ${\sf V}_{1}$ smooth, which is ramified along the irreducible components $E^{^{(d)}}$ for the fixed ramification data dictated by the local type $\tau$ of the representation. Let $\Gamma = \text{Gal} (h)$ be the Galois group of the covering $h$. The stabilizer subgroups of $\Gamma$ at each of the generic points of the rational components $R_{_j}$'s of $E^{^{(d)}}$ is precisely the cyclic group ${\sf\mu_{_d}}$. 

We now consider the two coverings $h:{\sf V}_{1} \to {\sf N}^{^{(d)}}$ and $f:D^{^{(d)}} \to  {\sf N}^{^{(d)}}$. The covering $f$ is \'etale over the complement of $E^{^{(d)}}$. Hence by \cite[Corollary 2.6]{vieweg}, we have a smooth finite covering $\tilde{\sf N}(\tau) \to D^{^{(d)}}$ such that there is an \'etale morphism $j: \tilde{\sf N}(\tau) \to {\sf V}_{1}$ which can be assumed to be Galois (by going to the canonical Galois closure if need be). The group scheme $E(\tau,G)$ on $D^{^{(d)}}$ pulled back by $\tilde{\sf N}(\tau) \to D^{^{(d)}}$ gives a group scheme $E(\tau,G)'$ on $\tilde{\sf N}(\tau)$ and hence on ${\sf V}_{1}$. Note that the composite  $h \circ j$ gives a Galois cover.

Thus we can again take the Weil restriction of scalars and invariants under the composite $h \circ j = \varphi:\tilde{\sf N}(\tau) \to {\sf N}^{^{(d)}}$ and we get a group scheme $({\tt Inv} \circ \varphi_{_*})\Big(E(\tau,G)'\Big)$. It is now easily checked that this group scheme is isomorphic to  $\qh^{^G}_{_{\tau,{\sf N}^{^{(d)}}}}$.
\epr

\subsubsection{Global constructions}\label{fixsomething}

Fix a $G$-torsor $E_{_A}^{^o}$ on $C_{_{A}} - c$ and let $E^{^o}_{_{d,A}}$ be the pull-back to $C_{_{d,A}} - c$. Let $\eT_{_{D_{_{\'et}}}}$ be the $(\smu,G)$-torsor  on $D_{_{\'et}} = \ba^{^{2}}$  be  as in \eqref{basictorsor} given by $\rho:\smu \to G$ of  {\em type} $\tau$ i.e. we have a $\smu$-equivariant trivialization $\eT_{_{D_{_{\'et}}}} \simeq D_{_{\'et}} \times^{^{\rho}} G$. We can again consider the adjoint group scheme $\eT_{_{D_{_{\'et}}}}(G)$ which is an equivariant group scheme on $D_{_{\'et}}$. Note that we can view $\eT_{_{D_{_{\'et}}}}$ as an $\eT_{_{D_{_{\'et}}}}(G)$-torsor as well. As in the analytic case, we see that we have a non-trivial  group scheme $E(G, \tau)_{_{D_{_{\'et}}^{^{(d)}} }}$ on $D_{_{\'et}}^{^{(d)}}$ and we can take the invariant direct image to get an affine group scheme $\qh^{^G}_{_{\tau,{\sf N_{_{\'et}}^{^{(d)}}}}}$ on $\sf N_{_{\'et}}^{^{(d)}}$. The same process also gives a $\qh^{^G}_{_{\tau,{\sf N_{_{\'et}}^{^{(d)}}}}}$-torsor $E_{_{\sf N_{_{\'et}}^{^{(d)}}}} := ({\tt Inv}\circ f_{_*})(q^{^*}(\eT_{_{D_{_{\'et}}}})$

Fixing a choice of $E_{_A}^{^o}$ as above and consider the adjoint group scheme $p_{_d}^{^*}(E_{_A}^{^o}(G))$ on ${S}^{^{(d)}} - F^{^{(d)}}$, $F^{^{(d)}}$ as in \eqref{closedfibre}. We now define $\qh^{^G}_{_{\tau,S^{^{(d)}}}}$,
as {\em a group algebraic space of ``local type" $\tau$} (in fact by a group scheme in the \'etale topology on $S^{^{(d)}}$) on the minimal smooth model $S^{^{(d)}}$, by {\em gluing} the group scheme $\qh^{^G}_{_{\tau,{\sf N_{_{\'et}}^{^{(d)}}}}}$ with $p_{_d}^{^*}(E_{_A}^{^o}(G))$ on ${S}^{^{(d)}} - F^{^{(d)}}$. We make a choice for this gluing. But having done that, we also observe that this choice does not affect the study since our main concern is with torsors under these group schemes.

The sheaf in the \'etale topology on $S^{^{(d)}}$ is represented by a group algebraic space (\cite[Theorem 1.5, page 165]{milne}). By faithfully flat descent (\cite[Section 6.5, Example D]{blr}), we can deduce that a  group algebraic space $\qh^{^G}_{_{\tau,S^{^{(d)}}}}$ on $S^{^{(d)}}$ is scheme-like over all height one prime ideals and in particular at the generic points of the components of  $C^{^{(d)}}$. Moreover, the restriction $\qh^{^G}_{_{\tau, C^{^{(d)}} }} := \qh^{^G}_{_{\tau,S^{^{(d)}}}} \bigg\vert_{_{C^{^{(d)}}}}$ of $\qh^{^G}_{_{\tau,S^{^{(d)}}}}$ to $C^{^{(d)}} \subset S^{^{(d)}}$ is a veritable group scheme which is also immediate since we are gluing group schemes on smooth curves  and (\cite[Section 6.5, Example D]{blr}) applies. Indeed, $\qh^{^G}_{_{\tau, C^{^{(d)}} }}$ is a group scheme obtained by gluing the closed fibre  $\qh^{^G}_{_{\tau, F^{^{(d)}} }}$ of $\qh^{^G}_{_{\tau,{\sf N}^{^{(d)}}}}$  with a semisimple group scheme on $C^{^{(d))}} - F^{^{(d)}}$.

We begin by observing that as in the analytic case we have a Kawamata cover ${\sf \tilde{N}_{_{\'et}}}(\tau) \to {\sf N_{_{\'et}}^{^{(d)}}}$ with the same ramification data as before.

\bcor\label{algspacelikeonsurfaces} The group algebraic space $\qh^{^G}_{_{\tau,S^{^{(d)}}}}$ of local type $\tau$ on $S^{^{(d)}}$ can be realized as an invariant direct image \eqref{nsgrpscheme1.5} of a group algebraic space $\mathcal G$ with constant fibres $G$, from a global smooth ramified covering $\tilde{S}(\tau)$ of  the smooth surface $S^{^{(d)}}$.\ecor
\bpr The closed fibre $C^{^{(d)}} \subset S^{^{(d)}}$ is a reduced divisor with normal crossing singularities. We can therefore get a smooth ramified covering $\tilde{S}(\tau) \to S^{^{(d)}}$ which is ramified over the precise locus in $C^{^{(d)}}$ (depending on $\tau$) which local analytically gives the Kawamata cover $\tilde{\sf N}(\tau) \to {\sf N}^{^{(d)}}$. In fact, this can be done in the \'etale setting and we get an equivariant group algebraic space over $\tilde{S}(\tau)$ whose invariant direct image is $\qh^{^G}_{_{\tau,S^{^{(d)}}}}$. We recall that Weil restriction of scalars sends algebraic spaces to algebraic spaces.
\epr

\subsubsection{The McKay correspondence revisited}\label{mckaystory} Before going to the salient feature of the group schemes $\qh^{^G}_{_{\tau, {\sf N}^{^{(d)}}}}$, we recall the geometric interpretation of the McKay correspondence given by Gonzalez-Springberg and Verdier (\cite{gonverd}).  Let $\text{Irr}^{^o}({\sf\mu_{_d}}) \subset \text{Irr}({\sf\mu_{_d}})$ be the nontrivial irreducible representations of ${\sf\mu_{_d}}$ and let $\text{Irr}(E^{^{(d)}})$ denote the set of irreducible rational components of  the exceptional divisor  of the minimal resolution $p_{_d}:{\sf N}^{^{(d)}} \to N_{_d}$ \eqref{keydiag1}.  Let $\psi$  be a non-trivial character of ${\sf\mu_{_d}} = \langle \gamma \rangle$. Then $\psi$ corresponds to $\zeta \mapsto \zeta^{^s}$, where $\zeta$ corresponds to a primitive $d^{^{th}}$-root of $1$ and $1 \leq s \leq d-1$. Let $L_{_{\psi}}$ be the equivariant line bundle  on $D$ where ${\sf\mu_{_d}}$ acts on $D \times k$ as $\gamma.(u,v, a) = (\zeta.u, \zeta^{{d -1}}.v, \zeta^{^s}.a),~~~ a \in k.$
A ${\sf\mu_{_d}}$-invariant section \text{\cursive s} of this line bundle is given by the relation $\text{\cursive s} (\gamma.(u,v))  = \zeta^{^s} \text{\cursive s}(u,v)$, and hence the ${\sf\mu_{_d}}$-invariant sections are generated by $u^{^s}$ and $v^{^{d-s}}$. From this, it is easily checked that the invariant direct image of $L_{_{\psi}}$ under $\sigma:D \to N_{_d}$ is given by:
\beqa\label{afterthis}
({\tt Inv} \circ \sigma_{_*})(L_{_{\psi}})  = (x, t^{^{d-s}})
\eeqa
where $(x, t^{^{d-s}})$ is an ideal sheaf on $N_{_d} = \spec~k\llbracket x,y,t \rrbracket/(x.y - t^{^{d}})$ \eqref{gammaaction}. 
Let  $\cl_{_\psi} := {\tt Inv} \circ f_{_*}(q^{^{*}}(L_{_\psi}))$ be the induced line bundle on $\sf N^{^{(d)}}$. This is a line bundle since $f$ is finite and flat.
\bth\label{mckay} (Gonzalez-Springberg, Verdier) There is a bijection $
\text{Irr}^{^o}({\sf\mu_{_d}}) \to \text{Irr}(E^{^{(d)}})$, $\psi \mapsto E_{_\psi}$,  such that for any $E_{_j} \in \text{Irr}(E^{^{(d)}})$, we have
\beqa
c_{_1}(\cl_{_\psi}).E_{_j} =   \left \{
\begin{array}{rcl}
0
~~{\rm if}~ E_{_j} \neq E_{_\psi} \\
~1
~{\rm if}~E_{_j} = E_{_\psi}.\\
\end{array}\right. \eeqa
\eeth
The statement in (\ref{mckay}) implies that the first Chern class $c_{_1}(\cl_{_\psi})$ can be represented by a divisor which meets $F^{^{(d)}}$ transversally at a unique point   which lies in $E_{_\psi}$. More precisely, consider the divisor $q^{^{-1}}(0) \subset D^{^{(d)}}$, with  $0 \in D$ as in \eqref{keydiag1}.  Consider the reduced fibre $\tilde{E} :=  q^{^{-1}}(0)_{_{red}}$. The group ${\sf\mu_{_d}}$-fixes the divisor $q^{^{-1}}(0)$ and hence its reduced subscheme $\tilde{E}$. Thus,  given $\psi \in \text{Irr}^{^o}({\sf\mu_{_d}})$, there is a unique component $\tilde{E_{_\psi}}$ of $\tilde{E}$  such that the line bundle $q^{^{*}}(\cl_{_{\psi}})$ gets a nontrivial linearization by the ${\sf\mu_{_d}}$ action at the generic point of $\tilde{E_{_\psi}}$.

\subsubsection{A brief description of the group scheme}\label{briefdescript} Suppose that we are given a  homomorphism $\rho:{\sf\mu_{_d}} \to T$ of {\em local type} $\tau$.  For the simple roots $\{\alpha_{_j}\}_{_{j = 1}}^{^\ell}$ of $G$, let 
\beqa\label{chainrule}
{\sf C}_{_\rho} := \{\alpha_{_j} \circ \rho \mid \alpha_{_j} \circ \rho \neq 1.\}\eeqa
Then  ${\sf C}_{_\rho}$ gives a subset of $\text{Irr}^{^o}({\sf\mu_{_d}})$. The McKay correspondence says that to each $\psi \in {\sf C}_{_\rho}$ we have a unique rational component $E_{_\psi} \subset F^{^{(d)}}$. Furthermore, the covering $f:D^{^{(d)}} \to {\sf N}^{^{(d)}}$ is ramified precisely over the rational curves $E_{_\rho}:= \{E_{_\psi} \mid \psi \in {\sf C}_{_\rho} \}$ with ramification index dictated by the number of $\psi$'s which give independent characters of ${\sf\mu_{_d}}$ and their multiplicities. 
\brem\label{typeandlength} By Bruhat-Tits theory, for each facet $\ul{a}$ in the apartment $\ca_{_T}$ of the Bruhat-Tits building of $G(k(\!(t)\!))$, there is a smooth group scheme $\sf {P}_{_{\ul{a}}}$ over $\spec~k\llbracket t \rrbracket$ with connected fibers whose generic fiber is $G \times _{_{\spec k}} \spec k(\!(t)\!)$. We call such  a $\sf{P}_{_{\ul{a}}}$ a Bruhat-Tits group scheme. Let ${\bf P}_{_{\ul{a}}}$ be the functor $R \to \sf{P}_{_{\ul{a}}}(R\llbracket t \rrbracket)$, which is representable by a pro-algebraic group over $k$. We call ${\bf P}_{_{\ul{a}}}(k)$ a parahoric subgroup of $G(k(\!(t)\!))$. The conjugacy classes of parahoric subgroups of $G(k(\!(t)\!))$ are classified by proper subsets of the nodes of the extended Dynkin diagram of $G$ or the facets of the Weyl alcove ${\bf a} \subset \ca_{_T}$.  These group schemes are indexed by the rational points of the alcove ${\bf a}$ which are in turn given by the types $\tau$ \eqref{localtypeandweights}.

In summary, the {\tt 2BT}-group scheme $\qh^{^G}_{_{\tau,{\sf N}^{^{(d)}}}}$ is such that it has non-trivial parahoric structures prescribed by the McKay correspondence. These are at precisely the generic points of the rational curves in $E_{_\rho}$ with further degeneration at the nodes on $E^{^{(d)}}$. The local type $\tau$ carries the information of ramification at these primes. This becomes the data for the Kawamata cover $\tilde{\sf N}(\tau) \to {\sf N}^{^{(d)}}$ \eqref{schemelikeonsurfaces} and gives the points $\ul{a}_{_j}$ in the Weyl alcove and the group schemes  $\sf{P}_{_{\ul{a}{_{_j}}}}$ at the generic points. {\em The comments after \eqref{chainrule} show that each group scheme $\qh^{^G}_{_{\tau,{\sf N}^{^{(d)}}}}$ is non-trivial on at most $\ell = \text{rank}(G)$ number of rational curves on the exceptional divisor of ${\sf N}^{^{(d)}}$}.\erem

\section{Group algebraic spaces on standard models}
The aim of this section is to build the basic group algebraic spaces on the standard models $W[d]$ defined in \cite{Li}.
\subsubsection{Standard models for a semistable curve}\label{expandeddeg} 
 
Let ${\hat\ba}^{^{d}} :=\spec~k\llbracket t_{_1}, \ldots, t_{_d} \rrbracket$ be given a $A$-scheme structure via the morphism: $t \mapsto t_{_1} \cdots t_{_d}.$  

Gieseker in \cite[Lemma 4.2, Proposition 4.1]{gies} constructs a miniversal family for the semistable curve $C^{^{(d)}}$ with fixed stable model $(C,c)$. We will however follow the detailed construction of the {\em expanded degenerations} in the paper of Jun Li \cite{Li}. 

We begin with the base family $C_{_{A}} \to A$. We let $B[0] := A$ to start the inductive construction. Thus, we have the family $C_{_{B[0]}} \to B[0]$.
Let 
\beqa\label{bdtob}
B[d-1]:=  {\hat\ba}^{^{d}}
\eeqa
with ${\hat\ba}^{^{d}}$ being given the $A$-scheme structure as above.  Jun Li  constructs the {\tt standard models} $W[d-1]$ over $B[d-1]$  inductively as a small resolution of the scheme $W[d-2]$ (see \cite[page 521]{Li} for the details). The scheme $W[d-1]$ comes with a tautological projection:
\beqa\label{junlistdmodel}
W[d-1]  \to  C_{_A} \times_{_A} \times B[d-1].
\eeqa
The fibre  of $W[d-1]$ over ${0} \in B[d-1]$ is denoted by $W[d-1]_{_{0}}$, which is isomorphic to the projective curve $C^{^{(d)}}$. 
\subsubsection{The special degeneration} For the main applications, we need the description of an \'etale neighbourhood  of $W[d-1]_{_{0}}$ in $W[d-1]$. This is done by looking at a special degeneration (see below \eqref{thekeyconstruct}).

 \cite[page 522]{Li}, Li constructs the {\em special degeneration}:
\beqa\label{specialdeg}
\Gamma[d-1] \to {\ba}^{^{d}}.
\eeqa
We briefly recall its description for our purposes (see \cite{acfw} for a nice exposition) . 

The first observation is that the fibres of $\Gamma[d-1] \to {\ba}^{^{d}}$ are not projective. Secondly, an  \'etale neighbourhood of the fibre of the origin of this morphism coincides with an  \'etale neighbourhood of $W[d-1]_{_{0}}$, which is the fibre of the origin of the morphism $W[d-1] \to B[d-1]$.

The fibre $\Gamma[d-1]_{_0} = E^{^{(d)}} \cup E'(1) \cup E'(2)$ over $0 \in B[d-1]$ is a chain of $d+1$ curves of which the first and last are $E'(1)$ and $E'(2)$ which are $\ba^{^1}$'s and the rest of the members are $\bp^{^1}$'s (\cite[Lemma 1.2, page 522]{Li}). The scheme $\Gamma[d-1]$ can be covered by $d$-open subsets $U_{_1}, \ldots, U_{_{d}}$, each of which is isomorphic to $\ba^{^{d+1}}$. If the coordinates of $U_{_\text{\cursive s}}$, which is isomorphic to the affine space $\ba^{^{d+1}}$, are denoted by $(u^{^{(\text{\cursive s})}}_{_1}, \ldots,u^{^{(\text{\cursive s})}}_{_{d+1}})$, the transition function  from $u^{^{(\text{\cursive s})}}$ to $u^{^{(\text{\cursive s}+1)}}$ on $U_{_\text{\cursive s}} \cap U_{_{\text{\cursive s}+1}}$ is given by $u^{^{(\text{\cursive s}+1)}}_{_j} = u^{^{(\text{\cursive s})}}_{_j}$ for $j = 1, \ldots \text{\cursive s}-1, \text{\cursive s}+3, \ldots, d+1$. For the three coordinates with indices $j = \text{\cursive s}, \text{\cursive s}+1, \text{\cursive s}+2$ we have the following transition relations:
\beqa\label{transrels}
u^{^{(\text{\cursive s}+1)}}_{_\text{\cursive s}} =  u^{^{(\text{\cursive s})}}_{_\text{\cursive s}}. u^{^{(\text{\cursive s})}}_{_{\text{\cursive s}+1}}
\\
u^{^{(\text{\cursive s}+1)}}_{_{\text{\cursive s}+1}} =  1/u^{^{(\text{\cursive s})}}_{_{\text{\cursive s}+1}}
\\
u^{^{(\text{\cursive s}+1)}}_{_{\text{\cursive s}+2}} = u^{^{(\text{\cursive s})}}_{_{\text{\cursive s}+1}}.u^{^{(\text{\cursive s})}}_{_{\text{\cursive s}+2}}.
\eeqa
\brem A fact which will be used later is  that the action of $G[d]:= \bg_{_m}^{^d}$ (see \cite[page 525]{Li}) on the open subsets $U_{_{\text{\cursive{s}}}}$ covering $Z[d-1]$ is given by:
\beqa\label{thegmaction0}
\sigma \cdot (u^{^{(\text{\cursive s})}}_{_1}, \ldots,u^{^{(\text{\cursive s})}}_{_{d+1}}) := (\ldots, \bar{\sigma}_{_{\text{\cursive s} - 2}} u^{^{(\text{\cursive s})}}_{_{\text{\cursive s} - 2}}, \bar{\sigma}_{_{\text{\cursive s} - 1}} u^{^{(\text{\cursive s})}}_{_{\text{\cursive s} - 1}}, {\sigma}_{_{\text{\cursive s} - 1}}^{^{-1}} u^{^{(\text{\cursive s})}}_{_{\text{\cursive s}}}, {\sigma}_{_{\text{\cursive s}}} u^{^{(\text{\cursive s})}}_{_{\text{\cursive s} +1}},\bar{\sigma}_{_{\text{\cursive s} + 1}} u^{^{(\text{\cursive s})}}_{_{\text{\cursive s} + 2}}, \dots)
\eeqa
where $\bar{\sigma}_{_{i}} = \sigma_{_i}/\sigma_{_{i-1}}$ with the convention that $\sigma_{_i} = 1$ for $i = 0,d$.\erem

\subsubsection{{\tt BT}-group schemes on the local standard model}\label{thekeyconstruct} From here onwards we work with {\sl the base change of the special degeneration $\Gamma[d-1] \to {\ba}^{^{d}}$ to the analytic neighbourhood of $0$ in $\ba^{^{d}}$}.  We  denote this by  $Z[d-1] \to {\hat\ba}^{^{d}}$  and call it the {\bf local standard model}. Thus, $Z[d-1]_{_0} = E^{^{(d)}} \cup E(1) \cup E(2) = F^{^{(d)}}$ as in \eqref{closedfibre} and $Z[d-1]$ can be identified with the analytic neighbourhood of the closed fibre $W[d-1]_{_0}$ of $W[d-1]$.

The scheme $Z[d-1]$ is smooth  with a {\sl divisor $\mathfrak D := Z[d-1] \times _{_{{A}}} {0} \subset Z[d-1]$ \eqref{keydiagzd}, with normal crossing singularities having $d+1$ irreducible components};  $\mathfrak D = \bigcup_{_{j = 0}}^{^{d}} \mathfrak D_{_j}$, where $\mathfrak D_{_{j-1}}$ and $\mathfrak D_{_{j}}$, $j = 1, \ldots, d$ are the smooth irreducible components of $\mathfrak D$ intersecting transversally along the $d+1$ disjoint,  smooth, codimension two subvariaties $\{\pmb D_{_j}\}_{_{ j = 1}}^{^{d}}$  respectively (see \cite[page 538]{Li} and \cite[page 207]{Li1}, where these are expressed  in the setting of logarithmic schemes).
This is precisely the configuration which matches the configuration of the nodal structure of the  central fibre $Z[d-1]_{_0} \simeq F^{^{(d)}}$. We have the following structure of $Z[d-1]$:
\beqa\label{keydiagzd}
\xymatrix{
{Z[d-1]_{_0}} \ar[r] \ar[d] &
 Z[d-1] \supset \mathfrak D \ar@/^/[ddr] \ar[d] \\
\{0\} \ar[r]^{origin} &  {\hat\ba}^{^{d}} \ar@{>}[dr]|{t\mapsto t_{_1}\ldots t_{_d}}\\
&& A \ni 0 }
\eeqa

Let $\qh_{_{\tau, {\sf N}^{^{(d)}}}}^{^G}$ on  be the group scheme  of local type $\tau$ on ${\sf N}^{^{(d)}}$ \eqref{keydiag1}.  Let 
\beqa\label{thegrpschemeonthecurve}
\qh_{_\tau}^{^G} := \qh^{^G}_{_{\tau,F^{^{(d)}}}}
\eeqa  
denote the restriction of $\qh_{_{\tau, {\sf N}^{^{(d)}}}}^{^G}$ to  $F^{^{(d)}}$ the closed fibre \eqref{closedfibre}.

In this subsection we prove the following result which will play a central role in what follows.

\bprop\label{localgrpscheme} Fix a group scheme $\qh_{_\tau}^{^G}$ on the fibre $Z[d-1]_{_0} \simeq F^{^{(d)}}$ of $Z[d-1] \to {\hat\ba}^{^{d}}  $. There exists a group scheme $\qh_{_{\tau, Z[d-1]}}^{^G}$ on $Z[d-1]$ which restricts to $\qh_{_\tau}^{^G}$ on $Z[d-1]_{_0}$ and which is isomorphic to {\sf the  split group scheme with fibre group $G$} on the complement of the divisor $\mathfrak D$ \eqref{keydiagzd}. \eprop 
\bpr We recall that the minimal smooth resolution  ${\sf N}^{^{(d)}} \to N_{_d}$ has a description analogous to the one given for $Z[d-1]$. The closed fibre \eqref{closedfibre} $F^{^{(d)}}$ is identified with $Z[d-1]_{_0}$.

The scheme ${\sf N}^{^{(d)}}$ can be covered by $d$-open subsets $B_{_1}, \ldots, B_{_{d}}$, each of which is isomorphic to $\ba^{^{2}}$. If the coordinates of $B_{_\text{\cursive s}}$ are denoted by $(a_{_\text{\cursive s}}, b_{_{\text{\cursive s}}})$ in its identification with $\ba^{^{2}}$, the transition functions  from $(a_{_\text{\cursive s}}, b_{_{\text{\cursive s}}})$ to $(a_{_{\text{\cursive s}+1}}, b_{_{{\text{\cursive s}+1}}})$ on $B_{_\text{\cursive s}} \cap B_{_{\text{\cursive s}+1}}$ are given by 
\beqa\label{transrels2}
a_{_{\text{\cursive s}+1}} =  a_{_\text{\cursive s}}^{^2} b_{_{\text{\cursive s}}}
\\
b_{_{\text{\cursive s}+1}} =  1/a_{_{\text{\cursive s}}}
\eeqa
We firstly show that there is a natural morphism 
\beqa\label{level d-1}
q_{_d}: Z[d-1] \to {\sf N}^{^{(d)}}
\eeqa
which has the following key property:  for each $j$ the intersection of the divisor $\mathfrak D \subset Z[d-1]$ with the open set $U_{_j}$ is  mapped by $q_{_d}$  precisely to the intersection of the  normal crossing divisor \eqref{closedfibre} $F^{^{(d)}} \subset {\sf N}^{^{(d)}}$, with the open set $B_{_j}$.

We begin by defining the morphism for $d = 2$:
\beqa\label{level1}
q_{_2}: Z[1] \to B_{_1} \cup B_{_2} = N^{^{(2)}}
\eeqa
by sending (for $\text{\cursive s} = 1,2$), 
\beqa
(u^{^{(\text{\cursive s})}}_{_1},u^{^{(\text{\cursive s})}}_{_{2}}, u^{^{(\text{\cursive s})}}_{_{3}}) \mapsto (a_{_\text{\cursive s}}, b_{_\text{\cursive s}})
\eeqa 
where
\beqa
a_{_1} = u^{^{(1)}}_{_{2}},  ~~
b_{_1} = u^{^{(1)}}_{_1} u^{^{(1)}}_{_3}, \\
a_{_2} = u^{^{(2)}}_{_1} u^{^{(2)}}_{_3},  ~~ b_{_2} = u^{^{(2)}}_{_{2}}
\eeqa

We now define the morphism \eqref{level d-1}
recursively as follows: if $\phi_{_\text{\cursive s}}:U_{_\text{\cursive s}} \to B_{_\text{\cursive s}}$ sends $\phi_{_\text{\cursive s}}(u^{^{(\text{\cursive s})}}_{_1}, \ldots,u^{^{(\text{\cursive s})}}_{_{d+1}}) = (a_{_\text{\cursive s}}, b_{_\text{\cursive s}})$, then recursively
$\phi_{_{\text{\cursive s} +1}}(u^{^{(\text{\cursive s}+1)}}_{_1}, \ldots,u^{^{(\text{\cursive s}+1)}}_{_{d+1}}) =  (a_{_{\text{\cursive s}+1}}, b_{_{\text{\cursive s} +1}})$, where 
\beqa
a_{_{\text{\cursive s}+1}}(u^{^{(\text{\cursive s}+1)}}_{_1}, \ldots,u^{^{(\text{\cursive s}+1)}}_{_{d+1}}) = a_{_\text{\cursive s}}^{^2}.b_{_\text{\cursive s}}\\ 
b_{_{\text{\cursive s}+1}}(u^{^{(\text{\cursive s}+1)}}_{_1}, \ldots,u^{^{(\text{\cursive s}+1)}}_{_{d+1}}) = 1/a_{_\text{\cursive s}}.
\eeqa The explicit formulae can be obtained by using the transition data \eqref{transrels}. For example, the maps $U_{_j} \to B_{_j}, j = 1,2,3$, look as follows: 
\beqa
a_{_1}(u^{^{(1)}}_{_1}, \ldots,u^{^{(1)}}_{_{d+1}}) = u^{^{(1)}}_{_2}, ~~~~~~
b_{_1}(u^{^{(1)}}_{_1}, \ldots,u^{^{(\ell)}}_{_{d+1}}) = u^{^{(1)}}_{_1} u^{^{(1)}}_{_3} \\
a_{_2}(u^{^{(2)}}_{_1}, \ldots,u^{^{(2)}}_{_{d+1}}) = u^{^{(2)}}_{_1} u^{^{(2)}}_{_3}, ~~~~ 
b_{_2}(u^{^{(2)}}_{_1}, \ldots,u^{^{(2)}}_{_{d+1}}) = u^{^{(2)}}_{_2}\\
a_{_3}(u^{^{(3)}}_{_1}, \ldots,u^{^{(3)}}_{_{d+1}}) = \big(u^{^{(3)}}_{_1}\big)^{^2} u^{^{(3)}}_{_2}. \big(u^{^{(3)}}_{_3}\big)^{^{-1}}, ~~~~
b_{_3}(u^{^{(3)}}_{_1}, \ldots,u^{^{(3)}}_{_{d+1}}) = \big(u^{^{(3)}}_{_1}\big)^{^{-1}}.u^{^{(3)}}_{_3}
\eeqa
That the morphism $q_{_d}$ has the stated property can be seen by a simple check using \cite[Lemma 1.2(i)]{Li} and the computations above.

We define:
\beqa
\qh_{_{\tau, Z[d-1]}}^{^G} := q_{_d}^{^*}(\qh_{_{\tau, {\sf N}^{^{(d)}}}}^{^G})
\eeqa
The group scheme $\qh_{_{\tau, Z[d-1]}}^{^G}$ can be easily seen to satisfy the stated properties. 
\epr

\subsubsection{The Kawamata cover of $Z[d-1]$}\label{sncinzd}

\bth\label{schemelikeonsurfaces1} The group scheme  $\qh_{_{\tau, Z[d-1]}}^{^G}$ on $Z[d-1]$ can be realized as invariant direct image \eqref{nsgrpscheme1.5} from a global smooth ramified covering $\tilde{Z}_{_\tau}[d-1] \to Z[d-1]$.
\eeth
\bpr
By \cite[Lemma 1.2(i)]{Li} and the computations above it follows by a simple check that for each $j$ the intersection of the divisor $\mathfrak D \subset Z[d-1]$ with the open set $U_{_j}$  maps precisely to the intersection of the  normal crossing divisor \eqref{closedfibre} $F^{^{(d)}} \subset {\sf N}^{^{(d)}}$, with the open set $B_{_j}$. By the general theory, \eqref{kawa} (\cite{kawamata} or \cite[Lemma 2.5]{vieweg}), there exist a Kawamata covering of $Z[d-1]$ with the prescribed ramification data. 

Let ${\tilde{\sf N}}(\tau) \to {\sf N}^{^{(d)}}$ be as in \eqref{schemelikeonsurfaces} (see also \eqref{briefdescript}). Pulling back by the morphism $Z[d-1] \to {\sf N}^{^{(d)}}$ and taking reduced scheme structure we get a covering $Z" \to Z[d-1]$. The scheme $Z"$ need not be smooth but which is \'etale over the complement of the divisor $\mathfrak D$.  This can be rectified by \cite[proof of Corollary 2.6]{vieweg}, where we can get the required Kawamata covering $\tilde{Z}_{_\tau}[d-1]$ of $Z[d-1]$ as a finite covering of $Z"$ (see \eqref{schemelikeonsurfaces} for similar arguments). The ramification data is controlled by the local type $\tau$. More importantly, we can conclude as in \eqref{schemelikeonsurfaces}, that we can realize the group scheme $\qh_{_{\tau, Z[d-1]}}^{^G}$ on $Z[d-1]$ as an invariant direct image \eqref{nsgrpscheme1.5} of an equivariant group scheme on the Kawamata covering $\tilde{Z}_{_\tau}[d-1]$.\epr
\subsubsection{Group algebraic spaces on $W[d-1]$\label{kawamatafortheversal}}  As in \eqref{sncinzd}, by \cite[page 538]{Li}, we see that $W[d-1]$ with its tautological projection \eqref{junlistdmodel} also has a canonically defined divisor $\mathfrak D := W[d-1] \times _{_{A}} 0$ defined by \eqref{keydiagzd}, with normal crossing singularities and we have a Kawamata covering $\tilde{W}_{_\tau}[d-1] \to W[d-1]$ with the precise ramification data (given by $\tau$) at the generic points of the smooth irreducible components $\{\mathfrak D_{_j}\}$ \eqref{kawa}. We can obtain {\tt group algebraic spaces $\qh_{_{\tau, W[d-1]}}^{^G}$ of local type $\tau$} by gluing the group scheme with constant fibre $G$ away from the divisor with $\qh_{_{\tau, Z[d-1]}}^{^G}$. The gluing is done using the pull-back of $E_{_A}^{^o}(G)$ on $C_{_A} - c$ as in \eqref{fixsomething}. It is not hard to check that there is a group algebraic space $\mathcal G_{_{\tilde{W}}}$ on the covering $\tilde{W}_{_\tau}[d-1]$ such that the group algebraic space $\qh_{_{\tau, W[d-1]}}^{^G}$ is obtained as an invariant direct image \eqref{nsgrpscheme1.5} of $\mathcal G_{_{\tilde{W}_{_\tau}}}$ as in \eqref{schemelikeonsurfaces1}.
\brem\label{typeandlength1} By \eqref{typeandlength} we see that the Kawamata cover $\tilde{W}_{_\tau}[d-1]$ of $W[d-1]$ is ramified at the generic points of at most $\ell$ of the divisors $\mathfrak D_{_j}$'s. Likewise, the group algebraic space $\qh_{_{\tau, W[d-1]}}^{^G}$ {\sl has non-trivial parahoric structure at the generic points of at most $\ell$ components of the divisor $\mathfrak D$}.\erem

\section{Admissible pairs on standard models}
The aim of this section is to define admissible pairs consisting of basic group algebraic spaces and a special class of torsors on the standard models $W[d]$. This is done as before by building it first on analytic surfaces $N^{^{(d)}}$ and from there to more general objects.  

\subsubsection{Quasi-admissibility of vector bundles and McKay correspondence} In this subsection we assume $G  = \text{GL}(n)$.  The aim is to understand the notion of {\em (quasi)admissibility} of vector bundles \eqref{Giesekervectorbundleappendix}  as an outcome of the McKay correspondence. This gives the new perspective on Gieseker's objects.

\bdefe\label{newadmissibility}  A vector bundle  $\eV$ of rank $n$ on the smooth surface ${\sf N}^{^{(d)}}$ is called {\em quasi-admissible} (see \eqref{Giesekervectorbundleappendix}) if the restriction $\eV\mid_{_{E^{^{(d)}}}}$ to the chain $E^{^{(d)}} \subset F^{^{(d)}}$ in the closed fibre \eqref{closedfibre}  is {\em  standard} \eqref{Giesekervectorbundleappendix},  and furthermore the direct image $p_{_{d,*}}(\eV)$ is torsion-free on  $N_{_d}$, where $p_{_{d}}$ is as in \eqref{keydiag2}.\edefe

We consider representations $\rho:{\sf\mu_{_d}} \to \text{GL}(n)$ of type $\tau$.  We fix a maximal torus $T_{_n} \subset \text{GL}(n)$ (which can be taken as the $(n \times n)$-diagonal matrices by choice of a basis).  Recall that the isomorphism classes of $({\sf\mu_{_d}}, \text{GL}(n))$-bundles on $D$ are classified by the equivalence classes of representations $\rho:{\sf\mu_{_d}} \to T_{_n}$, which we term the local type $\tau$ of $\rho$. In particular, this is given by writing $\rho$ as $\rho(\gamma) = \text{diag}(\zeta^{^{m_{_1}}}, \ldots, \zeta^{^{m_{_j}}})$, where $\gamma$ is a generator of the cyclic group ${\sf\mu_{_d}}$, and $\zeta$ is the primitive $d^{^{th}}$-root of unity defined by $\gamma \cdot z = \zeta.z$, with $z = (u,v)$ being local coordinates in $D$, and $0 \leq m_{_1} < m_{_2} \ldots < m_{_j} \leq d - 1$. Each $m_{_i}$ repeats $r_{_i}$-times so that $\sum_{_i} r_{_i} = n$.

The basic motivation for our study is to relate this definition to the notion of quasi-admissibility of vector bundles on the smooth surface ${\sf N}^{^{(d)}}$. Recall the basic diagram:
\beqa\label{keydiag2}
\xymatrix{
{D^{^{(d)}}} \ar[r]^{f} \ar[d]_{q} &
 {\sf N}^{^{(d)}}  \ar[d]_{p_{_d}} \\
D \ar[r]^{\sigma} &  N_{_d} \\
}
\eeqa
\bprop\label{mckayadm} ("Fourier-Mukai") A vector bundle $\mathcal V$ on  ${\sf N}^{^{(d)}}$  is quasi-admissible if and only if there is a ${\sf\mu_{_d}}$-equivariant vector bundle $V$ on $D$ arising from a  representation $\rho$ such that $\mathcal V \simeq ({\tt Inv}\circ f_{_*})(q^{^*}(V))$. In particular, the vector bundle $\mathcal V$ is non-trivial on {\em at most} $n$ rational curves in the exceptional divisor of ${\sf N}^{^{(d)}}$. \eprop
\bpr 

Let  $V$ be a  ${\sf\mu_{_d}}$-vector bundle on $D$ coming from a representation $\rho:\smu \to \text{GL}(n)$ of local type $\tau$. The invariant direct image ${\tt Inv}\circ \sigma_{_*}(V)$ of $V$ under $\sigma:D \to N_{_d}$ gives a reflexive sheaf $\cf$ on the normal surface $N_{_d}$ which takes the form:
\beqa\label{ideal}
\cf \simeq  \bigoplus_{i = 1}^{j} \big(x, t^{^{m_i}}\big)^{^{\oplus r_{_i}}}
\eeqa 
Each summand  is the invariant push-forward of the equivariant line bundle associated to a character of ${\sf\mu_{_d}}$ given by the $m_{_i}$'s with multiplicity $r_{_i}$ in terms of the representation $\rho$. 

In \cite[Proposition 4.2]{bbn} it was shown that quasi-admissible vector bundles on ${\sf N}^{^{(d)}}$ \eqref{newadmissibility} are {\em precisely} those with the property that the direct image $ (p_{_d})_{_*}(\eV)$ is the reflexive sheaf $\cf$ as in \eqref{ideal}.  Indeed, we can identify the quasi-admissible vector bundle with  $p_{_d}^{^*}(\cf)/(\text{tors})$.

On the other hand, the result of Gonzalez-Springberg and Verdier \cite[Theorem 2.2]{gonverd} proves that   any reflexive sheaf $\cf$ on $N_{_d}$ can be expressed as ${\tt Inv}\circ \sigma_{_*}(V)$, for $V$ a vector bundle on $D$ of local type $\tau$ and conversely (see also \cite{ns2} and \cite{bbn}). Moreover,  one has the isomorphism \cite[Theorem 2.2, Proposition 2.8]{gonverd}:
\beqa\label{gonspr}
({\tt Inv}\circ f_{_*})(q^{^*}(V)) \simeq p_{_d}^{^*}(\cf)/(\text{tors}).
\eeqa
Whence by \eqref{gonspr}, the process $({\tt Inv}\circ f_{_*})(q^{^*}(V))$ gives quasi-admissible  vector bundles on ${\sf N}^{^{(d)}}$ and conversely. \epr

\brem\label{spread} Let $V$ be a quasi-admissible vector bundle  on $F^{^{(d)}}$, i.e. $V$ is firstly standard and further, if $p:F^{^{(d)}} \to E_{_c}$ be the projection to the analytic neighbourhood of $c \in C$, the sheaf $p_{_*}(V)$ is torsion-free on the reducible curve $E_{_c} \subset N_{_d}$. At the node, $p_{_*}(V)$ gets a decomposition $\co^{^a} \oplus \mathfrak m^{^b}$ for some $a, b \geq 0$. Hence it is the restriction of a reflexive sheaf $\cf$ on $N_{_d}$ given as in \eqref{ideal}, where the $m_{_i}$ come from the rational component $E_{_\chi}$ labelled by the character $\chi:\zeta \mapsto \zeta^{^{m_{_i}}}$ as dictated by McKay correspondence and the multiplicity $r_{_i}$ is precisely the number of copies of $\co(1)$ in the restriction of $V$ to  $E_{_\chi}$. Whence, $V$ is isomorphic to the restriction of the quasi-admissible vector bundle $\eV$ on ${\sf N}^{^{(d)}}$. 
\erem

\subsubsection{Quasi-admissible vector bundles on $Z[d-1]$ and parabolic structures} Let $\eV$ be a vector bundle on the local standard model $Z[d-1] \to {\hat\ba}^{^d}$. Then $\eV$ is called quasi-admissible \eqref{Giesekervectorbundleappendix} if it is so on each closed fibre. The fibre over $0_{{{\hat\ba}^{^d}}}$ is the chain $F^{^{(d)}}$ with $d+1$-components, $d-1$ of which are $\bp^{^1}$'s and the deformations are by smoothing of nodes. A quasi-admissible bundle $\eV_{_0}$ has no deformations along the fibre since the summands are $\co$ or $\co(1)$. Therefore the lift $\eV$, in an analytic or \'etale neighbourhood of the central fibre, is uniquely determined by $\eV_{_0}$. 
By following the construction in \eqref{thekeyconstruct}, we can construct the vector bundle $\eV$ on $Z[d-1]$ knowing $\eV_{_0}$.  In other words, we can realize $\eV$ as the invariant direct image of an equivariant vector bundle of local type determined by $\eV_{_0}$,  on the Kawamata cover $\tZ_{_\tau}[d-1]$ of $Z[d-1]$. This model of local type determined by $\eV_{_0}$ on the local standard model then allows us to construct quasi-admissible bundles on the projective family given by the standard model $W[d-1]$ by gluing the local quasi-admissible ones with bundles on the complement of the divisor $\mathfrak D \subset W[d-1]$.  Whence all quasi-admissible vector bundles on $W[d-1]$ with fixed local type $\tau$ can also be obtained as invariant direct images of certain equivariant vector bundles on $\tW_{_\tau}[d-1]$.

\bcor\label{parabstr} Let $\eV$ be a quasi-admissible vector bundle on $W[d-1]$. Then $\eV$ gets a canonical parabolic structure at the generic points of the s.n.c divisor $\mathfrak D \subset W[d-1]$ with weights determined by the local type $\tau$ of $\eV_{_0}$.  The number of components at the generic points of which $\eV$ has non-trivial parabolic structure is {\em bounded above by the rank of $\eV$}. \ecor
\bpr Since  $\eV$ comes as an invariant direct image of an equivariant bundle on the Kawamata cover,  by \cite{biswas} it gets canonical parabolic structure at the generic points of $\mathfrak D_{_j}$ and the number of these components is bounded by the rank of $\eV$.\epr

\subsubsection{Admissible pairs for the $\text{GL}(n)$ case}
Our first task is to get a "group scheme" plus "torsor" equivalent of these definitions in the case when $G = \text{GL}(n)$. Fix a $({\sf\mu_{_d}}, \text{GL}(n))$-torsor  $\eT_{_{D}}$ of local type $\tau$, and let $\qh^{^{\text{GL}(n)}}_{_{\tau,{\sf N}^{^{(d)}}}}$ be as in \eqref{nsgrpscheme1.5} and  $\qh^{^{\text{GL}(n)}}_{_{\tau,S^{^{(d)}}}}$ the group algebraic space on $S^{^{(d)}}$ (coming from a fixed gluing). We keep the diagram \eqref{keydiag1} in mind for the next definitions. Let $E$ be a $\eT_{_{D}}(\text{GL}(n))$-torsor on $D$. Then the pull-back $q^{^*}(E)$ gets the structure of an $E(G,\tau)$-torsor on $D^{^{(d)}}$ \eqref{egtau}.

\bdefe\label{admissibletorsor}   Let $\eE$ be a $\qh^{^{\text{GL}(n)}}_{_{\tau,{\sf N}^{^{(d)}}}}$-torsor  on ${\sf N}^{^{(d)}}$. The pair $(\qh^{^{\text{GL}(n)}}_{_{\tau,{\sf N}^{^{(d)}}}}, \eE)$ is called {\em admissible}  if $\eE$ arises as $({\tt Inv}\circ f_{_*})(q^{^*}(E))$ for a $\eT_{_{D}}(\text{GL}(n))$-torsor $E$ on $D$. A pair $(\qh^{^{\text{GL}(n)}}_{_{\tau,S^{^{(d)}}}}, \eE)$ on  $S^{^{(d)}}$ is called admissible if its restriction to ${\sf N}^{^{(d)}}$ is so. \edefe

\brem\label{fromjag} By \cite[Theorem 4.1.6]{base}, the process of taking "invariant direct images" extends to torsors as well and gives an isomorphism of stacks. In other words, $({\tt Inv}\circ f_{_*})(q^{^*}(E))$ above does give a $\qh^{^{\text{GL}(n)}}_{_{\tau,S^{^{(d)}}}}$-torsor and this process can be carried out for all groups $G$, not merely $\text{GL}(n)$. \erem
\brem By using the Kawamata covering $h:\tilde{S}(\tau) \to S^{^{(d)}}$ (with local ramification given by $\tau$) \eqref{algspacelikeonsurfaces} we can check that in an admissible pair $(\qh^{^{\text{GL}(n)}}_{_{\tau,S^{^{(d)}}}}, \eE)$ the torsor $\eE$ can be realized as  $({\tt Inv}\circ h_{_*})(E'))$ for an equivariant torsor $E'$ on $\tilde{S}(\tau)$ for a group algebraic space with fibre $G$. Moreover, $E'$ agrees with $E$ on the inverse image of ${\sf N}^{^{(d)}}$. \erem

\noindent
{\it A parabolic vector bundle on ${\sf N}^{^{(d)}}$.}

\bdefe\label{assocvb} Let  $(\qh^{^{\text{GL}(n)}}_{_{\tau,{\sf N}^{^{(d)}}}}, \eE)$ be an admissible pair on ${\sf N}^{^{(d)}}$.
The vector bundle $\eE^{^{par}}(k^{^n})$  {parabolically}-associated to $\eE$ is defined as:
\beqa
\eE^{^{par}}(k^{^n}) := {\tt Inv}\circ f_{_*}(q^{^*}(E(k^{^n}))
\eeqa
Let $(\qh^{^{\text{GL}(n)}}_{_{\tau,S^{^{(d)}}}}, \eE)$ be an admissible pair on $S^{^{(d)}}$. Define the vector bundle $\eE^{^{par}}(k^{^n}) := {\tt Inv}\circ h_{_*}(E'(k^{^n}))$ with $E'$ on $\tilde{S}(\tau)$. 
\edefe
The vector bundle $\eE^{^{par}}(k^{^n})$ defined above does indeed come with canonical parabolic structures in the sense of Seshadri. These structures are at the generic points of the rational components of $C^{^{(d)}}\subset S^{^{(d)}}$. We note that it is this phenomenon  which makes sure that the line bundles in the decomposition along the chain of $\bp^{^1}$'s in $E^{^{(d)}}$ remains $\co$ or $\co(1)$.  The terminology comes from the similarity of the phenomenon with the process of taking "parabolic tensor products" of parabolic vector bundles on curves.  Taking usual tensor products will increase the degree of the line bundles on the rational components.

\subsubsection{Equivalence of various notions of admissibility}  
The next result ties  up the various notions  of {admissibility}.  
\bprop\label{admingies}({\em"Fourier-Mukai" on group schemes and torsors}) Let    $(\qh_{_{\tau,{\sf N}^{^{(d)}}}}^{^{\text{GL}(n)}},\eE)$ be an admissible pair on ${\sf N}^{^{(d)}}$ \eqref{admissibletorsor}. Then the {parabolically}-associated vector bundle $\eE^{^{par}}(k^{^n})$  is a quasi-admissible  vector bundle on ${\sf N}^{^{(d)}}$ \eqref{newadmissibility}. Conversely, if $\eV$ is a quasi-admissible  vector bundle of rank $n$ on ${\sf N}^{^{(d)}}$, then there exists  an admissible pair $(\qh_{_{\tau,{\sf N}^{^{(d)}}}}^{^{\text{GL}(n)}},\eE)$  on ${\sf N}^{^{(d)}}$ and an admissible torsor $\eE$ such that $\eV = \eE^{^{par}}(k^{^n})$.\eprop

\bpr This is immediate from \eqref{mckayadm} and the fact that giving a ${\sf\mu_{_d}}$-vector bundle on $D$ coming from a representation $\rho$ is equivalent to giving a  $\ce_{_D}(\text{GL}(n))$-torsor $E$ and $\eV = \eE^{^{par}}(k^{^n})$. \epr

\bcor\label{admingies1} A vector bundle $\eV$ on $S^{^{(d)}}$ is quasi-admissible if and only if $\eV \simeq \eE^{^{par}}(k^{^n})$ for an admissible pair $(\qh^{^{\text{GL}(n)}}_{_{\tau,S^{^{(d)}}}}, \eE)$. \ecor

\subsubsection{Torsors and vector bundles on $W[d-1]$}  Let $\qh_{_{\tau, Z[d-1]}}^{^{\text{GL}(n)}}$ be the group scheme constructed on the local standard model $Z[d-1]$ \eqref{thekeyconstruct}. By \eqref{schemelikeonsurfaces1}, these group schemes can be realized as invariant direct images \eqref{nsgrpscheme1.5} of group schemes $\mathcal G_{_{\tilde{Z}}}$ with fibres $G = \text{GL}(n)$, from the Kawamata cover $\kappa:\tilde{Z}_{_\tau}[d-1] \to Z[d-1]$. Let $\ce$ be a $\qh_{_{\tau, Z[d-1]}}^{^{\text{GL}(n)}}$-torsor. Then by \cite[Theorem 4.1.6, page 24]{base}, there exist a unique $\mathcal G_{_{\tilde{Z}}}$-torsor $E'$ on $\tilde{Z}_{_\tau}[d-1]$ such that the invariant push-forward of $E'$ is $\ce$. Hence we can define the {\em parabolically associated} vector bundle $\ce^{^{par}}(k^{^n}):=\text{Inv} \circ \kappa_{_*}(E'(k^{^n}))$. Observe that all these considerations make sense on the standard model $W[d-1]$.

Following \eqref{admingies}, we can therefore have the following definition. 
\bdefe\label{admonstdmodels}  The pair $(\qh_{_{\tau, Z[d-1]}}^{^{\text{GL}(n)}}, \ce)$  is called admissible if  the {parabolically}-associated vector bundle $\ce^{^{par}}(k^{^n})$ is quasi-admissible on $Z[d-1]$. Likewise,  for any group algebraic space $\qh_{_{\tau, W[d-1]}}^{^{\text{GL}(n)}}$ of local type $\tau$ on $W[d-1]$, a pair $(\qh_{_{\tau, W[d-1]}}^{^{\text{GL}(n)}}, \ce)$   is called admissible if the {parabolically}-associated vector bundle $\ce^{^{par}}(k^{^n})$ is quasi-admissible.\edefe

\subsubsection{Admissible pairs for general $G$}

Let $E$ be a $\eT_{_{D}}(G)$-torsor on $D$. Then the pull-back $q^{^*}(E)$ gets the structure of an $E(G,\tau)$-torsor on $D^{^{(d)}}$ \eqref{egtau}. Let $\qh^{^{G}}_{_{\tau,{\sf N}^{^{(d)}}}}$ be the  group scheme of local type $\tau$ on ${\sf N}^{^{(d)}}$.
\bdefe\label{admissibletorsor}   Let $\eE$ be a $\qh^{^{G}}_{_{\tau,{\sf N}^{^{(d)}}}}$-torsor  on ${\sf N}^{^{(d)}}$. A pair $(\qh^{^{G}}_{_{\tau,{\sf N}^{^{(d)}}}}, \eE)$ on ${\sf N}^{^{(d)}}$ is called {\em admissible}  if $\eE$ arises as $({\tt Inv}\circ f_{_*})(q^{^*}(E))$ \big(compare \eqref{fromjag} where $G = \text{GL}(n)$, and also \cite[Theorem 4.1.6]{base}\big). A pair $(\qh^{^{G}}_{_{\tau,S^{^{(d)}}}}, \eE)$ on  $S^{^{(d)}}$ is called admissible if its restriction to ${\sf N}^{^{(d)}}$ is so. \edefe

{\sl The notion of admissibility behaves well under extension of structure groups}. More precisely, 
let  $\eta:G \hra \text{GL}(n)$ be a faithful representation. Fix a maximal torus $T_{_n} \subset \text{GL}(n)$ such that $\eta:T \hra T_{_n}$.
Observe that the group scheme $\qh^{^{G}}_{_{\tau,{\sf N}^{^{(d)}}}}$ comes with a canonical inclusions 
\beqa\label{canonicalinclu}
\eta: \qh^{^{G}}_{_{\tau,{\sf N}^{^{(d)}}}} \hra \qh^{^{ \text{GL}(n)}}_{_{\tau,{\sf N}^{^{(d)}}}}, ~~~~
\eta: \qh^{^G}_{_{\tau,S^{^{(d)}}}} \hra \qh^{^{\text{GL}(n)}}_{_{\tau,S^{^{(d)}}}}.
\eeqa
by taking invariant push-forwards of inclusions of group schemes (with fibre $G$ and $\text{GL}(n)$) induced by $\eta$, via the Kawamata covering $\tilde{\sf N}(\tau) \to {\sf N}^{^{(d)}}$ (resp. $\tilde{S}(\tau) \to S^{^{(d)}}$ ) \eqref{schemelikeonsurfaces}, \eqref{algspacelikeonsurfaces} we get the following statement which is easily seen using \eqref{admingies}:
\blem\label{insight} A pair $(\qh^{^{G}}_{_{\tau,{\sf N}^{^{(d)}}}}, \eE)$ is admissible if and only if the associated pair $(\qh^{^{\text{GL}(n)}}_{_{\tau,{\sf N}^{^{(d)}}}}, \eta_{_*}(\eE))$ is so for any $\eta$. \elem

{\sl We now work with $Z[d-1]$ and $W[d-1]$ just as we did with ${\sf N}^{^{(d)}}$ and $S^{^{(d)}}$}. By the construction of the group scheme $\qh_{_{\tau, Z[d-1]}}^{^G}$  on  $Z[d-1]$, it is clear that $\eta$ induces a canonical inclusions 
\beqa\label{etaonstd}
\qh_{_{\tau, Z[d-1]}}^{^G} \hra \qh_{_{\tau, Z[d-1]}}^{^{\text{GL}(n)}}, ~~~~
\qh_{_{\tau, W[d-1]}}^{^G} \hra \qh_{_{\tau, W[d-1]}}^{^{\text{GL}(n)}}.
\eeqa

Modelling after \eqref{insight} and using \eqref{admonstdmodels}, we make the following:
\bdefe\label{etaonstd1}Say that a pair $(\qh_{_{\tau, W[d-1]}}^{^G}, \ce)$ is {\em admissible} if the associated pair $(\qh_{_{\tau, W[d-1]}}^{^{\text{GL}(n)}}, \eta_{_*}(\ce))$ is admissible \eqref{admonstdmodels}, i.e. the {parabolically}-associated vector bundle $\ce^{^{par}}(k^{^n})$ \footnote{more accurately $(\eta_{_*}(\ce))^{^{par}}(k^{^n})$} is a quasi-admissible vector bundle on $W[d-1]$.\edefe

\section{Admissible pairs on general Modifications}
Modifications or expanded degenerations of curves have been used by Gieseker and others to study degenerations of moduli spaces of vector bundles on smooth curves. 

In the previous section we defined and studied the notion of admissibility of pairs on {\sl standard models} $W[d]$ which are special examples of modifications but having certain local versal properties. The aim of  this section is to work with very general modifications of $C_{_A}$ and define admissible pairs on them. 
\bdefe\label{modification}(cf. \cite[Definition 1.9, page 531]{Li}, see also \cite[Definition 3.8]{kausz} and \cite{zollstocks}) For every \text{${A}$-scheme} $T$, a {\em modification or an expanded degeneration} of $C_{_T}$  over $T$ is a pair $(\tt M, \pi)$, where $\tt M$ is a flat family of projective curves over $T$ together with a $T$-projection:  
\beqa\label{modification1}
{\tt M} \stackrel{\pi} \to C_{_T} = C_{_{A}} \times_{_{A}} T
\eeqa
with the following property:  there is an open covering $T_{_\alpha}$ of  $T$ in the \'etale topology and morphisms $\zeta_{_\alpha}:T_{_\alpha} \to B[d_{_\alpha}]$ which induces an isomorphism $\xi_{_\alpha}: {\tt M}_{_{\alpha}}  \stackrel{\simeq} \to W[d_{_\alpha}] \times_{_{B[d_{_\alpha}]}} T_{_{\alpha}}$ compatible with the projection $\pi:{\tt M} \to C_{_T}$, where ${\tt M}_{_{\alpha}} := {\tt M} \times _{_T} T_{_{\alpha}}$.
\edefe
An arrow $\tt M \to \tt M'$ consists of a $A$-morphism $T \to T'$ and an $T$-isomorphism $\tt M \to \tt M' \times_{_{T'}} T$ which is compatible with  their tautological projections to $C_{_T}$.

Two modifications $\tt M_{_T}$ and $\tt M'_{_T}$ are isomorphic if there is a $T$-isomorphism $\tt M_{_T} \to \tt M'_{_T}$ compatible with the projections $\tt M_{_T} \to C_{_T}$ and $\tt M'_{_T} \to C_{_T}$. The groupoid of modifications is a {\em stack} (see \cite[Proposition 1.10]{Li}).

\subsubsection{Universal construction when $G = \text{\text{GL}(n)}$}\label{univconstructs}
We work in the setting of vector bundles on modifications. Here we rely on \cite{ns1} and \cite{kausz}. This has been summarized with some small variations in the appendix \eqref{GiesekerfunctorreltoR} below. The notations are as in the appendix.

For vector bundles of rank $n$, by the representability of the functor ${\eG}^{^\mathtt q}_{_{N}}$ (see  \eqref{GiesekerfunctorreltoR} and \cite[Proposition 8]{ns2}), we have a {\em flat}  $A$-scheme ${\mathtt Y} \to A$ which is quasi-projective and regular over $k$ and  which is obtained as a  $\text{PGL}(N)$-invariant open subscheme of a suitable Hilbert scheme. Its generic fibre is smooth and the closed fibre $\mathtt Y_{_o} \subset {\mathtt Y}$   over $o \in A$ is a reduced divisor   with normal crossing singularities. Indeed, $\mathtt Y_{_o}$ is known to be irreducible (see \cite[Proposition 8 and Page 200]{ns1}).

We also get a universal modification $\mathbb M \to {\mathtt Y}$ of chain length bounded by the rank $n$. Let $\{{\sf M}_{_\alpha} \to {\tt Y}_{_\alpha}\}$ togethere with the morphisms $\zeta_{_{\alpha}}, \xi_{_{\alpha}}$ provided by \eqref{modification} define the cover for the modification $\mathbb M$ in terms of the $W[d_{_{\alpha}}]$.

 In fact, at each closed point of $x \in \mathtt Y_{_o}$ corresponding to a semistable curve $C^{^{(d)}}$, there is an \'etale neighbourhood $T(x)$ and a smooth morphism $T(x) \to B[d-1]$ such that the pull back of the versal space $W[d-1]$ is isomorphic to the restriction of $\sf M$ to the neighbourhood. Since the $W[d]$'s have a divisor with simple normal crossing singularities, we see that the universal family of curves $\sf M$ also has a simple normal crossing divisor $\sf M_{_{\mathtt Y_{_o}}}$ which is \'etale locally the pull-back of the normal crossing divisor   on the $W[d]$ \eqref{kawamatafortheversal}.

\bprop\label{theunivconstructvb} Let $\mathbb V$ be the universal  quasi-admissible vector bundle  on $\mathbb M$ of rank $n$ (see \eqref{GiesekerfunctorreltoR} and the remarks there). There exists a group algebraic space ${\qh}^{^{\text{GL}(n)}}_{_{\tt univ}}$ together with a  torsor $\eP$ on $\mathbb M$ such that $({\qh}^{^{\text{GL}(n)}}_{_{\tt univ}}, \eP)$ gives the "universal" admissible pair, i.e.:
\beqa\label{universalgl} \eP^{^{par}}(k^{^n}) \simeq \mathbb V.
\eeqa
\eprop

\bpr By \cite[proof of Proposition 4.1, page 183]{gies}, the restriction ${\mathbb V}_{_{\alpha}}$ to $\mathbb M_{_{\alpha}}$ of the universal quasi-admissible bundle $\mathbb V$,   is the pull-back $= \xi_{_\alpha}^{^*}(V_{_\alpha})$ of an admissible vector bundle $V_{_{\alpha}}$ on $W[d_{_\alpha}]$. The local type of $V_{_{\alpha}}$ dictates the choice of the Kawamata cover $\tilde{W}_{_{\alpha}} \to W[d_{_\alpha}]$  which realizes $V_{_{\alpha}}$ as the invariant direct image of an equivariant vector bundle on  $\tilde{W}_{_{\alpha}}$ \eqref{kawamatafortheversal}. This can be done for each $\alpha$.

By \cite{gies}, in the case of the universal space $\mathbb M$,  the morphisms $\zeta_{_\alpha}:{\tt Y}_{_\alpha} \to B[d_{_{\alpha}}]$ and the induced $\xi_{_{\alpha}}:\mathbb M_{_{\alpha}} \to W[d_{_{\alpha}}]$ are smooth morphisms. Hence, 
the pull-backs $\tilde{\sf M}_{_{\alpha}} := {\xi_{_{\alpha}}}^{^*}(\tilde{W}_{_{\alpha}})$ are such that $\kappa_{_\alpha}:\tilde{\mathbb M}_{_{\alpha}} \to \mathbb M_{_{\alpha}}$ give finite covers of the smooth quasi-projective schemes $\mathbb M_{_{\alpha}}$ which are unramified away from the divisors in $\mathbb M_{_{\alpha}}$. 

Since $\mathbb V$ is defined on $\mathbb M$, the quasi-admissible vector bundles obtained by restricting $\mathbb V$ agree on the intersection $\mathbb M_{_{\alpha \beta}} :=  \mathbb M_{_{\alpha}} \cap \mathbb M_{_{\beta}}$. Hence we have two covers $\tilde{\mathbb M}_{_{\alpha \beta}} \to \mathbb M_{_{\alpha \beta}}$ and $\tilde{\mathbb M}_{_{\beta \alpha}} \to \mathbb M_{_{\alpha \beta}}$ with same ramification data. By \cite[Proof of corollary 2.6, page 56]{vieweg}, it follows that we can go to a larger cover $\tilde{\tilde{\mathbb M}}$ which is \'etale over both.

For each $\alpha$, the bundle ${\mathbb V}_{_{\alpha}}$ on $\mathbb M_{_{\alpha}}$  comes as an invariant direct image $({\tt Inv} \circ \kappa_{_\alpha})_{_*}(V'_{_\alpha})$ for  equivariant vector bundles $V'_{_\alpha}$ on $\tilde{\mathbb M}_{_{\alpha}}$. Whence, by going to $\tilde{\tilde{\mathbb M}}$ dominating $\tilde{\mathbb M}_{_{\beta \alpha}}$ we can identify the pull-backs of $V'_{_{\alpha}}$ and $ V'_{_{\beta}}$.

Now each $V_{_{\alpha}}$ comes as $\ce^{^{par}}_{_{\alpha}}(k^{^n})$ for admissible pairs $(\qh_{_{\tau, W[d_{_\alpha}]}}^{^{\text{GL}(n)}},  \ce_{_{\alpha}})$ on $W[d_{_\alpha}]$. Furthermore, each $\ce_{_{\alpha}}$ comes as invariant push-forwards of equivariant torsors $\ce'_{_{\alpha}}$ on $\tilde{W}_{_{\alpha}}$. Let ${\sf E}'_{_{\alpha}} = \xi_{_{\alpha}}^{^*}(\ce'_{_{\alpha}})$ on $\tilde{\mathbb M}_{_{\alpha}}$. Then, we have ${\sf E}'_{_{\alpha}}(k^{^n}) = V'_{_{\alpha}}$. 

Thus, we can glue together the invariant push-forwards $ ({\tt Inv} \circ \kappa_{_\alpha})_{_*}(\sf E'_{_{\alpha}})$ to construct the universal group algebraic space ${\qh}^{^{\text{GL}(n)}}_{_{\tt univ}}$ together with a  torsor $\eP$ on $\mathbb M$ such that $({\qh}^{^{\text{GL}(n)}}_{_{\tt univ}}, \eP)$ gives the "universal" admissible pair. This clearly has the property that:
\beqa\label{universalgl} \eP^{^{par}}(k^{^n}) \simeq \mathbb V.
\eeqa
\epr
\subsubsection{Universal constructions for $G$}\label{meta}
Fix a faithful representation $\eta:G \hra \text{GL}(n)$. We have a smooth $k$-scheme $\mathbb M$ with a group algebraic space ${\qh}^{^{\text{GL}(n)}}_{_{\tt univ}}$. Further, the inclusion $\eta$ gives a constant subgroup scheme with fibre type $G$
away from the divisor $\mathbb M_{_{\mathtt Y_{_o}}}$. By taking flat closure at the generic point of the divisor $\mathbb M_{_{\mathtt Y_{_o}}}$,  this subgroup scheme extends to a subgroup algebraic space over an open subscheme $\mathbb M_{_\eta} \subset \mathbb M$ with complement of codimension $\geq 2$. We may assume that $\mathbb M_{_\eta}$ is the maximal such open subset to which the subgroup scheme extends. Let ${\qh}^{^{G}}_{_{\tt univ}}$ denote this subgroup algebraic space over $\mathbb M_{_\eta}$ and \beqa\label{univincl} \eta_{_{univ}}: {\qh}^{^{G}}_{_{\tt univ}}  \hra {\qh}^{^{\text{GL}(n)}}_{_{\tt univ}}\eeqa the canonical inclusion over $\mathbb M_{_\eta}$.
\brem We wish to emphasize that at them moment we only have the group algebraic space ${\qh}^{^{G}}_{_{\tt univ}}$ on $\mathbb M_{_\eta}$ and not the universal torsor. \erem   

\subsubsection{Admissible pairs on modifications} We make {\sl two definitions}, the first one \eqref{admistorsoronT1} makes use of a faithful representation and is for practical purposes and applications. The second one \eqref{admistorsoronT2} is independent of the choice of any representations. However this involves some other choices, but the definition will be shown to be independent of these choices. 

The two will be reconciled in \eqref{admississame}.
Let $\mathbb M \to \mathtt Y$ be as in \eqref{univconstructs} and \eqref{theunivconstructvb}.
\bdefe\label{admistorsoronT1} Let $\eta:G \hra \text{GL}(n)$ be a faithful representation. An {\em $\eta$-admissible 
pair} $({\qh}^{^{G}}_{_{{\mathfrak t},\tt M}}, \eE)$ on a modification $q:{\tt M} \to T$ consists of a diagram: 
\beqa\label{keydiagT1}
\xymatrix{
{\tt M} \ar[r]^{\phi_{_\mathfrak t}} \ar[d]_{q} &
 {\mathbb M}  \ar[d]_{} \\
T \ar[r]^{\psi_{_\mathfrak t}} &  {\tt Y} \\
}
\eeqa
where $\phi_{_\mathfrak t}$ factors as $\phi_{_\mathfrak t}:\tt M \to \mathbb M_{_\eta} \subset \mathbb M$ and such that the following hold: 
\begin{itemize}
\item The group algebraic space ${\qh}^{^{G}}_{_{{\mathfrak t},\tt M}}$ is obtained as a pull-back, i.e. ${\qh}^{^G}_{_{{\mathfrak t},\tt M}} := \phi_{_\mathfrak t}^{^{*}}(\qh^{^{G}}_{_{\tt univ}})$.
\end{itemize}
Define the group space ${\qh}^{^{\text{GL}(n)}}_{_{{\mathfrak t},\tt M}} := \phi_{_\mathfrak t}^{^{*}}(\qh^{^{\text{GL}(n)}}_{_{\tt univ}})$
and let $\eta: {\qh}^{^G}_{_{{\mathfrak t},\tt M}} \hra {\qh}^{^{\text{GL}(n)}}_{_{{\mathfrak t},\tt M}}$ be the  inclusion obtained by pulling back the inclusion \eqref{univincl}. 
\begin{itemize} 
\item The ${\qh}^{^G}_{_{{\mathfrak t},\tt M}}$-torsor $\eE$ on $\tt M$ is 
such that $({\qh}^{^{\text{GL}(n)}}_{_{{\mathfrak t},\tt M}}, \eta_{_{*}}(\eE))$ is isomorphic to the pull-back $\phi_{_\mathfrak t}^{^{*}}({\qh}^{^{\text{GL}(n)}}_{_{\tt univ}},\eP)$, where $\eta_{_{*}}(\eE)$ is the ${\qh}^{^{\text{GL}(n)}}_{_{{\mathfrak t},\tt M}}$-torsor obtained by extension of structure group via $\eta$.
\end{itemize}
\edefe
\bdefe\label{admistorsoronT2} Let ${\tt M} \to T$ be a modification as in \eqref{modification}. A pair $({\qh}^{^{G}}_{_{{\mathfrak t},\tt M}}, \eE)$ is called {\em admissible} if for each $\alpha$, there is an admissible pair $(\qh_{_{\tau, W[d_{_\alpha}]}}^{^G}, \ce_{_\alpha})$ on the standard model $W[d_{_\alpha}]$ such that $\xi_{_\alpha}^{^*}(\qh_{_{\tau, W[d_{_\alpha}]}}^{^G}, \ce_{_\alpha}) \simeq ({\qh}^{^{G}}_{_{{\mathfrak t},\tt M}}, \eE)\mid_{_{\tt M_{_{\alpha}}}}$. \edefe
The fact that the definition \eqref{admistorsoronT2} is   {\bf independent} of  the choice of the covering $\{T_{_\alpha} \}$, is derived as a consequence from \eqref{admississame} below. Recall the notion of an {\em effective degeneration} of $C_{_A}$ (\cite[Page 527]{Li}), i.e. it is a modification ${\tt M} \to T$ as in \eqref{modification} such that there is a morphism $\zeta:{T} \to B[d]$ for a {\sl single $d$}, and such that isomorphism $\tt M \simeq \xi^{^*}(W[d])$ is compatible with the morphism to $C_{_T}$.
\bprop\label{junli1} Let ${\tt M} \to T$ be a modification which is made effective by two morphisms $\zeta_{_i}:T \to B[d_{_i}], i = 1,2$ and arrows $\xi_{_i}$.  Suppose further that we have admissible pairs $\big({\qh_{_{\tau_{_i}, W[d_{_i}]}}^{^G}, \ce_{_i}}\big), i = 1,2$ on the standard models $W[d_{_i}], i = 1,2$, which are such that $\xi_{_1}^{^*}\big({\qh_{_{\tau_{_1}, W[d_{_1}]}}^{^G}, \ce_{_1}}\big) \simeq \xi_{_2}^{^*}\big({\qh_{_{\tau_{_2}, W[d_{_2}]}}^{^G}, \ce_{_2}}\big)$ on $\tt M$. Then, given a faithful representation $\eta:G \hra \text{GL}(n)$, there is a unique $A$-morphism $\psi:T \to {\tt Y}$ and a corresponding morphism $\phi:\tt M \to \mathbb M_{_\eta} \subset \mathbb M$ (with a diagram \eqref{keydiagT1}), such that the pair $({\qh}^{^{G}}_{_{{\mathfrak t},\tt M}}, \eE) ~~~\big[{:=  \xi_{_1}^{^*}\big({\qh_{_{\tau_{_1}, W[d_{_1}]}}^{^G}, \ce_{_1}}\big) \simeq \xi_{_2}^{^*}\big({\qh_{_{\tau_{_2}, W[d_{_2}]}}^{^G}, \ce_{_2}}}\big) \big]$ is $\eta$-admissible  \eqref{admistorsoronT1}. \eprop
\bpr By \eqref{etaonstd1}, the faithful representation $\eta$ gives rise to  quasi-admissible vector bundles $\xi_{_i}^{^*}\big({\ce_{_i}^{^{par}}(k^{^n})}\big)$, $i = 1,2$ on $\tt M$, which in turn give unique morphisms $\psi_{_i}:B[d_{_i}] \to {\tt Y}$ and $\phi_{_i}:W[d_{_i}] \to {\mathbb M}$ such that $\phi_{_i}^{^*}(\mathbb V) \simeq  \xi_{_i}^{^*}\big({\ce_{_i}^{^{par}}(k^{^n})}\big)$, $i=1,2$.

By \cite[Lemma 1.8]{Li}, for each $p \in {T}$ there is an \'etale neighbourhood ${T}_{_o}$  such that the isomorphism between $\big({\xi_{_1}^{^*}(W[d_{_1}])}\big)_{_\alpha} \simeq \big({\xi_{_2}^{^*}(W[d_{_2}])}\big)_{_\alpha}$ on ${T}_{_o}$ is induced by a sequence of effective arrows \cite[page 527]{Li}, and this holds for each $\alpha$. As in the proof of \cite[Lemma 1.8]{Li}, we may assume that on ${T}_{_o}$, we have two morphisms $\zeta_{_i}:{T}_{_o} \to B[d]$ such that $\zeta_{_i}(p) = 0 \in B[d]$, which induces the isomorphism via the $\xi_{_i}$'s.  

The assumption, $\zeta_{_1}(p) = \zeta_{_2}(p) = 0 \in B[d]$, forces the effective arrow inducing the isomorphism between $\xi_{_1}^{^*}(W[d])$ and $\xi_{_2}^{^*}(W[d])$ on ${ T}_{_o}$ to be the one that is induced by an automorphism of the fibre $W[d]_{_0}$ which commutes with the canonical projection $W[d]_{_0} \to C$. This is induced by an element of the group $G[d+1]$ \eqref{thegmaction0} by \cite[Corollary 1.4]{Li}. By \cite[Lemma 1.2]{Li} this therefore lifts to an action on $W[d]$ preserving the configuration of the s.n.c divisor. 

The isomorphism $\xi_{_1}^{^*}\big({\qh_{_{\tau_{_1}, W[d]}}^{^G}, \ce_{_1}}\big) \simeq \xi_{_2}^{^*}\big({\qh_{_{\tau_{_2}, W[d]}}^{^G}, \ce_{_2}}\big)$ shows that restricted to the single semistable curve $C^{^{(d+1)}} \simeq W[d]_{_0}$, the admissible pairs are isomorphic. Recall that the Kawamata coverings $\tW_{_{\tau_{_j}}}[d] \to W[d]$ (with Galois group $\mathfrak G$) induced by the pairs $(\qh_{_{\tau_{_i}, W[d]}}^{^G}, \ce_{_i}), i = 1,2$ are completely determined by the admissible pair at the central fibre, i.e. the local types are completely determined, which in turn determines the ramification data. Thus, we conclude that the Kawamata covers are isomorphic. By pulling back using $\xi_{_i}$, we get a finite flat cover  $\tilde{\tt M}_{_o} \to {\tt M}_{_o}$ together with an action of  $\mathfrak G$. This allows us to take Weil restrictions and invariants.  

The isomorphism of the pairs $\xi_{_1}^{^*}\big({\qh_{_{\tau_{_1}, W[d_{_1}]}}^{^G}, \ce_{_1}}\big)$ and $\xi_{_2}^{^*}\big({\qh_{_{\tau_{_2}, W[d_{_2}]}}^{^G}, \ce_{_2}}\big)$ when restricted to ${\tt M}_{_o}$ is recovered as invariant push-forwards of isomorphism of pairs on $\tilde{\tt M}_{_o}$. By taking associated vector bundles and invariant push-forwards, we can firstly realize the admissible vector bundles $\xi_{_i}^{^*}({\ce_{_i}^{^{par}}(k^{^n})})$ on ${\tt M}_{_o}$ as invariant push-forwards from $\tilde{\tt M}_{_o}$.  Then the isomorphism of the pairs induces an isomorphism of these associated parabolic bundles $\xi_{_i}^{^*}({\ce_{_i}^{^{par}}(k^{^n})})$ over ${\tt M}_{_o}$. These are isomorphic quasi-admissible vector bundles on the modification ${\tt M}_{_o}$. Thus, they induce the same morphisms from $T_{_o}$ to the universal space $\tt Y$ (by uniqueness). All in all, we can conclude that on ${\tt M}_{_o} \to T_{_o}$, the restrictions of the morphisms $\psi_{_i}$ and $\phi_{_i}$ coincide  and hence they coincide everywhere on $T$. The remaining statements are straightforward.\epr

\bth\label{admississame} ("Tannakian") Let $\tt M \to T$ be a modification as in \eqref{modification}. A pair $({\qh}^{^{G}}_{_{{\mathfrak t},\tt M}}, \eE)$ is admissible  \eqref{admistorsoronT2} if and only if it is $\eta$-admissible \eqref{admistorsoronT1} for any faithful $\eta$. In particular, the definition \eqref{admistorsoronT2} does not depend on the covering $\{T_{_\alpha}\}$.\eeth
\bpr \eqref{admistorsoronT1} implies \eqref{admistorsoronT2} is easily seen. For the other direction, on each $T_{_\alpha}$ we have unique  morphisms $\psi_{_\alpha}, \phi_{_\alpha}$ with the added properties. The intersection $T_{_{\alpha,\beta}} := T_{_\alpha} \cap T_{_\beta}$ satisfies the assumption in \eqref{junli1} and hence the uniqueness forces that  $\psi_{_\alpha}, \phi_{_\alpha}$ agree on the intersections and hence glue to give unique morphisms $\psi,\phi$ which does the job.\epr

\section{The stack of Gieseker torsors}
Using the admissible pairs constructed on modifications  we define the stacks of Gieseker torsors and study the basic properties.

\bdefe\label{giesekertorsor}
For a scheme $T$ over $A$, a {\sf Gieseker torsor on $C_{_T}$ over $T$}  is a datum $\Big(\tt M,{\qh}^{^G}_{_{{\mathfrak t},\tt M}}, \ce\Big)$, consisting of a modification ${\tt M} \to T$,       and an  {\em admissible pair} $({\qh}^{^G}_{_{{\mathfrak t},\tt M}}, \ce)$  on ${\tt M}$ \eqref{admistorsoronT2}. 
\edefe

 Two Gieseker torsors $({\tt M}_{_1}, {\qh}^{^G}_{_{{\mathfrak t}_{_1},{\tt M}_{_1}}},\ce_{_1})$ and $({\tt M}_{_2}, {\qh}^{^G}_{_{{\mathfrak t}_{_2},{\tt M}_{_2},}}, {\ce}_{_2})$ on $C_{_{T_{_j}}}, j = 1,2$  are called {\em isomorphic} if there exists a $A$-isomorphism $\delta:T_{_1} \to T_{_2}$ and a diagram:
\beqa\label{automor1}
\xymatrix{
{\tt M}_{_1} \ar[r]^{\epsilon} \ar[d] &
 {\tt M}_{_2} \ar[d] \\
T_{_1} \ar[r]^{\delta} &  T_{_2}
}
\eeqa
compatible with the tautological projections $\tt M_{_j} \to C_{_{T_{_j}}}$, $j = 1,2$,
an isomorphism 
\beqa
\epsilon^{^*}({\qh}^{^G}_{_{{\mathfrak t}_{_2},{\tt M}_{_2}}}, {\ce}_{_2})  \simeq ({\qh}^{^G}_{_{{\mathfrak t}_{_1},{\tt M}_{_1}}},\ce_{_1}).
\eeqa
of admissible pairs on ${\tt M}_{_1}$. 

\bdefe
Let ${\text{Gies}}_{_G}(C_{_A})$   be the category over $\text{Sch}\big/_{{A}}$, 
whose objects  are Gieseker torsors $\Big(\tt M,{\qh}^{^G}_{_{{\mathfrak t},\tt M}}, \ce\Big)$.

The functor $\mathfrak f:{\text{Gies}}_{_G}(C_{_A}) \to \text{Sch}\big/_{{A}}$ which sends $\Big(\tt M,{\qh}^{^G}_{_{{\mathfrak t},\tt M}}, \ce\Big) \mapsto T$ realizes it as a fibered category. 
\edefe

An {\em arrow} between two objects $ \Upsilon_{_1} = \Big(\tt M,{\qh}^{^G}_{_{{\mathfrak t},\tt M}}, \ce\Big)$ and $\Upsilon_{_2} = \Big(\tt M',{\qh}^{^G}_{_{{\mathfrak t'},\tt M'}}, \ce'\Big)$ over $T$ and $T'$ consists of (1) a $A$-morphism $T \to T'$, (2) an isomorphism of modifications $\tt M \to \tt M' \times_{_T'} T$ and (3) an isomorphism over $T$ of Gieseker torsors $(\tt M, \ce)$ and $(\tt M' \times_{_T'} T, \ce' \times_{_{T'}} T)$.

 Since pull-backs of modifications (resp. group algebraic spaces, torsors) are modifications (resp. group algebraic spaces, torsors), and arrows between two objects are as defined above and are fiber diagrams, the category ${\text{Gies}}_{_G}(C_{_A})$ is fibered in groupoids under $\mathfrak f$. We have the following straightforward result.
\bprop\label{stackaxioms} The category  ${\text{Gies}}_{_G}(C_{_A})$   is a stack. \eprop
\bpr It suffices to show the following:
\begin{enumerate}
\item For any $T \in \text{Sch}/_{{A}}$ and two objects $\Upsilon_{_1},\Upsilon_{_2} \in {\text{Gies}}_{_G}(C_{_A})(T)$, the functor:
\beqa
\text{Isom}_{_T}(\Upsilon_{_1},\Upsilon_{_2}):\text{Sch}_{_T} \to \text{Sets}
\eeqa
which associates to any morphism $\phi:T' \to T$ the set of isomorphisms in ${\text{Gies}}_{_G}(C_{_A})$ between $\phi^{^*}(\Upsilon_{_1})$ and $\phi^{^*}(\Upsilon_{_2})$, is a sheaf in the \'etale topology.
\item ({\em Effective descent}) Let $\{T_{_i} \to T\}$ be a covering of $T$ in the \'etale topology. Let $\Upsilon_{_i}  \in {\text{Gies}}_{_G}(C_{_A})(T_{_i})$ and let $\phi_{_{ij}}:\Upsilon_{_i}\vert{_{_{T_{_i} \times_{_T} T_{_j}}}} \to \Upsilon_{_j}\vert{_{_{T_{_i} \times_{_T} T_{_j}}}}$ be isomorphisms in ${\text{Gies}}_{_G}(C_{_A})(T_{_i} \times_{_T} T_{_j})$ satisfying the cocycle condition. Then there is an $\Upsilon \in {\text{Gies}}_{_G}(C_{_A})(T)$ with isomorphisms $\psi_{_i}:\Upsilon\vert_{_{T_{_i}}} \to \Upsilon_{_i}$ so that 
\beqa
\phi_{_{ij}} = \psi_{_i} \circ \psi_{_j}^{^{-1}}.
\eeqa

\end{enumerate}
Each object $\Upsilon = \Big(\tt M,{\qh}^{^G}_{_{{\mathfrak t},\tt M}}, \ce\Big)$ consists of three components.  Modifications or expanded degenerations form a stack \cite[1.10, page 531]{Li} or \cite[Proposition 3.16]{kausz}.

For the second item, this follows  since the sheaf property and effective descent is automatic for morphisms.

The sheaf property of the third component, namely the torsor $\ce$ is immediate since isomorphisms of ${\qh}^{^G}_{_{{\mathfrak t},\tt M}}$-torsors $\ce_{_1}$ and $\ce_{_2}$ is given by a section of $(\ce_{_1} \times \ce_{_2})\Big(\tfrac{{\qh}^{^G} \times {\qh}^{^G}} {\Delta}\Big)$. Effective descent of torsors holds in the category of algebraic spaces, and   the action maps descend by \cite[Theorem 6, Section 6.1]{blr}. The admissibility of the descended pair is immediate since it holds on each  $T_{_i}$.\epr

\bth\label{balequalskausz} When $G = \text{\text{GL}(n)}$, we have  isomorphisms:

\beqa\label{oldisom} 
{\text{Gies}}_{_{\text{GL}(n)}}(C_{_A})
\simeq \text{GVB}_{_n}(C_{_{A}}).
\eeqa
In particular, ${\text{Gies}}_{_{\text{GL}(n)}}(C_{_A})$ is an algebraic ${A}$-stack, locally of finite type.
\eeth
\bpr We show that there is a canonical functor $\mathfrak e:{\text{Gies}}_{_{\text{GL}(n)}}(C_{_A})  \to \text{GVB}_{_n}(C_{_{A}})$ defined over $\text{Sch}\big/_{A}$, which is an equivalence of fibered categories. Let $T \in \text{Sch}\big/_{A}$. Let $\Big({\tt M},{\qh}^{^{\text{\text{GL}(n)}}}_{_{{\mathfrak t},\tt M}}, \ce\Big) \in {\text{Gies}}_{_{\text{GL}(n)}}(C_{_A})(T)$. By the definition of an admissible pair \eqref{admistorsoronT1}, we have the {parabolically}-associated vector bundle $\ce^{^{par}}(k^{^n}) := \phi_{_{\mathfrak t}}^{^*}(\eP^{^{par}}(k^{^n})) = \phi_{_{\mathfrak t}}^{^*}(\sf V)$  as in \eqref{universalgl}, and clearly $\ce^{^{par}}(k^{^n})$ is  a quasi-admissible vector bundle. Define:
\beqa
\mathfrak e\Big(\tt M,{\qh}^{^{\text{\text{GL}(n)}}}_{_{{\mathfrak t},\tt M}}, \ce\Big) := \Big(\tt M, \ce^{^{par}}(k^{^n})\Big)
\eeqa

The functor $\mathfrak e$ is essentially surjective for each $T \in \text{Sch}\big/_{A}$. Let $\eV$ be an admissible bundle on $\tt M$ of rank $n$. The universal property of $(\mathbb M, \mathbb V)$ shows that there is a morphism $\phi:\tt M \to \mathbb M$ such that $\phi^{^*}(\mathbb V) = \eV$. We get back $({\qh}^{^{\text{\text{GL}(n)}}}_{_{{\mathfrak t},\tt M}}, \ce)$ by pulling back $({\qh}^{^{\text{GL}(n)}}_{_{\tt univ}},\eP)$. That $\mathfrak e$ is fully-faithful is immediate from the definitions of isomorphisms of the objects.\epr

\bth\label{algebraicstack} The stack ${\text{Gies}}_{_G}(C_{_A})$   is an algebraic ${A}$-stack, locally of finite type. 
For the fixed nodal curve $(C,c)$ over $\bc$, ${\text{Gies}}_{_{G}}(C)$ is an algebraic $\bc$-stack, locally of finite type.
\eeth
\bpr  By  the definition of Gieseker torsors we have a  morphism: 
\beqa\label{heeheeeee}
\eta_{_*}: {\text{Gies}}_{_G}(C_{_A}) \to {\text{Gies}}_{_{\text{GL}(n)}}(C_{_A})\\
\Big(\tt M,{\qh}^{^G}_{_{{\mathfrak t},\tt M}}, \ce\Big) \mapsto \Big(\tt M,{\qh}^{^{\text{GL}(n)}}_{_{{\mathfrak t},\tt M}}, \eta_{_*}(\ce)\Big)
\eeqa
where $\eta: {\qh}^{^{G}}_{_{{\mathfrak t},\tt M}} \hra {\qh}^{^{\text{\text{GL}(n)}}}_{_{{\mathfrak t},\tt M}}$ is the canonical inclusion.

\noindent
{\em This morphism of stacks is representable, locally of finite presentation.} To see this we follow \cite{behrend}; let $T$ be a ${A}$-scheme and let $P$ be an admissible  ${\qh}^{^{\text{\text{GL}(n)}}}_{_{{\mathfrak t},\tt M}}$-torsor on a modification ${\tt M} \to T$ \eqref{admistorsoronT1}. 

The quotient ${\qh}^{^{\text{\text{GL}(n)}}}_{_{{\mathfrak t},\tt M}} \Big/ {\qh}^{^{G}}_{_{{\mathfrak t},\tt M}}$ exists as an algebraic space with a ${\qh}^{^{\text{\text{GL}(n)}}}_{_{{\mathfrak t},\tt M}}$ action. We identify the associated space $P\big({\qh}^{^{\text{\text{GL}(n)}}}_{_{{\mathfrak t},\tt M}} \Big/ {\qh}^{^{G}}_{_{{\mathfrak t},\tt M}}\big)$ with $P \big/ {\qh}^{^{G}}_{_{{\mathfrak t},\tt M}}$. Let $q:{\tt M} \to T$ be the arrow defining the modification. Then we have a $2$-cartesian diagram of ${A}$-stacks:
\beqa\label{2cartesian}
\xymatrix{
q_{_{*}}\big(P \big/ {\qh}^{^{G}}_{_{{\mathfrak t},\tt M}}\big) \ar[r]\ar[d] &
T \ar[d] \\
{\text{Gies}}_{_G}(C_{_A}) \ar[r] &  {\text{Gies}}_{_{\text{GL}(n)}}(C_{_A})
}
\eeqa
By  \cite[Corollary 2.17]{rydh}, $P \big/ {\qh}^{^{G}}_{_{{\mathfrak t},\tt M}}$ is an algebraic space of finite presentation. Thus, by using the theory of Hilbert schemes for algebraic spaces as in \cite[Section 6]{artin}, we see that $q_{_{*}}\big(P \big/ {\qh}^{^{G}}_{_{{\mathfrak t},\tt M}}\big)$ is also an algebraic $T$-space of finite presentation. Hence the morphism \eqref{heeheeeee} is locally of finite presentation.

By  (\ref{balequalskausz}) the stack ${\text{Gies}}_{_{\text{GL}(n)}}(C_{_A})$ is an algebraic ${A}$-stack locally of finite type. Hence by \eqref{2cartesian}, we conclude that ${\text{Gies}}_{_G}(C_{_A})$ is an algebraic stack locally of finite type over ${A}$. It is immediate that ${\text{Gies}}_{_{G}}(C)$ is an algebraic $k$-stack  and also locally of finite type being the closed fibre of ${\text{Gies}}_{_G}(C_{_A})$.
\epr	
\subsubsection{Deformations of Gieseker torsors}\label{deforgieseker}
We follow \cite{gies} and \cite[Appendix]{ns2}. We work in the setting of the diagram \eqref{keydiagT1} with $\eta:G \hra \text{GL}(n)$. Let  $\mathbb M \to {\tt Y}$ be the universal modification \eqref{univconstructs} and $\mathbb M_{_\eta}$ be as in \eqref{meta}. Let $\psi_{_\mathfrak t}:T \to {\tt Y}$ be a ${\tt Y}$-scheme and $q:{\tt M} \to T$ a modification such that the morphism ${\tt M} \to\mathbb M$ factors as   ${\tt M} \to {\mathbb M_{_\eta}} \subset \mathbb M$ as in \eqref{keydiagT1}. Let 
$\phi_{_\mathfrak t}:{\tt M} \to {\mathbb M_{_\eta}}$ be the induced morphism.

 Let $P_{_{\tt M}} = \phi_{_\mathfrak t}^{^*}(\eP)$ which is therefore a ${\qh}^{^{\text{\text{GL}(n)}}}_{_{{\mathfrak t},\tt M}}$-torsor on $\tt M$. Further, we have the group algebraic space ${\qh}^{^{G}}_{_{{\mathfrak t},\tt M}}$ and an inclusion of group algebraic spaces ${\qh}^{^{G}}_{_{{\mathfrak t},\tt M}} \hra {\qh}^{^{\text{\text{GL}(n)}}}_{_{{\mathfrak t},\tt M}}$ over $\tt M$.

\bdefe\label{yh} Define the functor ${\eG}^{^{G}}_{_{{A}}}:\text{Sch}_{_{{\mathtt Y}}} \to \text{Sets}$
\beqa
T \mapsto \left \{ (P_{_{\tt M}},\zeta) \Big\vert \begin{array}{l}
\mbox{ $\zeta \in \Gamma
({\tt M},\Big(P_{_{\tt M}}\Big({\qh}^{^{\text{\text{GL}(n)}}}_{_{{\mathfrak t},\tt M}} \Big/ {\qh}^{^{G}}_{_{{\mathfrak t},\tt M}}\Big)\Big)$}\\
\end{array} \right \} 
\eeqa
i.e. ${\eG}^{^{G}}_{_{{A}}} (T) $ consists of isomorphism classes of pairs $(P_{_{\tt M}},\zeta)$ on the modification ${\tt M} \to T$, where $\zeta$ is a  
reduction of structure group of $P_{_{\tt M}}$  to ${\qh}^{^{G}}_{_{{\mathfrak t},\tt M}}$.\edefe
We show that this functor is representable by a $\tt Y$-scheme, following \cite[Page 424-425]{ram2}, i.e. by embedding the homogeneous space ${\qh}^{^{\text{\text{GL}(n)}}}_{_{{\mathfrak t},\tt M}} \Big/ {\qh}^{^{G}}_{_{{\mathfrak t},\tt M}}$ in a vector bundle over $\tt M$. By Chevalley's theorem on semi-invariants we obtain an embedding:
\beqa
\text{GL}(n)/G \hra W
\eeqa 
in a $\text{GL}(n)$-module $W$. 

By following \eqref{univconstructs} and using the $\text{GL}(n)$-module $W$, it is straightforward to see that we can take the associated vector bundles ${\sf E}'_{_{\alpha}}(W)$ on $\tilde{\mathbb M}_{_{\alpha}}$ with its natural equivariant structure and define the invariant push-forward $\eW_{_{\alpha}} := ({\tt Inv} \circ \kappa_{_*})({\sf E}'_{_{\alpha}}(W))$. These glue up to give a vector bundle $\cw := \eP^{^{par}}(W)$ on $\mathbb M$. Since $\eP$ is a ${\qh}^{^{\text{GL}(n)}}_{_{univ}}$-torsor,  by restricting to $\mathbb M_{_\eta}$ \eqref{meta}, we can consider the algebraic space $\eP\big({\qh}^{^{\text{GL}(n)}}_{_{univ}}/{\qh}^{^{G}}_{_{univ}}\big)$ over $\mathbb M_{_\eta}$.

Restricting $\cw$ to $\mathbb M_{_\eta} \subset \mathbb M$, we get an embedding :
\beqa
\eP\big({\qh^{^{\text{GL(n)}}}_{_{\tt univ,\eta}}/\qh^{^{G}}_{_{\tt univ}}}\big) \hra \cw
\eeqa

Pulling back using the morphism $\phi_{_\mathfrak t}: \mathtt M \to \mathbb M_{_\eta}$, we see that we have the embedding:
\beqa
P_{_{\tt M}}\big({{\qh}^{^{\text{\text{GL(n)}}}}_{_{{\mathfrak t},\tt M}}/{\qh}^{^{G}}_{_{{\mathfrak t},\tt M}}}\big) \hra \cw_{_{\tt M}} 
\eeqa
where $\cw_{_{\tt M}} := \phi_{_\mathfrak t}^{^*}(\cw)$

In other words, we can realize the functor ${\eG}^{^{G}}_{_{{A}}}$ as a closed subfunctor of the functor $T \mapsto H^{^0}\Big(\cw_{_{\tt M}}\Big)$. By \cite[Proposition 8]{ns2} and \eqref{regularityofquasistuff}, $\mathtt Y$ is a reduced scheme and hence the functor $T \mapsto H^{^0}\Big(\cw_{_{\tt M}}\Big)$  is representable by a linear scheme; therefore, there exists a $\mathtt Y$-scheme ${\mathtt Y}^{^G}$  which represents ${\eG}^{^{G}}_{_{{A}}}$. 

We can also describe the $T$-points of ${\eG}^{^{G}}_{_{{A}}}(T)$
as $\big[(\tt M,\text{\cursive e}, {\qh}^{^{G}}_{_{{\mathfrak t},\tt M}},\ce)\big]$, where $(\tt M,\text{\cursive e})$ are as in \eqref{cursivee}, and $\ce$ is a ${\qh}^{^{G}}_{_{{\mathfrak t},\tt M}}$-torsor. Equivalently, one could describe it as $\big[(\tt M,\text{\cursive e}, \eta_{_*}(\ce),\zeta_{_T})\big]$, where $\big({{\qh}^{^{\text{GL}(n)}}_{_{{\mathfrak t},\tt M}}, \eta_{_*}(\ce)}\big)$ is an admissible pair and $\zeta_{_T}$ is a reduction of structure group to the subgroup scheme ${\qh}^{^{G}}_{_{{\mathfrak t},\tt M}}$. 

Notice that $\ce^{^{par}}(k^{^n})$ is a quasi-admissible vector bundle. Let
\beqa
\eV_{_{T}} := \ce^{^{par}}(k^{^n})
\eeqa
be the vector bundle on $\tt M$ which is {parabolically}-associated to the torsor $\eta_{_*}(\ce)$. Thus, giving the  representation $\eta:G \hra \text{GL}(n)$ also induces a morphism ${\eG}^{^G}_{_{{A}}} \to {\eG}_{_{N}}^{^{\tt q}}$ (more precisely, by taking the associated vector bundle plus a twisting of the vector bundles by a positive $m$ to ensure that the first cohomology vanishes and the sections generate the bundle, see \cite[Page 176, Remark 4 (ii)]{ns2}).

Let ${\eG}'_{_{{A}}}$ be the functor defined as: 
\beqa\label{as in GiesekerfunctorreltoR1}
{\eG}'_{_{{A}}}:\text{Sch}_{_{\mathtt Y}} \to \text{Sets}\\
T \mapsto \tt M
\eeqa
such that $\tt M \to C \times_{_{A}} T$ is a modification.

Thus, we have a morphism from ${\eG}_{_{N}}^{^{\tt q}}$ to  ${\eG}'_{_{{A}}}$ obtained by forgetting the condition (1) in Definition \ref{GiesekerfunctorreltoR} namely, the imbeddings into the Grassmannians (see \cite[Appendix, page 197]{ns2}).  Composing with the morphism ${\eG}^{^G}_{_{{A}}} \to {\eG}_{_{N}}^{^{\tt q}}$  we have the induced forget morphism:
\beqa\label{formalsmoothmap}
{\eG}^{^G}_{_{{A}}} \to {\eG}'_{_{{A}}}\\
\big[(\tt M,\text{\cursive e}, {\qh}^{^{G}}_{_{{\mathfrak t},\tt M}},\ce)\big] \mapsto \tt M 
\eeqa

The functors ${\eG}^{^G}_{_{{A}}}$, ${\eG}'_{_{{A}}}$ are defined with a fixed choice of the fibered surface $C_{_{A}} \to {A}$. Further, the  functor ${\eG}'_{_{{A}}}$ defined above parametrizes semistable curves with a fixed stable model, being the irreducible nodal curve $(C,c)$ with a single node. 
Gieseker \cite[page 183]{gies}(see also \cite[Appendix]{ns2}) shows that  the canonical map $B[d] \to {\eG}'_{_{{A}}}$ defined by the point $W[d] \in {\eG}'_{_{{A}}}(B[d])$ is formally smooth.

\bth\label{flatavecnc} The algebraic stack ${\text{Gies}}_{_G}(C_{_A})$ is regular and flat over ${A}$; further, ${\text{Gies}}_{_{G}}(C) \subset {\text{Gies}}_{_G}(C_{_A})$ is a divisor with normal crossings. More precisely, the morphism  \eqref{formalsmoothmap} is  formally smooth.  \eeth

\bpr Let $T$ be the spectrum of an Artin local ring, and $T_{_o} \subset T$ the subscheme defined by an ideal of  dimension $1$. Let $\tt M \in {\eG}'_{_{{A}}}(T)$ be such that the restriction $\tt M_{_o} \in {\eG}'_{_{{A}}}(T_{_o})$ can be lifted to an element of  
${\eG}^{^{G}}_{_{{A}}}(T_{_o})$, then we need to show that $\tt M$ itself can be lifted to an element of ${\eG}^{^{G}}_{_{{A}}}(T)$. Let $\tt M$ be defined by the family of curves $\tt M \to T$, and by the modification $\tt M \to C \times_{_{A}} T$. 

The lifting of the family $\tt M_{_o}\to T_{_o}$ to $\tt M \to T$ comes with information which we require: there is a morphism $\phi_{_o}:T_{_o} \to B[d]$ such that pull-backs  by $\phi_{_o}$ of the versal families $W[d] \to B[d]$  and $W[d] \to C_{_{B[d]}}$ coincide with the datum given by the point $\tt M_{_o} \in {\eG}'_{_{{A}}}(T_{_o})$. Gieseker then shows that we have a diagram:
\beqa\label{giesekersdiagram}
\xymatrix{
T_{_o} \ar[r]^{^\subset} \ar[dr]_{\phi_{_o}} &
 T \ar[d]^{\phi} \\
 &  B[d]
}
\eeqa 
such that the pull-backs of $W[d] \to B[d]$  and $W[d] \to C_{_{B[d]}}$ give the family $\tt M \to T$ and the point $\tt M \in {\eG}'_{_{{A}}}(T)$.

The lifting of $\tt M_{_o}$ to an element of  
${\eG}^{^{G}}_{_{{A}}}(T_{_o})$ defines an admissible pair $({\qh^{^G}_{_{\mathfrak t,\tt M_{_o}}}}, \ce_{_{T_{_o}}})$  on the restriction $\tt \tt M_{_o}$ of $\tt M$ to $T_{_o}$. The problem is:
\begin{enumerate}
\item to extend the pair $({\qh^{^G}_{_{\mathfrak t,\tt M_{_o}}}}, \ce_{_{T_{_o}}})$ to an admissible pair $({\qh^{^G}_{_{\mathfrak t,\tt M}}}, \ce_{_{T}})$ on $\tt M$  
\item to lift the morphism $\tt \tt M_{_o} \to T_{_o} \times  W(N,n)$ to a morphism $\tt M \to T \times W(N,n)$. 
\end{enumerate}

(1): By versality, the  group algebraic space ${\qh^{^G}_{_{\mathfrak t,\tt M_{_o}}}}$ is isomorphic to the  pull-back of a group algebraic space ${\qh^{^G}_{_{\tau, W[d]}}}$ on $W[d]$ by the morphism $\tt M_{_o} \to W[d]$ induced by $\phi_{_o}$. The diagram \eqref{giesekersdiagram} thus gives {\em the group algebraic space ${\qh^{^G}_{_{\mathfrak t,\tt M}}}$ on $\tt M_{_{T}}$ extending ${\qh^{^G}_{_{\mathfrak t,\tt M_{_o}}}}$},  namely the pull-back of ${\qh^{^G}_{_{\tau, W[d]}}}$ by the morphism $\tt M_{_{T}} \to W[d]$ induced by $\phi$.

Let $E$ be the restriction of  $\ce_{_{T_{_o}}}$ to the closed fibre $C^{^{(d)}}$ of $\tt M_{_{T_{_o}}}\to T_{_o}$. Let ${\qh}^{^G}_{_{\tau}}$ be the restriction of ${\qh^{^G}_{_{\mathfrak t,\tt M}}}$ to $C^{^{(d)}}$. It is standard that the obstruction to lifting $\ce_{_{T_{_o}}}$ to $ \ce_{_{T}}$  is simply the group $\text{H}^{^2}\Big(C^{^{(d)}}, E\Big({\text{Lie}}({\qh}^{^G}_{_{\tau}} )\Big)\Big)$ which vanishes since we are in  dimension $1$. 

(2): For proving the second item, by the definition of $W(N,n)$ \eqref{zN} and  by what has been discussed above regarding the versal property, it remains to extend  the sections of the  vector bundle $\eV_{_{T_{_o}}}(m)$ ($m \gg 0$) which defines the given map $\tt \tt M_{_o} \to T_{_o} \times \text{Grass}(N,n)$ to the sections of $\eV_{_T}(m)$ so as to define the lift $\tt M \to T \times \text{Grass}(N,n)$. The second item is therefore possible since the obstruction to lifting of sections lies in $H^{^1}(\tt M_{_t}, V_{_t}(m))$ and this group vanishes by \eqref{h1-vanish}). Thus we conclude that the morphism \eqref{formalsmoothmap} is formally smooth. It is shown in \cite{gies} that there is a formally smooth morphism from the versal space $B[d]$ to the functor ${\eG}'_{_{{A}}}$.

We deduce (using \ref{expandeddeg}) that the scheme ${\mathcal Y}^{^G}$ which represents ${\eG}^{^G}_{_{{A}}}$ has all the stated properties. One then concludes (following \cite[Proposition 3.24]{kausz}) that the stack ${\text{Gies}}_{_G}(C_{_A})$ has all the stated properties.\epr
\subsubsection{Some remarks on a weak properness}\label{somebounds}

Let $C_{_{d,A}} := C_{_A} \times_{_{\spec k\llbracket t \rrbracket}} \spec k\llbracket t \rrbracket$ via the map $t \mapsto t^{^d}$ \eqref{normalplusmodif}, and let $E_{_L}$ be a $G$-torsor on the generic fibre $C_{_{d,L}}$ of $C_{_{d,A}}$. It is not hard to see that the $G$-bundle extends to $C_{_{d,A}} - c$. Locally, we have a $G$-bundle on $N_{_d} - c$. By going to $N_{_0}$ ($N_{_0} = \spec~\frac{k\llbracket u,v \rrbracket}{(u.v)}$ with coordinates $u,v$), and using a Hartogs like argument, we can extend the $G$-torsor to a $(\smu,G)$-torsor of local type $\tau$ or equivalently a $\ce_{_D}(G)$-torsor \eqref{basictorsor}. Then it is straightforward to obtain a group algebraic space  ${\qh^{^G}_{_{\tau, S^{^{(d)}}}}}$ on the minimal desingularization $S^{^{(d)}} \to C_{_{d,A}}$ and an admissible pair $({\qh^{^G}_{_{\tau, S^{^{(d)}}}}}, \ce)$ such that over the generic fibre we have an isomorphism $\ce_{_L} \simeq E_{_L}$. As we have seen earlier, the group space  $({\qh^{^G}_{_{\tau, S^{^{(d)}}}}}, \ce)$ is non-trivial parahoric at the generic points of at most $\ell = \text{rank}(G)$ number of projective lines in the exceptional divisor. The closed fibre of the moduli stack is describable in terms of {\tt laced} parahoric torsors on $\tC$ with parahoric structures at two points \eqref{2-lacedbundle}. To work in families, we fix a faithful representation $\eta:G \hra \text{GL}(V)$, where ${\tt dim}(V) = n$ then the associated {\em parabolic bundle}, $\ce^{^{par}}(V)$ on $S^{^{(d)}}$ which is quasi-admissible is non-trivial in at most ${\tt dim}(V)$ rational curves on the exceptional divisor. Thus, the number ${\mathtt m}(G)$, being  the minimal dimension of a faithful representation of $G$, gives an upper bound for $d$. More precisely, the stack of Gieseker torsors can be obtained by taking torsors on semistable curves $C^{^{(d)}}$ with $d$ bounded above by ${\mathtt m}(G)$.

\section{The closed fibre of the stack and examples} 
By Bruhat-Tits theory, for each facet $\ul{a}$ in each apartment $\ca_{_T}$ of the Bruhat-Tits building of $G(k(\!(t)\!))$, there is a smooth group scheme $\sf {P}_{_{\ul{a}}}$ over $k\llbracket t \rrbracket$ with connected fibers whose generic fiber is $G \times _{_{\spec k}} \spec k (\!(t)\!)$. We call such  $\sf{P}_{_{\ul{a}}}$ a Bruhat-Tits group scheme. Let ${\bf P}_{_{\ul{a}}}(K) := {\sf P}_{_{\ul{a}}}(k\llbracket t \rrbracket)$. We call ${\bf P}_{_{\ul{a}}}$ a parahoric subgroup of $G(k(\!(t)\!))$. The conjugacy classes of parahoric subgroups of $G(k(\!(t)\!))$ are classified by proper subsets of the nodes of the extended Dynkin diagram of $G$ or the facets of the Weyl alcove ${\bf a} \subset \ca_{_T}$. Let $\sf {P}_{_{{\ul v}}}$ denote the group schemes associated to the maximal parahoric subgroups which are indexed by the vertices $\ul{\sf v}$ of ${\bf a}$. We summarize two results the details of which is work in progress.

\bth\label{components} 
\begin{enumerate}
\item The closed fibre ${\text{Gies}}_{_{G}}(C)$ of ${\text{Gies}}_{_{G}}(C_{_A})$ is a divisor with  simple normal crossing singuarities. It has $\ell + 1$ irreducible smooth components $\mathfrak G_{_j}$  indexed by the vertices of the extended Dynkin diagram. 
\item Let  $\sf P(\sf v)$ be the parahoric group scheme on the smooth projective curve $\tC$ which restricts to the maximal parahoric $\sf{P}_{_{\ul{{\sf v}}}}$ at the two marked points.
In each component $\mathfrak G_{_j}$, the open locus of Gieseker torsors can be identified with $\sf P(\sf v)$-torsors on $\tC$  for varying vertices ${\ul {\sf v}}$. 
\item The minimal stratum will be  torsors under ${\sf P}_{_I}$ on $\tC$  defined by the Iwahori group scheme at the two marked points on $\tC$.

\end{enumerate} 
\eeth
\bth\label{components1}  
\begin{enumerate}
\item Let ${\sf Q}_{_{\ul{{\sf v}}}}$ be the smooth group scheme (maximal parahoric) on the nodal curve  $(C,c)$ obtained by identifying the closed fibres of the group scheme $\sf P(\sf v)$ on $\tC$.
Then the stack $\text{Bun}^{^{\sf L}}_{_{C}}({\sf Q}_{_{\ul{{\sf v}}}})$ of laced $\sf P(\sf v)$ torsors \eqref{2-lacedbundle} \eqref{bttorsorlacedtorsor} on $\tC$ is a regular Artin stack. 
\item This stack is isomorphic to the stack of admissible torsors on $C^{^{(d)}}$ coming from representations $\rho:{\sf\mu_{_d}} \to G$ which give the maximal parahoric. 
\item The dimension of the stacks $\text{Bun}^{^{\sf L}}_{_{C}}({\sf Q}_{_{\ul{{\sf v}}}})$ and $\text{Bun}_{_{C}}(G)$ coincide and hence these constitute the components of the normal crossing divisor.
\end{enumerate}
\eeth

Let $e(\theta_{_\alpha})$ be as in \cite[7.2.1]{base} where $\theta_{_\alpha} = {\sf v}_{_\alpha}$ are the vertices of the Weyl alcove. Recall 
\beqa\label{explicitexpr}
 e({\theta_{\alpha}}) :=  2.\big({dim_{_\bc}(G/P_{_\alpha}) - \mu(\alpha)}\big)
\eeqa
where $P_{_\alpha}$ is the maximal parabolic subgroup of $G$ associated to $\alpha$ and 
\beqa\label{mualpha}
\mu(\alpha) = \#\{r \in R^{+} \mid r = c_{\alpha}.\alpha +  \sum_{\beta \neq \alpha} x_{\beta}.\beta \}.
\eeqa
Let $L_{_\alpha}$ be the Levi subgroup of the closed fibre of the Bruhat-Tits group scheme $\sf {P}_{_{\ul{v}}}$. Recall that $L_{_\alpha}$ are all semisimple. 
\blem We have the following relation: for each simple root $\alpha$,
\beqa
{\tt dim}(G) = e(\theta_{_\alpha}) + {\tt dim}(L_{_\alpha})
\eeqa
and hence
\beqa
{\tt dim}(\text{Bun}^{^{\sf L}}_{_{C}}({\sf Q}_{_{\ul{{\sf v}}}})) = {\tt dim}(G) (g-1).
\eeqa
\elem
\bpr This is computational and I have used the tables. I don't see any general argument for this. \epr
\brem The group scheme ${\sf Q}_{_{\ul{{\sf v}}}}$ on the irreducible nodal curve will be non-reductive for the case of non-hyperspecial $\ul {\sf v}$, while the hyperspecial cases will be semisimple group schemes which are not globally split. The groups of type {\tt A}$_{_n}$ will not give the exotic examples since all maximal parahorics are hyperspecial. \erem

\brem (a bit imprecise!)  The component $\mathfrak G_{_0}$ consists of torsors on $C^{^{(d)}}$ which are $G$-bundles with parabolic structures at the nodes. In terms of the surface $S^{^{(d)}}$, the admissible pairs $({\qh^{^G}_{_{\tau, S^{^{(d)}}}}}, \ce)$ are such that the representation $\rho:\smu \to G$ is of type $\tau$. Consider the set of Giseker torsors $\ce$ on  $C^{^{(d)}}$ for $1 \leq d \leq \ell$ such that $\ce \mid_{_{\tC}} \simeq \tC \times G$. This set is bijective to the stack theoretic compactification of $G$ constructed by Marten-Thaddeus \cite{mt}. By \cite[page 94]{mt} and \cite[Remark 5]{solis}, it follows that stable $G$-bundle chains on $C^{^{(d)}}$ correspond to admissible torsors $({\qh^{^G}}_{_{\tau}}, \ce)$ when $G = \text{GL}(n)$.  One can  assume firstly that $d \leq \ell$ and then with a bit more work, one can check that  the Martens-Thaddeus stack is a closed substack of the principal component $\mathfrak G_{_0}$. It seems likely that the principal component is isomorphic to a bundle on $\text{Bun}_{_{G}}(\tC)$ with fibres the Martens-Thaddeus stack for $G$. 
\erem

\subsubsection{Example: $G = \tsl(2)$} Let $\eta:\tsl(2) \hra \tgl(2)$ be the standard inclusion. We list the basic representation types. Let $\rho_{_d}:{\sf\mu_{_d}} \to \tsl(2)$ be given by
$ \zeta_{_d} \mapsto \left(\begin{array}{cc} \zeta_{_d} & 0 \\0 & \zeta_{_d}^{^{-1}} \\ \end{array}\right)$, where $\zeta = e^{^{2i\pi/d}}$. 

\noindent
$d = 2$: the representation  $ \zeta_{_2} \mapsto \left(\begin{array}{cc} \zeta_{_2} & 0 \\0 & \zeta_{_2} \\ \end{array}\right)$ is central and hence the group scheme $\qh^{^{G}}_{_{\tau, {\sf N}^{^{(1)}}}}$ on ${\sf N}^{^{(1)}}$, restricted to $C^{^{(1)}}$  is obtained by gluing a parahoric group scheme  $\sf P$ on $\tC$ (which is the maximal parahoric $\sf{P}_{_{\ul{{\sf v}}}}$ group scheme near the two marked points) with the constant group scheme $\sf{P}_{_{\ul{{\sf v}},o}} \times \bp^{^1}$ on the single rational component, where $\sf{P}_{_{\ul{{\sf v}},o}} $ is the closed fibre of the parahoric group scheme $\sf{P}_{_{\ul{{\sf v}}}}$.  Being hyperspecial, in these cases, $\sf{P}_{_{\ul{{\sf v}},o}}$ is isomorphic to $\text{PGl}(2)$.

Torsors are obtained from equivariant torsors $E$ on $\ba^{^2}$ for the $\mu_{_2}$-action by the process $\text{Inv} \circ f_{_*}(E)$. These in turn give {\em laced} torsors on the normalization $\tC$ \eqref{2-lacedbundle} \eqref{bttorsorlacedtorsor}. Since the parahoric is hyperspecial, the lacing is simply an isomorphism of the fibres of the torsors on the normalization which in turn give an object in $\text{Bun}^{^{\sf L}}_{_{C}}({\sf Q}_{_{\ul{{\sf v}}}})$. 

\brem If parahoric is non-hyperspecial, the identification of the fibres is via the centralizer of the image of $\rho$ and this translates as identification of the associated $L_{_\alpha}$-torsor and hence an element in the adjoint group of $L_{_\alpha}$. Since the $L_{_\alpha}$'s are semisimple  the dimension is that of $L_{_\alpha}$.\erem

\noindent
$d = 3$: the representation  $ \zeta_{_3} \mapsto \left(\begin{array}{cc} \zeta_{_3} & 0 \\0 & \zeta_{_3}^{^2} \\ \end{array}\right)$. This case gives the bundles on the surface ${\sf N}^{^{(2)}}$ with closed fibre $F^{^{(2)}}$ having two $\bp^{^1}$'s. The simple root $\alpha$ on the maximal torus of $\tsl(2)$ sends $\left(\begin{array}{cc} \zeta_{_3} & 0 \\0 & \zeta_{_3}^{^2} \\ \end{array}\right) \to \zeta_{_3}^{^2}$. The induced character $\chi$ of $\Gamma_{_3}$ is $\zeta_{_3} \mapsto \zeta_{_3}^{^2}$ and this corresponds to second $\bp^{^1}$ on $C^{^{(2)}}$.
In other words, the group scheme $\qh^{^{G}}_{_{\tau, {\sf N}^{^{(2)}}}}$ has non-trivial parahoric structure on $E_{_2}$ and is the constant group scheme on $E_{_1}$. The representation $ \zeta_{_3} \mapsto \left(\begin{array}{cc} \zeta_{_3}^{^2} & 0 \\0 & \zeta_{_3} \\ \end{array}\right)$ will produce the other case.

The associated rank 2 vector bundle has $\co \oplus \co(1)$ on each $\bp^{^1}$, where the $\co(1)$ on the first $\bp^{^1}$ is associated to the character $ \zeta_{_3} \mapsto \zeta_{_3}$ and on the second to the character $ \zeta_{_3} \mapsto \zeta_{_3}^{^2}$ using the McKay correspondence. Torsors on the closed fibre correspond to torsors on $\tC$ for the group scheme ${\sf P}_{_I}$, the Iwahori at the two marked points together with the lacing data. The flags are full flags and the lacing data disappear since the Levi is the maximal torus and "modulo centre" is trivial.

So the dimension of this stratum is: dimension of the space of  torsors  for the group scheme which is the Iwahori at two points on $\tC$ =  ${\tt dim}(G)(g-2) + 2 {\tt dim}(G/B)$ = ${\tt dim}(G)(g-1) - \ell$, where $\ell = n-1$ in this case, i.e. the stratum is of codimension $\ell$. In the case of $\tsl(2)$, the closed fibre of the stack is the union of two smooth components, which meet at the last stratum of the torsors for the group scheme on $\tC$ which is the Iwahori at the two points. The smoothness of the components can also be deduced by deformation arguments. The miniversal space for a Gieseker torsor $(C^{^{(d)}}, {\qh^{^G}}_{_{\tau}}, \ce)$ will be such that $d = \ell$ (in the Iwahori situation) and in the case of $\tsl(2)$ it will be $C^{^{(1)}}$, i.e. with a single $\bp^{^1}$. Hence, the number of components meeting the Iwahori-type bundles is $2$, corresponding to the two nodes on $C^{^{(1)}}$.
\section{On Mumford's toroidal realization of buildings} With notations as before, we work with the split group scheme $G_{_A} = G \times A$ and $T_{_A} \subset G_{_A}$ a split torus over $A$, i.e. $T_{_A} \simeq \bg_{_{m,A}}^{^{\ell}}$. In 
\cite{kkms} towards the very end Mumford gives a beautiful construction of the geometric realization of both the absolute and the relative case of Tits buildings via toroidal embeddings. We will talk of the relative case alone here. 

The choice of $T$ entails a choice of an apartment $\ca_{_T}$ and the choice of the root system entails a choice of an origin in $\ca_{_T}$ together with an alcove ${\bf a} \subset \ca_{_T}$. Let $\eR = \eR^{^+} \cup \eR^{^-}$ be the decomposition of the roots $\eR$ into positive and negative roots. An alcove is given by:
\beqa
{\bf a} := \{x \in \ca_{_T} \mid 0 \leq \langle \alpha,x \rangle \leq 1, \alpha \in \eR^{^+}  \}.
\eeqa
The alcoves give the top dimensional simplices of the polyhedral decomposition of $\ca_{_T}$ in terms of the affine hyperplanes. Define $\boldsymbol \sigma \subset \ca_{_T} \times \br$ to be the cone over ${\bf a} \times (1)$. This gives an affine torus embedding $T_{_A} \subset T_{_{\boldsymbol \sigma}}$. Hence we get the fibre bundle $G_{_A} \times ^{^{T_{_A}}} T_{_{\boldsymbol \sigma}}$ associated to the principal $T_{_A}$-bundle $G_{_A} \to G_{_A}/T_{_A}$. On the generic fibre we have the identification $T_{_K} = (T_{_{\boldsymbol \sigma}})_{_K}$ and hence $G_{_K} = (G_{_A} \times ^{^{T_{_A}}} T_{_{\boldsymbol \sigma}})_{_K}$. Mumford then defines the relative embedding:
\beqa\label{mumemb}
{\bar G}_{_A} := \bigcup_{_{x \in G(k(\!(t)\!)}} (G_{_A} \times ^{^{T_{_A}}} T_{_{\boldsymbol \sigma}}).x
\eeqa
where notation for the action of $G(k(\!(t)\!))$ on  $G_{_A} \times ^{^{T_{_A}}} T_{_{\boldsymbol \sigma}}$ stands for embeddings of $G_{_K}$ by a translation by $x$. These can be suitably glued to get a separated scheme over $A$ (see \cite[206]{kkms}).

The salient feature of the toroidal embedding $G_{_A} \subset {\bar G}_{_A}$ is that $G_{_K} = {\bar G}_{_K}$, for each $x \in G(K)$, the right multiplication extends to give an automorphism of ${\bar G}_{_A}$. Finally, 
the strata of ${\bar G}_{_A} \setminus G_{_A}$ are precisely the parahoric subgroups of $G(K)$. This bijection extends to an isomorphism of the graph of the embedding $G_{_A} \subset {\bar G}_{_A}$ with the Bruhat-Tits building of $G_{_A}$. The aim of the present section is to give one point of contact between this construction and the stack $\text {Gies}_{_{G}}(C_{_A})$ constructed earlier. 

Let $(C^{^{(\ell)}}, \qh_{_{\tau}}^{^{G}}, \ce)$ be a Gieseker torsor, where $\tau$ is the type of $\rho:\mu_{_{\ell}} \to G$ which is such that all the characters of $\mu_{_\ell}$ occur {\sl precisely once} in $\rho$. With the choice of the maximal torus $T \subset G$ and the root systen $\eR$ we see that $\tau$ gives a point $\theta_{_{\tau}} \in \bf a$. The group scheme $\qh_{_{\tau}}^{^{G}}$ when restricted to the normalization $\tC$ has the property that in the analytic neighbourhood of both the points $c_{_1}, c_{_2}$, the group scheme is the Iwahori group scheme. The Iwahori structure is a consequence of the distribution of characters of $\mu_{_\ell}$ in $\rho$. 

The Gieseker torsor $(C^{^{(\ell)}}, \qh_{_{\tau}}^{^{G}}, \ce)$ gives a point of the scheme ${\tt Y}^{^G}$ (\ref{deforgieseker}). Furthermore,  in an \'etale neighbourhood of this point we have a morphism to  $B[\ell]$. Note that the standard model $W[\ell]$ in our setting is such that $B[\ell] \simeq {\hat\ba}^{^{\ell +1}}$ but it can be defined over the affine space ${\mathbb A}^{^{\ell +1}}$. We work in the latter setting here. 

The base ${\mathbb A}^{^{\ell +1}}$ of the standard model $W[\ell] \to {\ba}^{^{\ell +1}}$ is affine toric variety and as a toric variety over $\bc$ we can identify it with $T_{_{\boldsymbol \sigma}}$. The big cell in $G_{_A}$ gives an open subset where the principal $T_{_A}$-bundle $G_{_A}\to G_{_A}/T_{_A}$ is a product. Thus, in an open subset  we can identify the associated fibre space $G_{_A} \times ^{^{T_{_A}}} T_{_{\boldsymbol \sigma}}$ as  a $\bc$-scheme with $U^{^-} \times {\ba}^{^{\ell +1}} \times U^{^+}$. The family $W[\ell]$ pulled-back by the projection thus gives a modification on the open subset $U^{^-} \times {\ba}^{^{\ell +1}} \times U^{^+}$. The Gieseker torsor $(C^{^{(\ell)}}, \qh_{_{\tau}}^{^{G}}, \ce)$ spreads to a formal neighbourhood of the origin in ${\ba}^{^{\ell +1}}$. This gives a morphism $U^{^-} \times {\hat\ba}^{^{\ell +1}} \times U^{^+} \to {\tt Y}^{^G}$ which send $0$ to $(C^{^{(\ell)}}, \qh_{_{\tau}}^{^{G}}, \ce)$. All in all this shows that ${\tt Y}^{^G}$ is formally smooth to the Mumford embedding \eqref{mumemb} at the point $(C^{^{(\ell)}}, \qh_{_{\tau}}^{^{G}}, \ce)$.

\begin{center}
\underline{\sc Part II}
\end{center}
\section{Twisted curves and torsion-free sheaves}\label{twistedstuff}
\subsubsection{Goals of Part II} In this part of the paper, we  focus on the single nodal curve $(C,c)$ and its normalization $\nu:(\tC, {\bf c}) \to (C,c)$ \eqref{normalization}. The final aim is to describe the Gieseker torsors on the nodal curves $C^{^{(d)}}$ in terms of its restriction to $\tC$ and thereby get a notion of "semistability"  for them. Recall that in \cite{ns1} and \cite{kausz}, the direct image under $p:C^{^{(d)}} \to C$ relates Gieseker vector bundles to torsion-free sheaves and gives a morphism of stacks. (Semi)stability for Gieseker vector bundles is then defined in terms of two ingredients, the (semi)stability of torsion-free sheaves and a  (semi)stability using a relative polarization for this morphism. Our approach is modelled after this one. The first goal is to understand the classical (semi)stability of torsion-free sheaves  from a new standpoint. {\sl The aim is to avoid going to sub-objects to test stability}. 

\subsubsection{The basic setup}
Let  $N_{_c} = \spec~\frac{k\llbracket x,y \rrbracket}{(x.y)}$ be the analytic neighbourhood at the node on $(C,c)$ and let  $N_{_0}$ denote $\spec~\frac{k\llbracket u,v \rrbracket}{(u.v)}$ with coordinates $u,v$. We can express $N_{_0} = D_{_0} \cup D'_{_0}$, where $D_{_0}$ and $D'_{_0}$ are identified with discs with $0$ as origin. Let  $\sf\mu_{_d}$  act on $N_{_0}$ by sending 
\beqa\label{chinought}
u \mapsto \zeta.u,~~~ v \mapsto \zeta^{^{-1}}.v
\eeqa 
where $\zeta$ is a primitive $d^{th}$--root of unity and $v$ is the local coordinate of $D_{_0}$ and $u$ for $D'_{_0}$. The quotient morphism 
\beqa\label{localformalquotient}
\sigma: N_{_0} \to N_{_c} = N_{_0}/\sf\mu_{_d} = D_{_0}/\sf\mu_{_d} \bigcup D'_{_0}/\sf\mu_{_d}
\eeqa
is given by  $\sigma(u)  = u^{^d} = x$, $\sigma(v)  = v^{^d} = y$.  

\subsubsection{Twisted curves} Let $d \leq {\mathtt m}(G)$ \eqref{somebounds}.
Let $\cc_{_d}$ be a twisted nodal curve  in the sense of \cite[Definition 2.1]{ollson}.  Note that $\cc_{_d}$ is an algebraic stack with $C$ as its coarse space, and we have the morphism $\sigma:\cc_{_d} \to C$ (see \cite{ollson}). Assume that analytic locally at the node $c \in C$ it is given by $N_{_0} \to [N_{_0}/{\sf\mu_{_d}}] \to N_{_c}$. 

Fix $\ce_{_0}$ the $(\sf\mu_{_d}, G)$-torsor on $N_{_0} \subset D$ obtained by restricting $\ce_{_D}$ \eqref{egtau}, given by a representation $\rho:\sf\mu_{_d} \to G$ {\em of local type} $\tau$. On $N_{_0} \times G$, the $\sf\mu_{_d}$-action is given by:
\beqa\label{localaction}
\gamma.(u,g) = (\zeta.u, \rho(\gamma).g), 
\\
\gamma.(v,g) = (\zeta^{^{-1}}.v, \rho(\gamma).g)
\eeqa
Let $\ce_{_0}(G)$ be the equivariant group scheme on $N_{_0}$ of type $\tau$. This is therefore fixed throughout.
\bdefe\label{torsoronDM} A $G$-torsor $\ce$ of local type $\tau$ on $\cc_{_d}$ is  the datum $(E', E_{_0}, \text{\cursive g}~)$, where
\begin{itemize}
\item $E'$ is a $G$-torsor  on the punctured curve $C-c$, 
\item $E_{_0}$ is a $\ce_{_0}(G)$-torsor,  
\item $\text{\cursive g}$ is a $\sf\mu_{_d}$-invariant gluing function
\beqa\label{glueabove1}
\text{\cursive g}: N_{_0}^{^*} \to G
\eeqa 
which gives an isomorphism: $\text{\cursive g}:\sigma^*(E'\mid_{_{N^{^*}_{_c}}}) \simeq  E_{_0} \mid_{_{N_{_0}^{^*}}}.$
\end{itemize} 
\edefe
We observe that $E_{_0}$ being a $\ce_{_0}(G)$-torsor encapsulates the statement that  $E_{_0}$ is a $(\sf\mu_{_d}, G)$-torsor of type $\tau$.

There is an obvious  notion of isomorphism of such torsors. {We also note that for sheaves on the twisted curve $\cc_{_d}$, we have the natural push-forward $\sigma_{_*}$. Locally on $N_{_0}$ this is the invariant push-forward ($\text{Inv} \circ \sigma_{_*}$)}.
\subsection{Torsion-free sheaves on $(C,c)$ as $\text{GL}(n)$-bundles on twisted curves $\cc_{_d}$} 
In this subsection we make some remarks on torsion-free sheaves on $C$ and their lifts to the stack $\cc_{_d}$. These remarks are key to the idea behind the general {\tt tf}-semistability of Gieseker torsors which we define later.

Let $\ca$ be a torsion-free sheaf on $C$ of rank $n$. We view $\ca$ as the datum $(V', F_{_0}, \cug')$, where $V'$ is a vector bundle of rank $n$ on $C-c$, $F_{_0}$ a torsion-free sheaf on $N_{_c}$ of rank $n$ and $\cug':V'\mid_{_{N_{_c}^*}} \simeq F_{_0}\mid_{_{N_{_c}^*}}$ a gluing function. We can express this equivalently as the datum of principal bundles as $(E', F_{_0}, \cug')$, where $E'$ is the frame bundle of $V'$ and $F_{_0}\mid_{_{N_{_c}^*}}$ (by an abuse of notation) stands for the  frame bundle of the vector bundle $F_{_0}\mid_{_{N_{_c}^*}}$ and $\cug'$ is an isomorphism of principal bundles on $N_{_c}^*$. The local type of $F_{_0}$ gives us a $(\sf\mu_{_d},\text{GL}(n))$-bundle $E_{_0}$ coming from a representation $\rho:\sf\mu_{_d} \to T_{_n} \subset \text{GL}(n)$ such that $\text{Inv} \circ \sigma_{_*}(E_{_0}(k^{^n})) = F_{_0}$. Thus we get a triple $(E', E_{_0}, \cug)$ with $\cug = \sigma^{^*}(\cug')$. Equivalently we get a $\text{GL}(n)$-torsor $\ce$ on $\cc_{_d}$ from the torsion-free sheaf $\cf$ and there is nothing unique about $\ce$ except that  $\sigma_{_*}(\ce(k^{^n})) \simeq \ca$. 

By the general remarks in the appendix below (see \eqref{parabtrio}), giving a $1$-PS produces weighted filtrations on vector spaces. From the standpoint of the datum above, given a filtration of $\ca$ by saturated subsheaves,  the associated graded sheaf can be recovered from data on the torsor $\ce$. More precisely, suppose we are given  the weighted filtration 
\beqa\label{satfiltontf}
{\ca}_{_\lambda}^{\bullet} : 0 \subsetneq \ca_{_{1}} \subsetneq \ldots \ca_{_{s}} \subsetneq \ca_{_{s+1}} = \ca,
\eeqa  
and let $\lambda:\bg_{_m} \to \text{GL}(n)$ be the $1$-PS coming from this datum. Further, let $P = P(\lambda) \subset \text{GL}(n)$ be the induced parabolic subgroup. Restricting the filtration \eqref{satfiltontf} to $C-c$, gives a weighted filtration of the locally free sheaf $E'(k^{^n}) = \ca\mid_{_{C-c}}$:
\beqa\label{satfiltontflf}
 0 \subsetneq A_{_{1}} \subsetneq \ldots A_{_{s}} \subsetneq A_{_{s+1}} = \ca\mid_{_{C-c}},
\eeqa
which is equivalent to giving a reduction of structure group $E'_{_P} \subset E'$ to $P$.
By saturating this filtration in the sheaf $\ca$ we get back the {\em weighted filtration} \eqref{satfiltontf} on $C$. In other words, a reduction of structure group $E'_{_P} \subset E'$ gives a saturated filtration of $\ca$ by saturated subsheaves and a tuple $\ul{\epsilon} = (\epsilon_1, \ldots, \epsilon_s)$ of positive rational numbers which are recovered from $\lambda$ \eqref{parabtrio}. 

How does one recover the associated graded sheaf ${\tt gr}_{_\lambda}({\ca}) = \bigoplus \ca_{_{j+1}}/\ca_{_j}$? Observe that this is also torsion-free, but {\sl its local type, i.e. its type restricted to $N_{_c}$ at the node, depends on the filtration \eqref{satfiltontflf}, which in turn comes by a process of saturation}. 

We proceed as follows: firstly,  locally, from the torsion-free sheaf ${\tt gr}_{_\lambda}({\ca})$, we can get a $\sf\mu_{_d}$-vector bundle $V_{_{0,P}}$ on $N_{_0}$ such that 
\beqa
\text{Inv} \circ \sigma_{_*}(V_{_{0,P}}) \simeq {\tt gr}_{_\lambda}({\ca})\mid_{_{N_{_c}}}.
\eeqa
The $\sf\mu_{_d}$-vector bundle $V_{_{0,P}}$ comes from a representation $\varphi:\sf\mu_{_d} \to T_{_n} \hra \text{GL}(n)$ and this gives a $(\sf\mu_{_d},T_{_n})$-bundle $E^{^{\varphi}}_{_0}$  on $N_{_0}$; since the $\sf\mu_{_d}$-action on $N_{_0}^{^*}$ is free, we have a canonical isomorphism:
\beqa\label{canoverpunctured}
E^{^{\varphi}}_{_0}\mid_{_{N_{_0}^*}} \simeq E_{_0}\mid_{_{N_{_0}^*}}.
\eeqa
The filtration \eqref{satfiltontf} on $\ca$ produces an obvious filtration on ${\tt gr}_{_\lambda}({\ca})$ by deleting at each step a summand. These coincide on $C-c$. i.e., this filtration, when restricted to $C-c$ produces a gluing function $\cug^{^{\varphi}}_{_P}:\sigma^{^*}(E'_{_P}\mid_{_{N_{_c}^*}}~) \simeq E^{^{\varphi}}_{_0}(P) \mid_{_{N_{_0}^*}}$, such that it induces $\cug$ when the structure group is extended from $P$ to $G$. Whence, we firstly get a $P$-torsor $(E'_{_P},E^{^{\varphi}}_{_0},\cug^{^{\varphi}}_{_P})$ from ${\tt gr}_{_\lambda}({\ca})$. 

{\em The recipe to recover ${\tt gr}_{_\lambda}({\ca})$ is as follows}. Let $H = H(\lambda)$ be the canonical Levi quotient of $P$. By the remarks made after \eqref{assocgradedlambda}, it follows easily that 
\beqa
\sigma_{_*}((E'_{_H}, E^{^{\varphi}}_{_{0,H}}, \cug^{^{\varphi}}_{_H})) = {\tt gr}_{_\lambda}({\ca})
\eeqa
where $(E'_{_H}, E^{^{\varphi}}_{_{0,H}}, \cug^{^{\varphi}}_{_H})$ is a $H$-torsor of local type $\varphi$ on $\cc_{_d}$, the objects being natural extension of structure groups to $H$.

A weighted slope is assigned by Schmitt (\cite{schmitt}) to the weighted filtration \eqref{satfiltontf} :
\beqa\label{wtdslope}
L\big(\ca_{_\lambda}^{\bullet}, \ul{\epsilon}\big) := \sum_{i = 1}^{^s} \epsilon_{_i}\bigg\{{\tt deg}_{_{C}}(\ca)\cdot \text{rk}~\ca_{i} - {\tt deg}_{_{C}}(\ca_{i}) \cdot \text{rk}~\ca)\bigg\}.
\eeqa
and an obvious weighted slope $L\big({\tt gr}_{_\lambda}({\ca}^{\bullet}), \ul{\epsilon}\big)$ for the associated graded sheaf ${\tt gr}_{_\lambda}({\ca})$. It is easy to see that 
\beqa\label{trivialeq}
L\big(\ca_{_\lambda}^{\bullet}, \ul{\epsilon}\big) = L\big({\tt gr}_{_\lambda}({\ca}^{\bullet}), \ul{\epsilon}\big)
\eeqa
We have the obvious but useful reformulation of semi(stablity) which circumvents going to sub-objects.
\blem
A torsion-free sheaf $\ca$ on $C$ is semi(stable) if and only if for every weighted filtration \eqref{satfiltontf}, $L\big({\tt gr}_{_\lambda}({\ca}^{\bullet}), \ul{\epsilon}\big) (\geq) 0$.\elem

\section{On semistability of $G$-torsors on twisted curves}
Let  $\ce = (E', E_{_0}, \text{\cursive g})$, be a  $G$-torsor on $\cc_{_d}$  of local type $\tau$ coming from a representation $\rho:\sf\mu_{_d} \to T \subset G$ (see \eqref{torsoronDM}). Thus, $E_{_0}$ can be taken as a $(\sf\mu_{_d},T)$-bundle coming from $\rho$ and $\cug:\sigma^*(E'\mid_{_{N_{_c}^*}}) \simeq E_{_0}\mid_{_{N_{_0}^*}}$.  The $(\sf\mu_{_d},T)$-bundle $E_{_0}$ gives a $\sf\mu_{_d}$-line bundle decomposition $\oplus E_{_{0}}(\alpha_{_i} \circ \rho)$ coming from the simple roots $\alpha_{_i}:T \to \bg_{_m}$ and we get a rank $\ell$  torsion-free sheaf 
\beqa
\cf_{_0} := \bigoplus \text{Inv} \circ \sigma_{_*}\big({E_{_{0}}(\alpha_{_i} \circ \rho)}\big)
\eeqa
on $N_{_c}$. Let $\lambda:\bg_{_m} \to T \subset G$ be a $1$-PS  which gives a filtration:
\beqa\label{satfiltontfell1}
 0 \subsetneq F_{_{1}} \subsetneq \ldots F_{_{s}} \subsetneq F_{_{s+1}} = \cf_{_0}\mid_{_{N_{_c}^{^*}}}
\eeqa
of the locally free sheaf $\cf_{_0}\mid_{_{N_{_c}^{^*}}}$.
By saturation in $\cf_{_0}$ we get a filtration:
\beqa\label{satfiltontfell2}
 0 \subsetneq \cf_{_{1}} \subsetneq \ldots \cf_{_{s}} \subsetneq \cf_{_{s+1}} = \cf_{_0},
\eeqa
by subsheaves (with quotients torsion-free) and we also get the associated graded sheaf ${\tt gr}_{_\lambda}(\cf_{_0})$ on $N_{_c}$. We now follow the earlier procedure. Firstly, the local type of ${\tt gr}_{_\lambda}(\cf_{_0})$ gives a new $\varphi:\sf\mu_{_d} \to T $ (dependent on the $1$-PS $\lambda$) and a $(\sf\mu_{_d},T)$-bundle $E^{^{\varphi}}_{_0}$ on $N_{_0}$. The induced bundle $\oplus E^{^{\varphi}}_{_0}(\alpha_{_i} \circ \varphi)$ recovers ${\tt gr}_{_\lambda}(\cf_{_0})$ as $\bigoplus \big(\text{Inv}\circ \sigma_{_*}( E^{^{\varphi}}_{_0}(\alpha_{_i} \circ \varphi))\big)$.

Let $\ce = (E', E_{_0}, \text{\cursive g}) $ be a $G$-torsor on $\cc_{_d}$ with $E_{_0}$ an $\ce_{_0}(G)$-torsor coming from  $\rho:\sf\mu_{_d} \to T \subset G$ \eqref{torsoronDM}. So $E_{_0}$ comes with a choice of reduction of structure group to $T$.  A reduction   of structure group $\ce_{_P}$ of $\ce$  to $P = P(\lambda)$ comprises of  the datum $(E'_{_P},E_{_0}, { \text{\cursive g}_{_{P}}})$, where $E'_{_P}$ is a reduction of structure group of $E'$ to the parabolic subgroup $P$ on $C - c$ and    ${\text{\cursive g}_{_P}}: N_{_c}^* \to P$ is a function which gives an identification:
\beqa\label{phew1}
{\text{\cursive g}_{_P}}:\sigma^*\big({E'_{_P} \Big\vert _{_{N_{_c}^*}}}\big) \simeq (E_{_0} \times^{^T} P) \Big\vert _{_{N_{_0}^*}}
\eeqa
which glues the $P$-torsors along $N_{_0}^*$ with the constraint that the composite ${\text{\cursive g}_{_P}}:{N_{_0}^*} \to P \hra G$ equals  ${\text{\cursive g}}.$

The $1$-PS $\lambda$  gave rise to a new  $\varphi:\sf\mu_{_d} \to T$ and a new $(\sf\mu_{_d},T)$-bundle $E^{^{\varphi}}_{_0}$ on $N_{_0}$. Since the $\sf\mu_{_d}$-action is free away from $0$, it follows that we have canonical identification:
\beqa\label{phew}
E_{_0} \mid_{_{N_{_0}^*}} \simeq E^{^{\varphi}}_{_0}\mid_{_{N_{_0}^*}}.
\eeqa
Since $T \subset P$, by extending structure groups via $P \to H$ and by \eqref{phew}, we get an $H$-torsor 
\beqa\label{ephih}\ce^{^{\varphi}}_{_H} := (E'_{_H}, E^{^{\varphi}}_{_0}, { \text{\cursive g}_{_{H}}})\eeqa 
of local type $\varphi$ on $\cc_{_d}$, where
\beqa\label{phew2}
{\text{\cursive g}_{_H}}:\sigma^*\big({E'_{_H} \Big\vert _{_{N_{_c}^*}}}\big) \simeq E^{^{\varphi}}_{_0}(H) \Big\vert _{_{N_{_0}^*}}
\eeqa
is the one induced by ${\text{\cursive g}_{_P}}$. Moreover, as we saw above,
\beqa
{\tt gr}_{_\lambda}(\cf_{_0}) := \bigoplus_{_i} \big(\text{Inv} \circ \sigma_{_*}( E^{^{\varphi}}_{_0}(\alpha_{_i} \circ \varphi))\big)
\eeqa
\bdefe\label{genericreductiondatum} We call $\ce^{^{\varphi}}_{_H}$ \eqref{ephih} the twisted $H$-torsor of local type $\varphi$ associated to $\ce_{_P}$.
\edefe
Let $\eta:G \hra \text{GL}(W)$ be a faithful representation and let $T_{_W} \subset \text{GL}(W)$ be a maximal torus such that $\eta:T \hra T_{_W}$. Let $\ce$ be a $G$-torsor on $\cc_{_d}$. Then we get the associated $\text{GL}(W)$-torsor $\ce(\text{GL}(W))$ and associated vector bundle  $\ce(W)$. We can get torsion-free sheaves on $C$ via $\sigma:\cc_{_d} \to C$. Let $\cf_{_W} := \sigma_{_*}(\ce(W))$.

Let $\lambda:\bg_{_m} \to T \subset G$ be a one-parameter subgroup, which gives parabolic subgroups $P = P(\lambda)$ and Levi $H$ and induced $1$-PS $\eta \circ \lambda:\bg_{_m} \to T \hra T_{_W}$ and canonical inclusions $P \subset P_{_W}$ and $H \hra H_{_W}$.

When the torsion-free sheaf $\cf_{_W}$ is restricted to $C-c$, the one-parameter subgroup $\eta \circ \lambda$ gives a weighted filtration 
\beqa
0 \subsetneq F_{_{1}} \subsetneq \ldots F_{_{m}} \subsetneq F_{_{m+1}} = \cf_{_W} \Big\vert _{_{C-c}}
\eeqa
of the locally free sheaf $\cf_{_W} \Big\vert _{_{C-c}}$  which by saturation gives a filtration 
\beqa\label{filtonD1}
0 \subsetneq \cf_{_{1}} \subsetneq \ldots \cf_{_{m}} \subsetneq \cf_{_{m+1}} = \cf_{_W},
\eeqa
by subsheaves (with torsion-free quotients). This gives the associated graded torsion-free sheaf:
\beqa
{\tt gr}(\cf_{_W}) := \bigoplus_{_{j = 1}}^{^{m}} \cf_{_{j+1}}/\cf_{_j}.
\eeqa
Note that, as $H_{_W}$-modules (and hence as $H$-modules), we have $W \simeq {\tt gr}(W)$ (filtered via $\eta \circ \lambda$). Note further that the weighted filtration on $W$ induced by $\lambda$ also gives a canonical weighted filtration $({\tt gr}(W)^{^{\bullet}}, \ul{\epsilon})$ on ${\tt gr}(W)$. Thus we get $\ce^{^{\varphi}}_{_{H}}({\tt gr}(W)) = \bigoplus_{_{j = 1}}^{^{m}} \ce^{^{\varphi}}_{_{H}}(W_{_{j+1}}/W_{_j})$. Therefore, for each $\varphi$  we have an associated vector bundle $\ce^{^{\varphi}}_{_H}({\tt gr}(W))$. By taking push-forward by $\sigma:\cc_{_d} \to C$ we have the isomorphism 
\beqa\label{leviidently}
 \sigma_{_*}\big({\ce^{^{\varphi}}_{_H}({\tt gr}(W))}\big) =  \sigma_{_*}\big(\bigoplus_{_{j = 1}}^{^{m}} \ce^{^{\varphi}}_{_{H}}(W_{_{j+1}}/W_{_j})\big) \simeq \bigoplus_{_{j = 1}}^{^{m}} \cf_{_{j+1}}/\cf_{_j} = {\tt gr}(\cf_{_W}),
 \eeqa
of torsion-free sheaves on $C$. We also get the obvious 
weighted filtration matching \eqref{filtonD1} term by term:
\beqa\label{wtdfiltongmodules0}
 0 \subsetneq \ce_{_{H,1}} \subsetneq  \ldots \ce_{_{H,m}} \subsetneq \ce_{_{H,m+1}} = \ce^{^{\varphi}}_{_{H}}({\tt gr}(W))
\eeqa
If we have a notion of "degree" of associated vector bundle on the twisted curve $\cc_{_d}$, then following Schmitt (\cite{schmitt}) we also have a  weighted slope:
\small
\beqa
L\big({\ce^{^{\varphi}}_{_H}({\tt gr}(W)})\big) := \sum_{i = 1}^{^s} \epsilon_{_i}\bigg\{{\tt deg}_{_{\cc_{_d}}}({\ce^{^{\varphi}}_{_H}({\tt gr}(W))})\cdot \text{rk} \ce_{_{H,i}} - {\tt deg}_{_{\cc_{_d}}}~\ce_{_{H,i}} \cdot \text{rk}({\ce^{^{\varphi}}_{_H}({\tt gr}(W))})\bigg\}.
\eeqa
\normalsize
By \eqref{parabquadro}, associated to $\lambda$, we have an anti-dominant character $\chi_{_\lambda}$ of the parabolic subgroup $P$ (and the Levi $H$) "dual" to $\lambda$. Suppose further that the "degree"  satisfies the  equality (see \eqref{formalequalityforpar}):
\beqa\label{formalequality}
{\tt deg}_{_{\cc_{_d}}}({\ce^{^{\varphi}}_{_H}}(\chi_{_\lambda})) = L\big({\ce^{^{\varphi}}_{_H}({\tt gr}(W)})\big)
\eeqa
By \eqref{schmittbook1}, for any anti-dominant character $\chi$ of $P$, there is a positive rational $r$ such that $\chi_{_\lambda} = r.\chi$ and hence by \eqref{formalequality}, 
we deduce that for each anti-dominant $\chi$, $L\big({\ce^{^{\varphi}}_{_H}({\tt gr}(W)})\big)$ and ${\tt deg}_{_{\cc_{_d}}}({\ce^{^{\varphi}}_{_H}}(\chi))$ have the same sign. Thus, with our hypothetical "degree" plus \eqref{formalequality}, (see \eqref{degreeoftorsoroncc}) we have the definitions:
\bdefe\label{newtfss} 
\begin{enumerate}
\item A $G$-torsor $\ce$ of local type $\tau$ on the twisted curve $\cc_{_d}$ is called {\tt tf}-semi(stable) if  for every $1$-PS $\lambda:\bg_{_m} \to T$, and reduction of structure group $\ce_{_P}$, we have ${\tt deg}_{_{\cc_{_d}}}({\ce^{^{\varphi}}_{_H}}(\chi)) (\geq) 0$.
\item $\ce$ is $\eta$-semi(stable) if for every $1$-PS $\lambda:\bg_{_m} \to T \subset T_{_W}$, we have $$L\big({\ce^{^{\varphi}}_{_H}({\tt gr}(W)})\big) ~(\geq)~0$$
\end{enumerate}
\edefe
Our observations show then that we have a theorem analogous to the classical theorem of Ramanathan:
\bth\label{goodshow} A $G$-torsor $\ce$ of local type $\tau$ on the twisted curve $\cc_{_d}$ is {\tt tf}-semi(stable) if and only if it is $\eta$-semi(stable) for every $\eta:G \hra \text{GL}(W)$. \eeth
The next task is to show that there is a well-defined notion of "degree" on $\cc_{_d}$ \eqref{degreeoftorsoroncc} with some properties like \eqref{formalequality} (see \eqref{formalequalityforpar}), and that the above notion is geometric invariant theoretic. To achieve this we traverse a {\em "parabolic path"} via a Fourier-Mukai from bundles on $\cc_{_d}$ to {\tt laced} ones on the normalization $\tC$ which comes with a {\em balanced} parabolic structure.

\section{Laced torsors via {\sl Fourier-Mukai}  for torsors on twisted curves}\label{curvecase}    

Let  $q:(\tilde{\cc_{_d}}, {\bf z}) \to (\cc_{_d}, o)$ be the  normalization, i.e. the maximal reduced substack of $\cc_{_d} \times_{_C} \tC$. Then $\tilde{\cc_{_d}}$ is a twisted curve with two markings, the normalization $\tC$ of $C$ is its coarse space with the canonical morphism $f:(\tilde{\cc_{_d}}, {\bf z}) \to (\tC, {\bf c})$, and we have a diagram:
\beqa\label{stackynormalization}
\xymatrix{
\tilde{\cc_{_d}} \ar[r]^{^f} \ar[d]_{q} &
 \tC  \ar[d]^{\nu} \\
\cc_{_d} {\ar[r]_{_\sigma}} &  C
}
\eeqa
and the corresponding local picture:
\beqa\label{stackynormalizationloc}
\xymatrix{
[\tN/\smu] \ar[r]^{^f} \ar[d]_{q} &
 \tU  \ar[d]^{\nu} \\
[N_{_0}/\smu] {\ar[r]_{_\sigma}} &  N_{_c}
}
\eeqa
An analytic neighbourhood of $\tilde{\cc_{_d}}$ at $z_{_1}$ (resp. $z_{_2}$) gets identified with $[D_{_0}/\smu]$ (resp. $[D'_{_0}/\smu])$ with action given by \eqref{chinought} and similarly analytic neighbourhood of $\tC$ at $c_{_1}$ (resp. $c_{_2}$) gets identified with $D_{_0}/{\sf\mu_{_d}}$ (resp. $D'_{_0}/{\sf\mu_{_d}}$). Given the ramification data at  $c_{_i}$ we can get a smooth projective Kawamata cover 
\beqa\label{kawaZ}
f':Z \to \tC
\eeqa 
with same local ramification data as $f$. The pull-back $q^{^*}(\ce)$ gives a $G$-torsor on $\tilde{\cc_{_d}}$ plus a "descent datum'', i.e., a  $(\sf\mu_{_d},G)$-isomorphism 
\beqa\label{dpgbundle}
\text{\cursive i}:\eP_{_{z_{_1}}} \simeq \eP_{_{z_{_2}}}^{^\dagger},
\eeqa
where $\eP^{^\dagger}$ is the $(\sf\mu_{_d},G)$-torsor in a neighbourhood of $z_{_2}$ given by the local type $\tau$. Equivalently, we can take the "adjoint" group scheme $\ce(G)$ on $\cc_{_d}$   and let 
\beqa\label{egtau}
\ce(G,\tau) = q^{^*}(\ce(G)).
\eeqa 
Then $q^{^*}(\ce)$ is simply a $\ce(G,\tau)$-torsor on the twisted curve $\tilde{\cc_{_d}}$.

\brem\label{ondescentdatum} ({\em On the descent datum}) Observe that a  $(\sf\mu_{_d},G)$-isomorphism $\text{\cursive i}:\eP_{_{z_{_1}}} \simeq \eP_{_{z_{_2}}}^{^\dagger}$ gives a $G$-isomorphism of the quotients $\ce_{_{z_{_1}}}/{\tt I}_{_\rho} \simeq \ce_{_{z_{_2}}}/{\tt I}_{_\rho}$, where ${\tt I}_{_\rho} := \text{Im}(\rho) \subset G$.

The descent datum \eqref{dpgbundle} in the case of the group scheme $\ce_{_0}(G,\tau)$ translates as a $G$-isomorphism:
\beqa
\iota_{_h}:G/{\tt I}_{_\rho}\to G/{\tt I}_{_\rho},
\eeqa
given by an inner automorphism induced by an element $h \in \text{Cent}_{_G}(I_{_\rho})$ modulo the center of $\text{Cent}_{_G}(I_{_\rho})$. 

For instance, when $\rho:\sf\mu_{_d} \to G$ sends a generator $\zeta$ to a Borel-de Seibenthal element $g_{_\alpha}$, we see that $\text{Cent}_{_G}(I)$ is precisely the Levi quotient of the closed fibre of the maximal parahoric group scheme ${\sf P}_{_{\theta_{_{\alpha}}}}$. Hence a {\tt laced} maximal parahoric group scheme on $\tC$ \eqref{bttorsorlacedtorsor} is given by a maximal parahoric group scheme ${\sf P}_{_{\theta_{_{\alpha}}}}$ at the two marked points with an isomorphism of the Levi quotients of the closed fibres, modulo the center. Instead, if $\rho$ gives  the Iwahori structure, then descent datum becomes trivial since the Levi is abelian.
\erem

\bdefe\label{balancedparahoricgrpsch} For $\theta_{_\tau} \in \ca_{_T}$, we define the  {\tt balanced} parahoric group scheme on  $(\tC, {\bf c})$ as:
$\eG_{_{(\theta_{_{\tau}})}} := f_{_*}(\ce(G,\tau))$, with $\ce(G,\tau)$ as in \eqref{egtau}.
\edefe
We could use the Kawamata cover \eqref{kawaZ} to see that $\eG_{_{(\theta_{_{\tau}})}} := \text{Inv} \circ f'_{_*}\big[(\ce(G,\tau)\big]$. Clearly, the group scheme $f_{_*}(\ce(G,\tau))$ on $\tC$ is obtained by 
gluing two Bruhat-Tits group schemes ${\sf P}_{_{\theta_{_{\tau}}}}$, ${\sf P}_{_{\theta_{_{\bar\tau}}}}$ at $c_{_1}, c_{_2}$ (with $\bar\tau$ as in Remark \ref{balancedaction}) together with the datum of an isomorphism of the Levi factors of the closed fibres ${\sf P}_{_{c_{_i}}}$.

\bprop\label{bttorsorlacedtorsor} Let ${\qh}^{^G}_{_{\tau}}$ be a group scheme   over $C^{^{(d)}}$ coming from a representation $\rho:\sf\mu_{_d} \to G$  together with a gluing datum. Then the restriction 
\beqa\label{2tbala}
\eG_{_{(\theta_{_{\tau}})}} := {\qh}^{^G}_{_{\tau}} \bigg\vert _{_{\tC}}
\eeqa
to $\tC \subset C^{^{(d)}}$ is a balanced group scheme. \eprop

\bpr The question is obviously local along the components of $C^{^{(d)}}$. 
The process of taking invariant direct images of group schemes, i.e., Weil restriction followed by invariants, commutes with base change. Hence by \eqref{egtau} and \eqref{nsgrpscheme1.5} we get:
\beqa\label{justabove}
({\tt Inv} \circ f_{_*})\Big(E(G,\tau)\bigg\vert_{_{\tN}}\Big) \simeq f_{_*}\big(\ce(G,\tau)\big) \bigg\vert_{_{\tU}}.
\eeqa
The proof of \eqref{2tbala}  follows immediately from \eqref{nsgrpscheme1.5} and \eqref{justabove}. \epr

We are in the setting of \eqref{stackynormalization}.
\bdefe\label{2-lacedbundle} ("Fourier-Mukai") Let $\eG_{_{(\theta_{_{\tau}})}}$ be a balanced parahoric group scheme. A  $\eG_{_{(\theta_{_{\tau}})}}$-torsor $E_{_\wp}$ on $(\tC, {\bf c})$  is called {\tt laced} if  ${E_{_{\wp}}} \simeq f_{_*}(q^{^*}(\ce))$ for a $G$-torsor $\ce$ of local type $\tau$ on $\cc_{_d}$ for some $d$.
\edefe

\subsubsection{{\tt Twisted} Levi torsor associated to {\tt laced} $\eG_{_{(\theta_{_{\tau}})}}$-torsors}\label{admlev}
Let $E_{_\wp}$ be the laced torsor associated to $\ce$. Given a reduction of structure group $\ce_{_P}$ of $\ce$ on $\cc_{_d}$, by pulling back $\ce^{^{\varphi}}_{_H}$ \eqref{genericreductiondatum}, using $q:\tilde{\cc_{_d}} \to \cc_{_d}$, we obtain a balanced $H$-torsor $q^{^*}(\ce^{^{\varphi}}_{_H})$ on $\tilde{\cc_{_d}}$. We can then take the push-forward ("invariant push-forward") by the morphism $f:\tilde{\cc_{_d}} \to \tC$ and we get the {\tt laced} torsor 
\beqa\label{lacedhtorsor}
E_{_{\wp,H}}^{^{\varphi}} := f_{_*}\big(q^{^*}(\ce^{^{\varphi}}_{_H})\big)
\eeqa
which is a torsor under the invariant direct image group scheme with generic fibre $H$. We term this the {\tt laced} Levi torsor associated to the $P$-torsor $\ce_{_P}$.
Let $\chi:P \to \bg_{_m}$ be any anti-dominant character. This is indeed a character of $H$ as well. Given a reduction of structure group of $\ce_{_P}$ to $P$, to the twisted $H$-torsor $\ce^{^{\varphi}}_{_H}$ we get line bundles $\ce_{_H}^{^{\varphi}}(\chi)$ on $\cc_{_d}$. These define parabolic line bundles on $\tC$ with a balanced parabolic structure at $\{c_{_1}, c_{_2}\}$ as below:
\beqa
E_{_{\wp,H}}^{^{\varphi}}(\chi):=  f_{_*}\big({q^{^*}(\ce_{_H}^{^{\varphi}}(\chi))}\big)
\eeqa
\subsubsection{Equivariant degree of line bundles on $\cc_{_d}$ and semistability} We now have  a well-defined {\tt tf}-(semi)stability of $G$-torsors $\ce$ on $\cc_{_d}$ by combining \eqref{goodshow} by using \eqref{degreeoftorsoroncc} for ${\tt deg}_{_{\cc_{_d}}}$.
\bdefe\label{degreeoftorsoroncc} Let $L$ be a line bundle on the twisted curve $\cc_{_d}$.   Define the equivariant degree ${\tt deg}_{_{\cc_{_d}}}(L) := {\tt par.deg}_{_{\tC}}~f_{_*}\big(q^{^*}(L)\big)$.\edefe

\section{{\tt tf}-semistability via GIT}  Let $\eta:G \hra \text{GL}(W)$ be a faithful representation and let $T_{_W} \subset \text{GL}(W)$ be a maximal torus such that $\eta:T \hra T_{_W}$. Let $\ce$ be a $G$-torsor on $\cc_{_d}$. Then we get the associated $\text{GL}(W)$-torsor $\ce(\text{GL}(W))$ and associated vector bundle  $\ce(W)$. Likewise, if $E_{_\wp} = f_{_*}(q^{^*}(\ce))$, then  $\eta$ gives the associated {\tt laced} $\text{GL}(W)$-torsor $E_{_\wp}(\text{GL}(W)) := f_{_*}(q^{^*}(\ce(\text{GL}(W))))$  and similarly the associated {\tt laced} vector bundle $E_{_\wp}(W)$ on $(\tC, {\bf c})$. 
Given a $1$-PS $\lambda$, we get graded versions and isomorphisms:
\small
\beqa\label{leviidently1}
 f_{_*} \circ q^{^*}\Big[\ce^{^{\varphi}}_{_H}({\tt gr}(W))\Big] =  f_{_*} \circ q^{^*}\Big[\bigoplus_{_{j = 1}}^{^{m}} \ce^{^{\varphi}}_{_{H}}(W_{_{j+1}}/W_{_j})\Big]  \simeq \bigoplus_{_{j = 1}}^{^{m}} E^{^{\varphi}}_{_{\wp,H}} (W_{_{j+1}}/W_{_j}) = E^{^{\varphi}}_{_{\wp,H}}({\tt gr}(W))
\eeqa \normalsize
Weighted filtrations \eqref{wtdfiltongmodules0} on associated objects by an application of $f_{_*}(q^{^*})$ term by term give an obvious filtration for the associated graded:
\beqa\label{wtdfiltongmodules1}
0 \subsetneq E_{_{\wp,1}} \subsetneq \ldots E_{_{\wp,m}} \subsetneq E_{_{\wp,m+1}} = E^{^{\varphi}}_{_{\wp,H}}({\tt gr}(W))
\eeqa
\blem\label{formalequalityforpar} Let $\chi_{_{\lambda}}$ be as in \eqref{parabtrio}. Then we have the equality:
\beqa
{\tt par.deg}(E_{_{\wp,H}}^{^{\varphi}}(\chi_{_\lambda})) = L\big(E^{^{\varphi}}_{_{\wp,H}}({\tt gr}(W))\big).
\eeqa
In particular, \eqref{formalequality} holds for degree as defined in \eqref{degreeoftorsoroncc}.
\elem
\bpr This follows immediately from \eqref{theabovefortheparabcase}.\epr
We can get torsion-free sheaves on $C$ via $\sigma:\cc_{_d} \to C$ as well as via $\nu:\tC \to C$ by the functor $\nu_{_*}:{\text{Vect}}^{^{\tt L}}_{_{(\tC, {\bf c})}} \to \text{Tf}_{_C}$. Then, it is routine to check that:
\beqa\label{assoctfhere}
\nu_{_*}(E_{_\wp}(W)) \simeq \sigma_{_*}(\ce(W)) \simeq \cf_{_W} 
\eeqa
In fact, the parabolic degree of the laced vector bundle $E_{_\wp}(W)$ coincides with the degree of the torsion-free sheaf $\cf_{_W}$ \eqref{disting}. 
\blem\label{leviident} We have isomorphisms :
\beqa
{\tt gr}(\cf_{_W}) = \bigoplus_{_{j = 1}}^{^{m}} \cf_{_{j+1}}/\cf_{_j}  \simeq \bigoplus_{_{j = 1}}^{^{m}}\nu_{_*}\big[E^{^{\varphi}}_{_{\wp,H}} (W_{_{j+1}}/W_{_j})\big] =  \nu_{_*}\big({E_{_{\wp,H}}^{^{\varphi}}({\tt gr}(W))}\big) .
\eeqa
of torsion-free sheaves on $C$. Furthermore, for each $j$, 
\beqa
{\tt par.deg}_{_{\tC}}~\big[E^{^{\varphi}}_{_{\wp,H}} (W_{_{j+1}}/W_{_j})\big] = {\tt deg}_{_C}\big(\cf_{_{j+1}}/\cf_{_j}\big).
\eeqa
\elem
\bpr The first statement is a summary of the previous discussion. The last statement follows from \eqref{disting}.\epr

The next proposition is the primary reason why we go to the laced category. The balanced parabolic structure of laced bundles makes their parabolic degree coincide with the push down torsion free sheaf \eqref{disting}. 
\bprop\label{keyformoduli} We have the equality:
\beqa
L(\cf_{_{W}}^{^\bullet},\epsilon_{_\bullet}) = L\big(E^{^{\varphi}}_{_{\wp,H}}({\tt gr}(W))\big) \stackrel{\eqref{formalequalityforpar}}= {\tt par.deg}(E_{_{\wp,H}}^{^{\varphi}}(\chi_{_\lambda})).
\eeqa
where $L(\cf_{_{W}}^{^\bullet},\epsilon_{_\bullet})$ is as in \eqref{wtdslope}.  
\eprop
\bpr  By \eqref{leviident} applied to each summand, it follows easily that $$L(\cf_{_{W}}^{^\bullet},\epsilon_{_\bullet}) \stackrel{\eqref{trivialeq}} = L({\tt gr}_{_\lambda}(\cf_{_W})) = L\big(E^{^{\varphi}}_{_{\wp,H}}({\tt gr}(W))\big).$$ \epr
\bcor Let $\ce$ be a $G$-torsor on $\cc_{_d}$ and let $\ce_{_P}$ be a reduction of structure group to $P(\lambda) \subset G$. Let $\eta:G \hra \text{GL}(W)$.
Then we have the equality: 
\beqa
L(\cf_{_{W}}^{^\bullet},\epsilon_{_\bullet}) =  L\big({\ce^{^{\varphi}}_{_H}({\tt gr}(W)})\big) = {\tt deg}_{_{\cc_{_d}}}({\ce^{^{\varphi}}_{_H}}(\chi_{_\lambda})).
\eeqa
\ecor
\bpr By \eqref{leviidently1} $L\big({\ce^{^{\varphi}}_{_H}({\tt gr}(W)})\big) = L\big(E^{^{\varphi}}_{_{\wp,H}}({\tt gr}(W))\big)$  and the rest follows immediately.\epr

We recall from \cite{schmittnewstead}. The torsion-free sheaf $\cf_{_W}$ when restricted to $C-c$ (or more precisely, the underlying frame bundle), comes with a reduction of structure group $\upsilon$ to $G$. Further, since $G$ is semisimple ${\tt deg}_{_{C}}(\cf_{_W}) = 0$. Thus the pair $(\cf_{_W}, \upsilon)$ is a {\em singular principal bundle} (\ref{spb}). By \cite[3.5.1]{schmittnewstead} $(\cf_{_W}, \upsilon)$ is semi(stable) if $L(\cf_{_{W}}^{^\bullet},\epsilon_{_\bullet}) (\geq) 0$. This definition is a GIT notion and we arrive at the following:
\bth\label{etassandschmittss} ({\em {\tt tf}-semistability is a GIT notion}) Let $\eta:G \hra {\text{GL}(W)}$. A $G$-torsor $\ce$ of local type $\tau$ on $\cc_{_d}$  is {\tt tf}-semi(stable if and only if it is $\eta$-(semi)stable and hence  if and only if the associated singular principal bundle $(\cf_{_W}, \upsilon)$ \eqref{assoctfhere} is (semi)stable. \eeth
\brem The importance of this notion stems from the fact that the moduli space of semi(stable) singular principal bundles $\mathcal M_{_\eta}(C)$ is a projective scheme \cite[Theorem 3.5.2]{schmittnewstead}. Moreover, if $\eta$ is chosen carefully, there is a projective scheme $\mathcal M \to A$, {\em albeit non-flat}, such that $\mathcal M_{_K}$ is the Ramanathan moduli space over $C_{_K}$ and $\mathcal M_{_o} \simeq \mathcal M_{_\eta}(C)$. \erem

{\subsubsection{ Semistability of $G$-torsors on nodal curves}\label{allintrinsicstuff} Recall that if $E$ is a $G$-torsor on $(C,c)$ then $\sigma^{^*}(E)$ is a $G$-torsor on $\cc_{_d}$. The notion of {\tt tf}-(semi)stability of a $G$-torsor on $\cc_{_d}$ therefore gives an intrinsic notion of (semi)stability of $G$-torsors on $(C,c)$ thereby answering a long-standing question.

\subsubsection{Gieseker torsors on $C^{^{(d)}}$ and laced torsors on $\tC$} Let ${\text{Lac}}_{_{\tC}}(\eG_{_{(\theta)}})$ denote the set of isomorphism classes of  {\tt laced} $\eG_{_{(\theta_{_{\tau}})}}$-torsors on $(\tC, {\bf c})$.
\bprop\label{lacedtoGieseker}
We have a set-theoretic map from the isomorphism classes of Gieseker torsors on $(C,c)$ to the {\tt laced} $\eG_{_{(\theta_{_{\tau}})}}$-torsors on $\tC$ for varying types $\tau$:
\beqa\label{onthegiesekercurve}
|{\text{Gies}}_{_{G}}(C)| \to \bigsqcup_{_{\tau}} {\text{Lac}}_{_{\tC}}(\eG_{_{(\theta_{_{\tau}})}}).
\eeqa

\eprop
\bpr
It is immediate from \eqref{bttorsorlacedtorsor} that a Gieseker torsor on $(C,c)$ when restricted to $\tC$ gives a {\tt laced} $\eG_{_{(\theta_{_{\tau}})}}$-torsors on $\tC$. \epr

A laced torsor $E_{_\wp}$ \eqref{2-lacedbundle} is called {\tt tf}-semi(stable)  if the $G$-torsor $\ce$ of local type $\tau$ is so. Let $\bigsqcup_{_{\tau}} {\text{Lac}}_{_{\tC}}(\eG_{_{(\theta)}})^{^{{\tt tf}-ss}}$ be the subset of {\tt tf}-semi(stable) laced torsors.
\bdefe\label{ssBTtorsors}  A Gieseker torsor on the semistable curve $C^{^{(d)}}$ (for some $d > 0$) is called {\tt tf}-(semi)stable if its image under \eqref{onthegiesekercurve} lies in  $\bigsqcup_{_{\tau}}{\text{Lac}}_{_{\tC}}(\eG_{_{(\theta)}})^{^{{\tt tf}-ss}}$.

We have an obvious notion of families of {\tt tf}-(semi)stable Gieseker torsors.
Let 
${\text{Gies}}_{_G}(C_{_A})^{^{{\tt tf}-ss}}$
denote the substack of {\tt tf}-semi(stable) Gieseker torsors on $C_{_{A}}$.\edefe
For the openness property of this notion of 
{\tt tf}-semistability, see \eqref{opennessofss} below.

\begin{center}
\underline{\sc Part III}
\end{center}
\section{The Moduli construction}
The aim of this part is to combine the results of Parts I and II to construct flat degenerations of the Ramanathan moduli space of slope (semi)stable $G$-torsors. When $G = \text{\text{GL}(n)}$, these are the precise analogues of  degenerations via Gieseker bundles on modifications of the nodal curve (\cite{gies}, \cite{ns2}, \cite{kausz}).

\subsubsection{ The Bhosle-Schmitt spaces and associated Gieseker bundles}  Fix a faithful representation $G \hra \text{GL}(W)$ and let $\eta:G \hra \text{GL}(W \oplus W^{*}).$ Let $2w:= \text{ {\tt dim}}(W \oplus W^*)$. In \cite{schmitt2} and \cite{schmitt3}, Schmitt studies the algebraic ${A}$-stack  ${\text{Bun}}_{_{G}}^{^{\eta, \text{Sing}}}(C_{_A})$ whose generic fibre is the stack ${\text{Bun}_{_{C_{_K}}}}(G)$ of $G$-torsors on $C_{_K}$ and whose closed fibre has  $T$-points which are  families of {\em singular principal $G$-bundles} i.e. of pairs $(\cf, \text{\cursive s})$, 
\begin{itemize}\label{spb}
\item A torsion-free $\co_{_{C}}$-module $\cf$ with generic fibre type $W \oplus W^*$.
\item A pseudo-$G$-structure $\text{\cursive s}$ which gives a  reduction of structure group   of the principal $\text{GL}(W \oplus W^*)$-bundle on $\tC^* = \tC - \{c_{_1}, c_{_2} \}$ underlying the locally free sheaf $\cf \Big\vert_{_{\tC^*}}$ to the subgroup $\eta: G \hra \text{GL}(W 
\oplus W^{*})$. 
\end{itemize}
Note that, $G$-torsors on the generic fibre $C_{_K}$ are also viewed as a pair $(\cf_{_K}, \text{\cursive s}_{_K})$ where $\cf_{_K}$ is a locally free sheaf of rank $2w$ on $C_{_K}$ and $\text{\cursive s}_{_K}$ a reduction of structure group of the frame $\text{GL}(2w)$-bundle of $\cf_{_K}$ to $G$.
In particular, there is a forget ${A}$-morphism 
\beqa
{\small\text{\cursive j}}: {\text{Bun}}_{_{G}}^{^{\eta, \text{Sing}}}(C_{_A}) \to \text{Tfs}_{_{2w}}(C_{_{A}}) ~~~~~,
(\cf, \text{\cursive s}) \mapsto \cf
\eeqa 
into the algebraic stack $\text{Tfs}_{_{2w}}(C_{_{A}})$ of relative torsion-free sheaves of rank $2w$ on the surface $C_{_{A}}$. On the other hand, there is the morphism of stacks over ${A}$:
\beqa\label{nagarajseshadri}
{\tt p}_{_*}: \text{GVB}_{_{2w}}(C_{_{A}}) \to \text{Tfs}_{_{2w}}(C_{_{A}})
\eeqa
which is obtained by taking direct images under the canonical morphism ${\tt M} \to C_{_T}$ for varying ${A}$-schemes $T$. Thus we have a fibre square:
\beqa\label{bhosleschmittgieseker1}
\xymatrix{
{\small\text{\cursive j}}^{^*}(\text{GVB}_{_{2w}}(C_{_{A}})) \ar[r]^{} \ar[d]_{{\tt p}_{_*}} &
\text{GVB}_{_{2w}}(C_{_{A}}) \ar[d]^{{\tt p}_{_*}} \\
{\text{Bun}}_{_{G}}^{^{\eta, \text{Sing}}}(C_{_A}) {\ar[r]_{_{\small\text{\cursive j}}}} &  \text{Tfs}_{_{2w}}(C_{_{A}})
}
\eeqa
The stack  ${\small\text{\cursive j}}^{^*}(\text{GVB}_{_{2w}}(C_{_{A}}))$ parametrizes pairs $(V, \text{\cursive s})$ of Gieseker vector bundles of rank $2w$ and a generic reduction of structure group $\text{\cursive s}$ of the underlying principal $\text{GL}(2w)$-bundle $V_{_{\text{GL}}}$ to the subgroup $G$ via $\eta:G \hra \text{GL}(2w)$.

By combining this diagram with the isomorphism (\ref{balequalskausz}) and the morphism $\eta_{_*}$ from \eqref{heeheeeee} (obtained via extension of structure groups), we get a commutative diagram:
\beqa\label{bhosleschmittgieseker2} 
\xymatrix{
{\text{Gies}}_{_G}(C_{_A}) \ar[r]^{\eta_{_*}} \ar[d]_{} &
{\text{Gies}}_{_{{\text{GL}(2w)}}}(C_{_A}) \ar[d]^{\simeq, (\ref{balequalskausz})} \\
{\small\text{\cursive j}}^{^*}(\text{GVB}_{_{2w}}(C_{_{A}})) \ar[r]^{f} \ar[d]_{{\tt p}_{_*}} &
\text{GVB}_{_{2w}}(C_{_{A}}) \ar[d]^{{\tt p}_{_*}} \\
{\text{Bun}}_{_{G}}^{^{\eta, \text{Sing}}}(C_{_A}) {\ar[r]_{_{\small\text{\cursive j}}}} &  \text{Tfs}_{_{2w}}(C_{_{A}})}
\eeqa
Loosely speaking, we may view the stack ${\text{Gies}}_{_G}(C_{_A})$ as parametrizing pairs $(\ce, \text{\cursive s})$ with $\ce \in {\text{Gies}}_{_{C_{_{A}}}}({\text{GL}(2w)})(T)$ together with a reduction of structure group $\text{\cursive s}$. 

Thus, by composition with the direct image, we get a morphism ${\text{Gies}}_{_G}(C_{_A})  \to \text{Tfs}_{_{2w}}(C_{_{A}}).$ This morphism factors via a vertical morphism  ${\text{Gies}}_{_G}(C_{_A})  \to {\text{Bun}}_{_{G}}^{^{\eta, \text{Sing}}}(C_{_A})$. To see this, we need to note the following on the torsion-free sheaf $p_{_*}(\ce)$. The reduction of structure group $\text{\cursive s}$ to $G$ away from the singularity on the normal surface $C_{_{A}}$ extends to give a pseudo $G$-bundle $(p_{_*}(\ce), \tau)$ (with notation as in \cite[Page 1428, Section 1.1]{schmitt2}). This follows by the normality of the surface and a simple Hartogs argument (see for example \cite[Remark 2.5]{balajidu}). 

On the other hand, the image under the unlabelled vertical morphism from ${\text{Gies}}_{_G}(C_{_A})  \to {\small\text{\cursive j}}^{^*}(\text{GVB}_{_{2w}}(C_{_{A}}))$ consists of pairs   $(V_{_{\text{GL}}}, \text{\cursive s}')$ such that the $G$-torsor given by the generic reduction of structure group $\text{\cursive s}'$ extends to a full Gieseker torsor over ${\tt M}$ for a group algebraic space   ${\qh}^{^G}_{_{\mathfrak t(T), \tt M}}$ (which extends the semi-simple group scheme with fibre $G$). We denote the composite vertical ${A}$-morphism by:
\beqa
\text{\cursive p}_{_{G}}: {\text{Gies}}_{_G}(C_{_A})  \to {\text{Bun}}{_{_G}}^{^{\eta, \text{Sing}}}(C_{_A}).
\eeqa
which is  analogous to taking direct images in the case of locally free sheaves.
\subsubsection{A properness result} We recall the following definitions from \cite[Page 15]{bbn}.
\bdefe\label{horizontal properness} {\em (Horizontal properness)} Let $F,G: \text{Sch}_{_S} \to \text{Sets}$ be two functors with $S = \spec~A$ for a discrete valuation ring $A$ and quotient field $K$. Let $f:F \to G$ be a $S$-morphism. We say, $f$ is {\em horizontally proper} if  the following property holds:
let $B$ be a discrete valuation ring with function field $L$ such that $L$ is a finite extension of $K$ and $\spec~ B \to \spec~ A$ is surjective. Then for every map $\alpha \in F(L)$, if the composite $f(\alpha) \in G(L)$ extends to an element $G(A)$, then $\alpha$ also extends to an element in $F(A)$. \edefe 
This definition becomes significant because of the following observation.

\blem\label{quasiprojmorphs} Let $f:F \to G$ be a projective $S$-morphism of schemes of finite type such that $f_{_\zeta}:{F}_{_\zeta} \to {G}_{_\zeta}$ over the generic point is proper. Suppose further that the structure morphisms $F \to S$ and $G \to S$ are surjective, that $F$ is $S$-flat and that $f$ is horizontally proper. Then $f$ is proper. \elem
Using a stability parameter $\delta$, Schmitt (\cite[page 340]{schmittbook}) defines an open substack ${\text{Bun}}_{_C}^{^{\eta, \text{Sing}}}(G)^{^{ss}}$ of $\delta$-(semi)stable pairs $(\cf, \text{\cursive s})$ and by using GIT (of decorated objects), he then goes on to construct the coarse moduli space $\cm_{_C}^{^{\eta, \text{Sing}}}(G)^{^{ss}}$ of this stack (as a subspace of a certain "big'' moduli space constructed by Bhosle \cite{bhosle1}).

When $\delta$ is chosen "large'', Schmitt, in \cite[Theorem 1.1]{schmitt2}, shows that the generic fibre $\cm_{_{C_{_K}}}^{^{\eta, \text{Sing}}}(G)^{^{ss}}$ is independent of the faithful representation $\eta$ and is in fact the moduli space of $S$-equivalence classes of slope (semi)stable $G$-torsors (in the sense of Ramanathan) and the special fibre is the moduli space $\cm_{_C}^{^{\eta, \text{Sing}}}(G)^{^{ss}}$ of semistable singular principal bundles. A small check using \eqref{etassandschmittss} shows that 
\beqa\label{veryclose}
\text{\cursive p}_{_{G}}^{^\ast}\big( {\text{Bun}}_{_{G}}^{^{\eta, \text{Sing}}}(C_{_A})^{^{ss}}\big) \simeq {\text{Gies}}_{_G}(C_{_A}) ^{^{{\tt tf}-ss}}
\eeqa 
\brem\label{opennessofss} The identification \eqref{veryclose} in particular shows that the substack ${\text{Gies}}_{_G}(C_{_A}) ^{^{{\tt tf}-ss}}$  of  ${\text{Gies}}_{_G}(C_{_A}) $  of {\tt tf}-(semi)stable Gieseker torsors on $C_{_{A}}$  is an {\em open} substack, thereby verifying the openness of the notion of {\tt tf}-semistability as defined in (\ref{ssBTtorsors}).\erem

\bth\label{verticalproper} The  ${A}$-morphism:
\beqa\label{tobhossch}
\text{\cursive p}_{_{G}}:{\text{Gies}}_{_G}(C_{_A}) ^{^{{\tt tf}-ss}} \to {\text{Bun}}_{_{G}}^{^{\eta, \text{Sing}}}(C_{_A})^{^{ss}}
\eeqa
is {\em horizontally} proper and an isomorphism over $C_{_K}.$ 
Over the closed point $a \in {A}$ it induces  a morphism:
\beqa\label{tobhossch1}
\text{\cursive p}:{\text{Gies}}_{_{G}}(C)^{^{{\tt tf}-ss}}
 \to {\text{Bun}}_{_G}^{^{\eta, \text{Sing}}}(C)^{^{ss}}.
\eeqa
\eeth
\bpr  Let $A' = \spec~k\llbracket z \rrbracket$ with $z^{^{d}} = t$  let $L = \spec~k((z))$ be the function field of $A'$. Let $E_{_L}$ be a family of semistable $G$-bundles on the smooth curve $C_{_L}$ degenerating to a semistable singular principal $G$-bundle $E_{_o}$ on $(C,c)$. Equivalently, on the nodal curve $(C, c)$, there exists a pair $(\cf, \text{\cursive s})$, with a torsion-free $\co_{_C}$-module $\cf$ of rank $2w$ together with a  reduction of structure group $\text{\cursive s}$ on $C-c$ to $G$ such that the family $E_{_L}$ degenerates to $(\cf, \text{\cursive s})$ which is $\delta$-semistable. By definition, the reduction of structure group $\text{\cursive s}$ gives rise to a $G$-torsor $E_{_{A' -c}}$ on $C_{_{A'}} - c$. The local type of the torsion-free sheaf gives also a diagram as in \eqref{keydiag1}. 
 
The diagram is  also such that the pull-back $\sigma_{_{A'}}^{^*}(E_{_{A' -c}})$ gives a $({\sf\mu_{_d}},G)$-torsor $P$ on $D-o$ (since the frame  $\text{GL}(2w)$-bundle of $\cf\mid_{_{C-c}}$ comes with a reduction of structure group to $G$ in the complement of $c$ in ${\sf N}^{^{(d)}}$). By the smoothness of $D$, and an application of Hartogs theorem, or by a theorem of Colliot-Th\'el\`ene and
Sansuc (\cite[Theorem 6.7]{collio}), it follows that we get an extension of $P$ to a $({\sf\mu_{_d}},G)$-torsor $\overline{P}$ on  ${D}$. This  is of some local type $\tau$.

Therefore, $q^{^*}(\overline{P})$ gives a $E(G, \tau)$-torsor on $D^{^{(d)}}$ and by taking invariant direct images  we get $\eP = \text{Inv}\circ f_{_*}\big(q^{^*}(\overline{P})\big)$ such that $({\qh^{^G}_{_{\tau, {\sf N}^{^{(d)}}}}}, \eP)$ is an admissible pair on ${\sf N}^{^{(d)}}$. This can be glued  to the $G$-torsor pulled back from $C_{_{A'}} - c$ to get an admissible pair $({\qh^{^G}_{_{\tau,C_{_{A'}}^{^{(d)}}}}},\eP_{_{A'}})$ on $C_{_{A'}}^{^{(d)}}$.

By definition  the parabolically associated vector bundle ${\eP}^{^{par}}_{_{A'}}(W \oplus W^*)$ via the homomorphism $\eta:G \to \text{GL}(W \oplus W^*)$ is a quasi-admissible vector bundle on $C_{_{A'}}^{^{(d)}}$ and furthermore, we have:
\beqa
p_{_*}\big ({\eP}_{_o}(W \oplus W^*)\big) = \cf
\eeqa
Since $(\cf, \text{\cursive s})$ is a semistable singular bundle, by \eqref{etassandschmittss}, this implies that the $\eP_{_o}$ is {\tt tf}-semistable. Clearly, $\eP_{_L} \simeq E_{_L}$ and hence, the family $({\qh^{^G}_{_{\tau,C_{_{A'}}^{^{(d)}}}}},\eP_{_{A'}})$ gives the required point in   ${\text{Gies}}_{_G}(C_{_A}) ^{^{{\tt tf}-ss}}(A')$ proving the horizontal properness.\epr
\subsubsection{The coarse moduli}
Let ${\text{Bun}}_{_{G}}^{^{\eta, \text{Sing}}}(C_{_A})^{^{ss}} \to \cm_{_{_{C_{_{A}}}}}^{^{\eta, \text{Sing}}}(G)^{^{ss}} $
be the canonical morphism to the coarse space. Let $R_{_{_{C_{_{A}}}}}^{^{\eta}}(G)$ denote the total family so that the coarse space 
\beqa
\cm_{_{_{C_{_{A}}}}}^{^{\eta, \text{Sing}}}(G)^{^{ss}} = R_{_{_{C_{_{A}}}}}^{^{\eta}}(G) \parallelslant \text{PGL}_{_{A}}
\eeqa  
is realized as a GIT quotient of $R_{_{_{C_{_{A}}}}}^{^{\eta}}(G)$ by a suitable reductive group $\text{PGL}_{_{A}}$.

We recall the definition of the functor ${\eG}^{^{G}}_{_{{A}}}$ as also the $\mathtt \tC_{_{}}$-scheme $\mathtt \tC^{^G}$ which represents it (see (\ref{yh})). The diagram \eqref{bhosleschmittgieseker2} gives a diagram of total families:
\beqa
\xymatrix{
\mathtt \tC^{^G} \ar[r]^{} \ar[d]_{\text{\cursive p}_{_{G}}} &
\mathtt \tC_{_{}} \ar[d]^{p_{_*}} \\
R_{_{_{C_{_{A}}}}}^{^{\eta}}(G) {\ar[r]_{_{\small\text{\cursive j}}}} &  R^{^{tf}}_{_{C_{_{A}}}}
}
\eeqa
where the morphism $\text{\cursive p}_{_{G}}$ is {\em horizontally} proper (\ref{verticalproper}). Further,  since the objects are projective, by (\ref{quasiprojmorphs}, \ref{flatavecnc}) it follows that  $\text{\cursive p}_{_{G}}$ is  proper.

At the level of total families we also have the identification:
\beqa
(\mathtt \tC^{^G})^{^{{\tt tf}-ss}} = \text{\cursive p}_{_{G}}^{^*} \big(R_{_{_{C_{_{A}}}}}^{^{\eta}}(G)^{^{ss}} \big)
\eeqa
as an open subscheme of ${\mathtt \tC}^{^G}$.

We are now in the precise setting of the theme in Nagaraj-Seshadri \cite[Page 180 and Remark 6, page 180, 184]{ns2}. We consider the polarization 
\beqa
{\tt L}_{_{\epsilon}} := \text{\cursive p}_{_{G}}^{^*}\Big(\co_{_{R_{_{_{C_{_{A}}}}}^{^{\eta}}(G)^{^{ss}}}}(1)\Big) \otimes \Big(\co_{_{\mathtt \tC^{^G}}}(\epsilon)\Big)
\eeqa 
where $\co_{_{\mathtt \tC^{^G}}}(1)$ is the relative polarization for $\text{\cursive p}_{_{G}}$. By choosing a "small'' $\epsilon$ for the relative polarization,  we get a natural notion of (semi)stability (which we term {\tt L}-(semi)stability), on $(\mathtt \tC^{^G})^{^{{\tt tf}-ss}}$(and hence on ${\text{Gies}}_{_G}(C_{_A})^{^{{\tt tf}-ss}}$) which is such that {\tt L}-(semi)stability implies {\tt tf}-(semi)stability and via \eqref{veryclose} we get an inclusion:
\beqa\label{elss}
{\text{Gies}}_{_G}(C_{_A})^{^{{\tt L}-ss}} \subsetneq   {\text{Gies}}_{_G}(C_{_A})^{^{{\tt tf}-ss}}
\eeqa
(an inclusion which is a proper one in general (see \cite[page 185]{ns2})). In fact, by GIT, the {\tt L}-(semi)stability constructs the actual separated coarse space $\cm_{_{C_{_{A}}}}^{^{{\tt L}-ss}}({\qh^{^G}_{_{A}}}) = (\mathtt \tC^{^G})^{^{{\tt L}-ss}}\parallelslant \text{PGL}_{_{A}}$, for the  Artin stack ${\text{Gies}}_{_G}(C_{_A})^{^{{\tt L}-ss}}$

We summarise this discussion, using (\ref{flatavecnc}) and by following the arguments in \cite{ns2}, to arrive at the following main theorem: 
\bth\label{the big theorem} \begin{enumerate} 
\item The  stack ${\text{Gies}}_{_G}(C_{_A})^{^{{\tt L}-ss}} \subset {\text{Gies}}_{_G}(C_{_A})$ is an algebraic stack, which is locally of finite type and flat over ${A}$. 
\item The generic fibre ${\text{Gies}}_{_{G}}(C_{_K})^{^{{\tt L}-ss}}$ is isomorphic to the algebraic stack $\text{Bun}_{_{G}}(C_{_K})^{^{\mu-ss}}$ of $\mu$-(semi)stable $G$-torsors on the smooth projective curve $C_{_K}$ and the closed fibre ${\text{Gies}}_{_{G}}(C)^{^{{\tt L}-ss}}$ is a divisor with normal crossings.
\item The closed fibre has an open subscheme comprising of (semi)stable $G$-torsors on the nodal curve $(C,c)$, where (semi)stability is the intrinsic one from \eqref{etassandschmittss}. 
\item The coarse space $\cm_{_{C_{_{A}}}}^{^{{\tt L}-ss}}({\qh^{^G}_{_{A}}}) = (\mathtt Y^{^G})^{^{{\tt L}-ss}}\parallelslant \text{PGL}_{_{A}}$, for the Artin stack ${\text{Gies}}_{_G}(C_{_A})^{^{{\tt L}-ss}}$ therefore provides a proper and flat degeneration of the  moduli space of $\mu$-(semi)stable $G$-torsors on smooth curves degenerating to a simple nodal curve.

\item There is a morphism of coarse spaces:
\beqa
\tt p_{_{*}}: \cm_{_{C_{_{A}}}}^{^{{\tt L}-ss}}({\qh^{^G}_{_{A}}}) \to \cm_{_{_{C_{_{A}}}}}^{^{\eta, \text{Sing}}}(G)^{^{ss}}
\eeqa
The closed fibre $\cm_{_{C}}^{^{{\tt L}-ss}}({\qh^{^G}_{_a}})$ of the coarse moduli scheme $\cm_{_{C_{_{A}}}}^{^{{\tt L}-ss}}({\qh^{^G}_{_{A}}})$ parametrizes $S$-equivalence classes of {\tt L}-semistable Gieseker torsors. This scheme contains an open dense subscheme of {\tt tf}-semistable $G$-torsors on the underlying nodal curve $(C,c)$.
\end{enumerate} 
\eeth

\brem Recall that the moduli scheme $\cm_{_{_{C_{_{A}}}}}^{^{\eta, \text{Sing}}}(G)^{^{ss}} \to S$ provides a degeneration of the moduli spaces $\cm_{_{C_{_K}}}(G)$, but the drawback with this construction is that  $\cm_{_{_{C_{_{A}}}}}^{^{\eta, \text{Sing}}}(G)^{^{ss}}$ is {\em not} ${A}$-flat. \erem
\subsection{The orthogonal and symplectic case}\label{faltingsstuff} Recall that Faltings \cite{faltings} has constructed the moduli space of semistable orthogonal and symplectic torsion-free sheaves on nodal curves and gets a flat degeneration of the moduli space of semistable orthogonal and symplectic bundles on smooth curves when the curves degenerates to a simple nodal curve. By the comments in  \cite[Page 1430, 1436]{schmitt2}, we see that under the faithful representation $\eta:G \hra \text{GL}(W \oplus W^{*})$ the image lies in the orthogonal (resp. symplectic) group $O(W \oplus W^{*},q)$ (resp. $Sp(W \oplus W^{*},q)$) where $W \oplus W^{*}$ is seen to be equipped with a canonical non-degenerate symmetric (resp. alternating) form $q$. One could more generally have worked with any pair $(W,q)$ where  $W$ is equipped with a non-degenerate symmetric (resp. alternating) form $q$ and carried out the entire construction of the coarse spaces $\cm_{_{C_{_{A}}}}^{^{{\tt L}-ss}}({\qh^{^G}_{_{A}}})$. It is now easy to conclude from the main results that in the case when $G$ is either orthogonal or symplectic, then  the moduli space $\cm_{_{_{C_{_{A}}}}}^{^{\eta, \text{Sing}}}(G)^{^{ss}}$ is the moduli space of \cite{faltings} and the  morphism  $\tt p_{_{*}}: \cm_{_{C_{_{A}}}}^{^{{\tt L}-ss}}({\qh^{^G}_{_{A}}}) \to \cm_{_{_{C_{_{A}}}}}^{^{\eta, \text{Sing}}}(G)^{^{ss}}$ is a surjection.

\section{Appendix to Part I}

\subsubsection{Quasi-Gieseker bundles} In the first part of the appendix we will outline a small variant of the theme developed in \cite{gies}, \cite{ns1} and \cite{kausz}.   

We recall the notion of an admissible vector bundle $V$ on a curve $C^{^{(d)}}$ 
(\cite[Definition 3.11]{kausz}, \cite[Definition 3.6]{bbn}) and add a variant, namely the notion of a {\em quasi-admissible bundle}. In fact, Kiem and Li in \cite[Lemma 1.2(a)]{kiem-li} just call these {\em admissible bundles}. In \cite[Definition 1, page 167]{ns1}, we have the notion of a {\em standard} vector bundle on $C^{^{(d)}}$ as a preliminary notion.
\bdefe\label{Giesekervectorbundleappendix} Let $V$ be a vector bundle  of rank $n$ on a chain $E^{^{(d)}}$. Let $V\mid_{_{R_i}} = \oplus_{j=1}^{n} \co(a_{ij})$, where the $R_{_j}$ are the $\bp^{^1}$'s on the chain $E^{^{(d)}}$. Say that $V$ is {\em standard} if the $a_{ij}$ are $0$ or $1$.  The bundle $V$ is called  {\em strictly standard}  if moreover, for every $i$ there is an index $j$ such that $a_{ij} = 1$.

 A vector bundle $V$ on $C^{^{(d)}}$ of rank $n$ is called {\em admissible} (resp. {\em quasi-admissible}), if, for $d \geq 1$, the restriction $V\lvert{_{_{E^{^{(d)}}}}}$ is {\em strictly standard} (resp. {\em standard}) and  
the direct image $(p_{_d})_{_*}(V)$ is a torsion-free  $\co_{_C}$-module, where $p_{_d}:C^{^{(d)}} \to C$ is the canonical morphism which contracts the chain to the node. \edefe

{\sl The notion of admissibility (resp. quasi-admissibility) extends obviously to vector bundles on any modification $\tt M$ over $T \in \text{Sch} \big/{A}.$} Let $V$ be a standard vector bundle on $C^{^{(d)}}$ of $\text{rank}(V) = n$. Then, by the discussions in \cite[Page 168-171]{ns1}, after twisting the vector bundles sufficiently to ensure the vanishing of the first cohomology and ensure generation by sections, {\bf we get a canonical  morphism} $\phi_{_V}:C^{^{(d)}} \to \text{Grass}(\text{H}^{^0}(V), n)$. This morphism contracts the $R_{_j} = \bp^{^1}$'s on the chain $E^{^{(d)}} \subset C^{^{(d)}}$ such that the restriction $V {{\mid_{_{R_{_j}}}}}$ is trivial. The condition that $V$ is strictly standard is shown to be equivalent to the morphism $\phi_{_V}$ being a {\em closed immersion}. 

Let $N := {\tt dim}(\text{H}^{^0}(V))$. Let $W[j]$  be the $j$-th standard model with  $C^{^{(j)}} \subset W[j]$ as the central fibre  \eqref{expandeddeg}. Recall that this a smooth quasi-projective scheme with a tautological morphism $W[j] \to C \times_{_{A}} B[j]$.
For each $j \leq n$, fix the coordinate plane embedding $\ba^{^{j+1}} \subset \ba^{^{n+1}}$ by the first coordinates. This gives an identification $W[j] \simeq W[n]\times_{_{B[n]}} B[j]$ compatible with the tautological morphism (\cite[page 526]{Li}).  Define
\beqa\label{zN}  
W(N,n) :=  W[n] \times_{_{A}} \text{Grass}(N, n)
\eeqa
If $V$ is a standard bundle on $C^{^{(j)}}$ for some $j$ we get a closed immersion:
\beqa
C^{^{(j)}} \subset W(N,n)
\eeqa
via the inclusions $\text{Graph}(\phi_{_V}) \subset C^{^{(j)}} \times \text{Grass}(N, n) \subset W(N,n)$.

Following \cite[page 179]{gies} and \cite[Definition 7, page 185]{ns1} we have the definition.
\bdefe\label{GiesekerfunctorreltoR}  
Let  ${\eG}^{^{\mathtt q}}_{_{N}}:\text{Sch}_{_{A}} \to \text{Sets}$, be the functor defined as: \small
\beqa\label{as in GiesekerfunctorreltoR}
{\eG}^{^{{\mathtt q}}}_{_{N}}(T) = (\tt M, \text{\cursive e})
\eeqa
\normalsize
where 
\beqa\label{cursivee}
\text{\cursive e}: {\tt M} \hra C_{_{A}} \times_{_{A}} T\times_{_k} W(N,n)
\eeqa 
is a  closed embedding in the product and  such that,
(a) the projection $j:\tt M \to T\times_{_k} W(N,n)$ is a closed immersion,(b) the projection $\pi:\tt M \to C \times_{_{A}} T$ is  a {modification} as in Definition \ref{modification}, and (c) the projection $q_{_T}:\tt M \to T$ is a flat family of curves $\tt M$, $t \in T$ as in Definition \ref{modification}. (d) Moreover, {\tt the chain lengths $d$ occuring in $\tt M$ is bounded above by} $n$.

Further, if $V$ is the tautological quotient bundle of rank $n$  on $\text{Grass}(N,n)$ and $V_{_T}$ its pull-back to $T\times \ W(N,n)$, then the pull-back  $\eV_{_T}:= j^*(V_{_T})$ is such that, $\eV_{_T}$ is a {\emph {quasi-admissible}} vector bundle of rank $n$\eqref{Giesekervectorbundleappendix}  for the modification $\tt M \to C \times_{_{A}} T$. 
\edefe
By the definition of $\eV_{_T}$, for each $t \in T$ we get  a quotient morphism $\co^{^N}_{_{\tt M}} \twoheadrightarrow \eV_{_t}$, and  we assume that this map induces an isomorphism: $H^{^0}(\co^{^N}_{_{\tt M}}) \simeq H^{^0}(\eV_{_t})$. In particular, we have $ {\tt dim}(H^{^0}(\eV_{_t})) = N$ and it follows that
\beqa\label{h1-vanish}
H^{^1}(\eV_{_t}) = 0
\eeqa
As  in \cite{gies} and \cite[Proposition 8]{ns2}, it is easily seen that this new functor ${\eG}^{^{{\mathtt q}}}_{_{{N}}}$ is also represented by a $\text{PGL}(N)$-invariant open subscheme $\mathtt Y$ of the Hilbert scheme $\text{Hilb}(C_{_{A}} \times_{_{A}} W(N,n))$ for the natural polarization on $W(N,n)$. 

Let $\mathbb M \subset C_{_{A}} \times_{_{A}} \mathtt Y\times_{_k} \ W(N,n)$ be the {\em universal object} defining the functor ${\eG}^{^{\mathtt q}}_{_{{N}}}$. This defines a {\em universal modification}  $\mathbb M \to  \mathtt Y$ together with a {\em universal quasi-admissible vector bundle} $\mathbb V$ on $\mathbb M$. The representability of the functor ${\eG}^{^{\mathtt q}}_{_{{N}}}$ implies that for any quasi-admissible vector bundle $V$ on a modification $\tt M_{_{T}}$ there exists a unique morphism $\psi:T \to 
\mathtt Y$ and  $\phi:\tt M_{_{T}} \to \mathbb M$ so that $\phi^{^{*}}(\mathbb V)$ is $V$.

\noindent 
{\ul{The stack}} $\text{GVB}^{^{\mathtt q}}_{_n}(C_{_{A}})$ (cf. \cite[Definition 3.11]{kausz}):  $(\tt M, \mathcal V) \in \text{GVB}^{^{\mathtt q}}_{_n}(C_{_{A}})(T)$ is such that (1) $\mathcal V$ is a quasi-admissible vector bundle on the modification $\tt M$ and (2) $d \leq n$ for chains $E^{^{(d)}}$ in $\tt M$ . We may call $(\tt M, \mathcal V)$ a {\em quasi-Gieseker bundle}. Modifications with bounded chain lengths is easily seen to be a stack and $\text{GVB}^{^{\mathtt q}}_{_n}(C_{_{A}})$ is easily checked to be an Artin stack.

As in \cite[Definition 3.22]{kausz}, if we fix a very ample sheaf on $C$. Then for a quasi-Gieseker vector bundle $(\tt M, \cv)$ for $T$ a $A$-scheme and for an integer $N'$ we have the quasi-admissible bundle $\cv(N')$ and for every pair of integers $N \geq n, N' \geq 0$, we have a canonical morphism of $A$-groupoids:
\beqa
{\eG}^{^{\mathtt q}}_{_{N}} \to \text{GVB}^{^{\mathtt q}}_{_n}(C_{_{A}})
\eeqa
Analogous to \cite[Lemma 3.23]{kausz}, given a quasi-Gieseker bundle $({\tt M} \to C \times_{_{A}} T, \eV)$, we again have an open subschemes $T_{_{N,N'}} \subset T$ which has properties (1) and (2) in \cite[Lemma 3.23]{kausz}, with the added observation that the scheme $W(N,n)$, which replaces the Grassmannian in {\em loc cit}, ensures that $\forall t \in T_{_{N,N'}}$, the induced morphism ${\tt M_{_t}} \to W(N,n)$ is a closed immersion. 

For the analogue of \cite[Proposition 3.24]{kausz}, we need to do a bit more.
\bprop\label{prop3.24} The morphism of $A$-groupoids:
\beqa
\coprod_{_{N \geq n, N'\geq 0}} {\eG}^{^{\mathtt q}}_{_{N}} \to  \text{GVB}^{^{\mathtt q}}_{_n}(C_{_{A}})
\eeqa
is smooth and surjective.
\eprop
\bpr Let $T$ be a $A$-scheme and let $T \to \text{GVB}^{^{\mathtt q}}_{_n}(C_{_{A}})$ a $T$-point on $\text{GVB}^{^{\mathtt q}}_{_n}(C_{_{A}})$ given by a quasi-Gieseker bundle $({\tt M} \to C \times_{_{A}} T, \eV)$. Let $Z$ be the $A$-groupoid defined by the cartesian square:
\beqa\label{atlasshrugged}
\xymatrix{
Z \ar[r] \ar[d] &
{\eG}^{^{\mathtt q}}_{_{N}} \ar[d] \\
T \ar[r] &  \text{GVB}^{^{\mathtt q}}_{_n}(C_{_{A}})
}
\eeqa
Let $\{T_{_\alpha}\}$ be en \'etale cover of $T$ so that we have a morphism $T_{_\alpha} \to B[d_{_\alpha}]$ and modification $\tt M_{_\alpha}$ comes as a pull-back. For each quasi-Gieseker bundle $(\tt M_{_\alpha}, \eV_{_\alpha})$, we again have open subschemes $T_{_{N,N',\alpha}} \subset T_{_\alpha}$ with properties as stated above. We in fact have a morphism $\tt M \mid _{_{T_{_{N,N',\alpha}}}} \longrightarrow T_{_{N,N',\alpha}} \times W[d_{_\alpha}] \times \text{Grass}(N,n)$ and hence a morphism $\tt M \mid _{_{T_{_{N,N',\alpha}}}} \longrightarrow T_{_{N,N',\alpha}} \times W(N,n)$. This morphism is proper and for each $\forall t \in T_{_{N,N',\alpha}}$, the induced morphism $\tt M_{_t} \to W(N,n)$ is a closed immersion. Hence by \cite[Lemma 3.13]{kausz}, we get a closed immersion $\tt M \mid _{_{T_{_{N,N',\alpha}}}} \hra T_{_{N,N',\alpha}} \times W(N,n)$.

Let $Z_{_\alpha} = T_{_\alpha} \times_{_T} Z$.Then, following the arguments in \cite[Page 4913]{kausz}, we again have the identification $Z_{_\alpha} = \text{Isom}(\co^{^N}_{_{T_{_{N,N',\alpha}}}}, \pi_{_*}(\eV_{_\alpha})(N') \mid_{_{T_{_{N,N',\alpha}}}})$, where $\pi:\tt M_{_\alpha} \to T_{_\alpha}$. Thus, $Z_{_\alpha}$ is smooth and surjective over $T_{_{N,N',\alpha}}$ and since the $T_{_{N,N',\alpha}}$ cover $T_{_\alpha}$ for each $\alpha$ we are done. \epr
\brem\label{regularityofquasistuff} The analogues of \cite[Theorem 3.21]{kausz} hold without any serious difficulty. In particular, the deformation theory works to show that $\mathtt Y$ is regular, its generic fibre over $A$ is smooth while its special fibre $\mathtt Y_{_{\tt o}}$ is a divisor with normal crossings. The proof of \eqref{flatavecnc} gets easily adapted to this case. \erem

\subsubsection {Kawamata Coverings}\label{kawa}
Let $X$ be a smooth quasi-projective variety and let $D = \sum_{_{i=1}}^{^r} D_i$
be the decomposition of the simple or reduced normal crossing divisor
$D$ into its smooth components (intersecting transversally).
The "Covering Lemma'' of Y. Kawamata
(see \cite[Lemma 2.5, page 56]{vieweg}, and \cite[Theorem 17]{kawamata}) says
that, given positive integers $N_{_1}, \ldots, N_{_r}$, there is a connected smooth quasi-projective
variety $Z$ over $\bc$ and a Galois covering morphism
\beqa\label{kawamatacm}
\kappa:Z \to  X 
\eeqa
such that the reduced divisor $\kappa^{^*}{D}:= \,({\kappa}^{^*}D)_{_{\text{red}}}$
is a normal crossing divisor on $Z$ and furthermore,
${\kappa}^{^*}D_{_i}= N_{_i}.({\kappa}^{^*}D_{i})_{_{\text{red}}}$. Let $\Gamma$ denote the Galois group
for the covering map $\kappa$.

The isotropy group of any point $z \in Z$, for the
action of $\Gamma$ on $Z$, will be denoted by ${\Gamma}_{_z}$.  It is easy to see that the stabilizer at generic points of the irreducible components of $(\kappa^{^*}D_i)_{_{\text{red}}}$ are cyclic of order $N_{_i}$. By an equivariant principal $G$-torsor $P$ on $Z$ of local type $\tau = \{\tau_{_i}\}_{_{i=1}}^{^r}$ we mean: 
\begin{enumerate}

\item{} The restriction of the $G$-torsor $P_{_{U_{_z}}}$ to an \'etale neighbourhood at a generic point $z$ of an irreducible
component of $(\kappa^{^*}D_i)_{_{\text{red}}}$ is given by a representation $\rho_{_i}:\Gamma_{_z} \to G$;

\item{} for a general point $y$ of an irreducible
component of a ramification divisor for $\kappa$
not contained in $(\kappa^{^*}D)_{_{\text{red}}}$,
the action of ${\Gamma}_{_y}$ on $P$ is the trivial action.

\end{enumerate}
Such a $P$ will always exists as an algebraic space with a $G$-action and can be obtained by gluing trivial $(\Gamma_{_z},G)$-torsors given by $\rho_{_i}$, in $U_{_z}$ for the generic point $z$ of $(\kappa^{^*}D_i)_{_{\text{red}}}$  with pull-backs of $G$-torsors on $X \setminus D$ to $Z$. By a Hartogs type argument, it is easily checked that equivariant $G$-torsors are uniquely defined on $Z$ once given on a subscheme of codimension bigger than $1$.

\section{Appendix to Part II}
\subsubsection{Laced vector bundles}
In this subsection we analyse the special case of {\tt laced} torsors when $G$ is the linear group. Much of the early material in this subsection is adapted from \cite{sesnewstead}. 
\bnot Let ${\text{Vect}}_{_{(\tilde{\cc_{_d}},\bf z)}}^{^{\bf d}}$ denote the category of vector bundles $W$ on $(\tilde{\cc_{_d}}, \bf z)$, i.e. balanced vector bundles on $(\tilde{\cc_{_d}}, \bf z)$ with descent datum (\ref{dpgbundle}) which translates as an isomorphism $V_{_{z_{_1}}} \simeq V^{^*}_{_{z_{_2}}}$. \enot

\noindent
\bdefe\label{balancedparabolicbundle}  A {\tt balanced parabolic structure} on a vector bundle $V$ of rank $n$  on a doubly marked curve $(\tC, \bf c)$ is given by the following datum: 
\begin{enumerate}
\item For $1 \leq s \leq n$, weights, $(\alpha_{_1}, \ldots, \alpha_{_s})$, which are rational numbers such that 
\beqa\label{alpha}
0 \leq \alpha_{_1} < \alpha_{_2} <\cdots < \alpha_{_s} < 1.
\eeqa
and "dual weights'' $(\beta_{_1}, \ldots, \beta_{_s})$:
\beqa\label{dualweights}
 (\beta_{_1}, \ldots, \beta_{_s}) = \left \{\begin{array}{l}  (1-\alpha_{_s},  1-\alpha_{_{s-1}}, \cdots  ,
1-\alpha_{_1})~~{\rm if}~   \alpha_{_1}  \neq   0 \\  
(0, ~~~~  1-\alpha_{_s},  \cdots,1-\alpha_{_2}) ~~~{\rm if}  ~\alpha_{_1}=0
\end{array} \right. 
\eeqa
\item A {\tt balanced parabolic structure} on $V$ at 
$c_{_j}$, $j = 1,2$, i.e.,  {\it strictly} decreasing flags
\beqa
V_{c_{_j}}\, =\, \cf_{c_{_j}}^{^1}\, \supset\, \cf_{c_{_j}}^{^2}\, \supset\, \cdots
\,\supset\, \cf_{c_{_j}}^{^2}\,\supset\, \cf_{c_{_j}}^{^{s+1}}\,=\, 0, ~~j=1,2
\eeqa 
together with weights given as follows:
\begin{itemize}
\item The weight of $\cf_{c_{_1}}^{^m}$ is $\alpha_{_m}$, where $\alpha_{_1},\cdots ,\alpha_{_s}$  as in \eqref{alpha}. 
\item The weight of $\cf_{c_{_2}}^{^m}$ is $\beta_{_m}$, where $\beta_{_1}, \ldots, \beta_{_s}$ are as in \eqref{dualweights}:
\end{itemize}
\end{enumerate}
Let $\text{PVect}^{^{bal}}_{_{(\tC, {\bf c})}}$ denote the category of vector bundes on $(\tC, \bf c)$ with balanced parabolic  structure.
\edefe

 Let $V$ be an object in $\text{PVect}^{^{bal}}_{_{(\tC, {\bf c})}}$.
\begin{enumerate}
\item The flag  $\cf_{c_{_2}}$ and $V_{c_{_2}}$
induces on the  dual $V_{c_{_2}}^*$ of  $V_{c_{_2}}$, the natural  {\it dual flag} $\cf_{c_{_2}}^*$ and the weights  of $\cf_{c_{_2}}^*$
are   "dual''  to   those  of   $\cf_{c_{_2}}$  i.e., they   coincide  with
$\alpha_{_1},\cdots, \alpha_{_s}$, the weights associated to $\cf_{y_{_1}}$ .  

\item For $i = 1,2$, define
\beqa
{\tt gr}(V_{c_{_i}}) := \bigoplus_{_m} {\tt gr}^{^m} \cf_{c_{_i}},~with\\
{\tt gr}^{^m} \cf_{c_{_i}} := \cf_{c_{_i}}^{^m}/\cf_{c_{_i}}^{^{m+1}}.
\eeqa
The graded pieces, ${\tt gr}(V_{c_{_2}}^*)$ gets identified with ${\tt gr}(V_{c_{_2}})$ by a shifting of degrees as
follows:
\beqa 
\left \{\begin{array}{l}
{\tt gr}^{^m} \cf^*_{c_{_2}} = {\tt gr}^{^{s+1-m}} \cf_{c_{_2}}~for~ ~1 \leq m \leq s,~{\rm if}
~ \alpha_{_1} \neq 0 \\
{\tt gr}^{^1} \cf^*_{c_{_2}} = {\tt gr}^{^1}\cf_{c_{_2}}~and~ {\tt gr}^{^m} \cf_{c_{_2}}^* = {\tt gr}^{^{s+2-m}}
\cf_{c_{_2}} ~{\rm if} ~ \alpha_{_1}=0.
\end{array} \right . \eeqa
\end{enumerate}
\bdefe\label{dps}  Let $V$ be an object in $\text{PVect}^{^{bal}}_{_{(\tC, {\bf c})}}$. A {\tt lacing} on $V$ (or more precisely a {\tt s-lacing})  is a {\tt s}-tuple 
\beqa\label{gradedisom}
\wp: = \biggl\{\wp_{_m}: {\tt gr}^{^m} \cf_{c_{_1}} \to {\tt gr}^{^m} \cf^*_{c_{_2}}\biggr\}^{{s}}_{{m = 1}}
\eeqa
of linear isomorphisms. \edefe

\bdefe A balanced parabolic vector bundle endowed with a {\tt lacing} $\wp$ will be called a {\tt laced vector bundle}, i.e., given by the datum:
\beqa
V_{_\wp} := (V_{_\star}, \wp)
\eeqa 
where $V_{_\star}$ is a balanced parabolic bundle on $(\tC, {\bf c}).$ 
\edefe
\bdefe The parabolic degree of a {\tt laced} bunde $V_{_\wp}$ is defined as:
\beqa
{\tt par.deg}_{_{\tC}}(V_{_\wp}):= {\tt par.deg}_{_{\tC}}(V_{_\star}).
\eeqa
\edefe

\blem\label{onpardegs} 
Let $V_{_\wp}$ be a {\tt laced} bundle  on $(\tC, {\bf c})$ and let $k = k_1$ denote the multiplicity of the weight $\alpha_1$. Let $l = (n-k)$.  Then
\beqa\label{pardeg1}
\mbox{${\tt par.deg}_{_{\tC}}$}~V_{_\wp}  = \deg V + (n-k) = \deg V + l. 
\eeqa
{\em As a consequence, the parabolic degree of a {\tt laced} bundle  on $(\tC, {\bf c})$ does not depend on the choice of the parabolic weights}. 
\elem 
\bpr (see \cite{sesnewstead}) By the definition of parabolic degree, 
we see that
\small
\beqa
{\rm  {\tt par.deg}_{_{\tC}}~}V_{_\wp} =  \left \{
\begin{array}{rcl}
\deg V + \sum_{_{m=1}}^{^s} k_{_m} \alpha_{_m} + \sum_{_{m=1}}^{^s} k_{_m}(1-\alpha_{_m})
~{\rm if}~\alpha_{_1} \neq 0 \\
\deg V + \sum_{_{m=2}}^{^s} k_{_m} \alpha_{_m} + \sum_{_{m=2}}^{^s} k_{_m}(1-\alpha_{_m})
~{\rm if}~\alpha_{_1} = 0.\\
\end{array}\right. \eeqa
Hence\small
\beqa
{\rm  {\tt par.deg}_{_{\tC}}}~V_{_\wp} =  \left \{\begin{array}{lcl}
\deg V+n {\rm ~if~} \alpha _1 \neq 0 \\
\deg V + (n-k_1) {\rm ~if~} \alpha_1 = 0.\\
\end{array}\right.
\eeqa
which give the equation \eqref{pardeg1}.
\epr
We summarize the following from \cite{sesnewstead}.
\bprop\label{disting} Let $V_{_\wp}$ be a laced bundle on $\tC$. Then the direct image $\nu_{_*}(V_{_\wp}) = \cf$ is a torsion-free sheaf on $C$ and conversely, $\nu^{^*}(\cf)/{tors}$ recovers the underlying vector bundle of $V_{_\wp}$. Moreover, 
\beqa
{\tt par.deg}_{_{\tC}}(V_{_{\wp}}) = {\tt deg}_{_{C}}(\cf).
\eeqa
\eprop

\subsubsection{Some remarks on parabolic subgroups}
\brem\label{parabtrio}  Let $\lambda:\bg_{_m} \to G$ be a one-parameter subgroup and $P(\lambda)$ be the associated parabolic subgroup and $H(\lambda)$ the Levi quotient which canonically defines a Levi subgroup $L(\lambda)$ as the centralizer of $\lambda$. Let $\eta:G \hra \text{GL}(W)$ be a faithful representation. Then the one-parameter subgroup given by the composition $\eta \circ \lambda:\bg_{_m} \to \text{GL}(W)$ defines a parabolic and Levi subgroups $P(\eta \circ \lambda)$ and $L(\eta \circ \lambda)$ of $\text{GL}(W)$. 

We can view the parabolic subgroup $P(\eta \circ \lambda)$ as the stabilizer of the flag:
\beqa\label{wtdfiltongmodules}
(W_{_\bullet}(\lambda), \epsilon_{_\bullet}) : 0 \subsetneq W_{_1} \subsetneq W_{_2} \ldots W_{_s} \subsetneq W_{_{s+1}} = W\subsetneq W
\eeqa
where $W_{_i}:= \oplus_{j = 1}^{^{i}} W^{^j}$, with $W^{^j}$ being the eigenspace of the $\bg_{_m}$-action via $\lambda$ for the character $z \mapsto z^{^{\gamma_j}}$, and $\gamma_{_1} < \ldots < \gamma_{_{s+1}}$ are the distinct weights which occur. Set $\epsilon_{_i} := (\gamma_{_{i+1}} - \gamma_{_i})/ {\tt dim}(V)$, $i = 1, \ldots, s$. The pair $(W_{_\bullet}(\lambda), \alpha_{_\bullet})$ is called the {\em associated weighted filtration} of $\lambda$. The weighted filtration  $(W_{_\bullet}(\lambda), \alpha_{_\bullet})$ has an associated graded:
\beqa\label{assocgradedlambda}
{\tt gr}_{_\lambda}(W) :=
\bigoplus_{_{j = 1}}^{^{s}} W_{_{j+1}}/W_{_j} = \bigoplus_{_{j = 1}}^{^{s}} W^{^{j+1}}
\eeqa
and it is easy to see that as $H = H(\lambda)$-modules, $W \simeq {\tt gr}_{_\lambda}(W)$. Further,  $H$ fixes the $\lambda$-eigenspaces $W^{^{j+1}}$, i.e., the above decomposition is a decomposition of $H$-modules. We also have an obvious weighted filtration $({\tt gr}(W)^{^{\bullet}}, \ul{\epsilon})$ with the same weights $\ul{\epsilon}$:
\erem

The $1$-PS $\eta \circ \lambda$ also defines a canonical anti-dominant character $\chi_{_{\eta \circ \lambda}}:P(\eta \circ \lambda) \to H(\eta \circ \lambda) \to \bg_{_m}$ dual to $\eta \circ \lambda$ \cite[2.4.9]{schmittbook}. For instance, if $\ul{m} := (m_{_1}, \ldots, m_{_{s+1}})$ is a point of the Levi $H(\eta \circ \lambda)$ as block matrices, then $\chi_{_\lambda}(\ul(m)) := \otimes_{_j} \text{det}(m_{_j})^{^{\gamma_{_j}}}$. This  which restricts to an anti-dominant character $\chi_{_\lambda}$ of $P(\lambda)$. We recall the following result.
\blem\label{schmittbook1} (\cite[Proposition 2.4.9.1]{schmittbook}) Let  $\chi:P(\lambda) \to \bg_{_m}$ be any anti-dominant character. Then there is a positive rational number $r$ such that $\chi = r \cdot \chi_{_\lambda}.$
\elem

\brem\label{parabquadro} Let  $E$ be a $G$-torsor and suppose that we are given  a reduction of structure group $E_{_P} \subset E$ to the parabolic $P$. There is a  canonical anti-dominant character $\chi_{_\lambda}:P \to \bg_{_m}$ (\ref{schmittbook1}) which defines a line bundle $E_{_P}(\chi_{_\lambda})$ on $Y$.

Again, the representation $\eta$ gives a weighted filtration \eqref{wtdfiltongmodules} stabilized by $P(\eta \circ \lambda)$.  We can take the associated vector bundle $E_{_P}(W)$ which comes with its weighted filtration:
\beqa\label{wtdfiltongmodules1}
(E_{_P}(W))_{_\bullet}, \ul{\epsilon}) : 0 \subsetneq E_{_P}(W_{_1}) \subsetneq \ldots E_{_P}(W_{_s}) \subsetneq E_{_P}(W_{_{s+1}}) = E_{_P}(W)
\eeqa
and the weighted slope defined by Schmitt (\cite{schmitt}):
\small
\beqa
L\big((E_{_P}(W)_{_\bullet}, \epsilon_{_\bullet})\big) := \sum_{i = 1}^{^s} \epsilon_{_i}\bigg\{{\tt deg}_{_{C}}(E_{_P}(W))\cdot rk E_{_P}(W_{_{i}}) - {\tt deg}E_{_P}(W_{_{i}}) \cdot rk(E_{_P}(W))\bigg\}.
\eeqa
\normalsize
Claim:
\beqa\label{formalequalityclassiccase}
{\tt deg}(E_{_P}(\chi_{_\lambda})) = L\big((E_{_P}(W)_{_\bullet}, \epsilon_{_\bullet})\big)\\
{\tt deg}(E_{_H}(\chi_{_\lambda})) = L\big((E_{_H}({\tt gr}(W))_{_\bullet}, \epsilon_{_\bullet})\big)
\eeqa
To see this, note that the line bundle $E_{_H}(\chi_{_\lambda}) \simeq \bigotimes \text{det}(E_{_H}(W_{_i}))^{^{-\epsilon_{_i}}}$ with $\epsilon_{_i} := (\gamma_{_{i+1}} - \gamma_{_i})/ {\tt dim}(V)$ as above. 
(see \cite[Exercise 2.4.9.2, page 209]{schmittbook}).
\erem
\brem\label{theabovefortheparabcase} Let $Y$ be a smooth projective curve
and let $\sf P_{_{\ul{v}}}$ be a parahoric group scheme generically split with fibre $G$, with parahoric structures $\ul{v} := (v_{_j})$ at points $y_{_j} \in Y$ given by a tuple of points in the affine apartment $\ca_{_T}$ \cite{base}. Given a faithful representation $\eta:G \hra \text{GL}(W)$, we get a corresponding group parahoric group scheme ${\sf P}_{_{\text{GL}(W)}}$ with generic fibre $\text{GL}(W)$. If $E$ is a $\sf P_{_{\ul{v}}}$-torsor then we get an associated parabolic vector bundle $E(W)_{_*}$ with parabolic structures at $y_{_j}$. If $\lambda:\bg_{_m} \to G$ is a $1$-PS, and the setting be as in the previous paragraph, then we have a parahoric Levi-type torsor $E_{_H}$ for a parahoric group scheme with generic fibre isomorphic to $H$ and associated parabolic line bundles $E_{_H}(\chi_{_\lambda})_{_*}$. The standard properties of degrees of direct sum of vector bundles in terms of the determinants obviously go through in the parabolic setting by replacing degrees with parabolic degrees and tensor products with parabolic tensor products. This follows by expressing parabolic bundles in terms of orbifold bundles and push-forwards. Thus the entire formalism goes through and we get a relation ${\tt par.deg}(E_{_H}(\chi_{_\lambda})_{_*}) = L\big((E_{_H}({\tt gr}(W))_{_\bullet}, \epsilon_{_\bullet})\big)$ with parabolic degrees everywhere.  

We apply it in the main paper for the {\tt laced} bundle $E_{_\wp}$ on $\tC$ which has an underlying parahoric structure at the two points $c_{_i}$.     
\erem
{\subsubsection{ A counter example to a simplistic generalization of Ramanathan's definition in the nodal case}\label{counterexamples} 
Let $E$ be a principal $G$-bundle on the nodal curve $C$. {A na{\"i}ve generalization of the usual definition   along the lines of A. Ramanathan's definition  turns out to be false even when $G = \text{GL}(2).$ \\}

For every maximal parabolic subgroup $P \subset G$ and for every reduction of structure group $E_{_P}$ of $E$ over $C -c$, consider the Lie algebra sub-bundle $E_{_P}(\mathfrak p) \subset E(\mathfrak g)|_{_{C - c}}.$ Let $\overline{E_{_P}(\mathfrak p)}$ be the torsion-free sheaf which is the saturation of the sub-bundle $E_{_P}(\mathfrak p)$ in  $E(\mathfrak g)$ over $C$. The bundle $E$ is ``conjecturally'' (semi)stable if 
\beqa
\text{deg}(\overline {E_{_P}(\mathfrak p)}) < 0 (\leq 0)
\eeqa
For the failure of this ``conjectural definition'' of (semi)stability of $G$-torsors on nodal curves even when $G = \text{GL}(2)$, we give the following counter-example which essentially comes from a remark due to Seshadri.

Let $L,M$ be torsion-free sheaves on $C$ of rank $1$ and degree $0$ which are not locally free. In particular, they are of local type $\mathfrak m$. Consider the group $Ext^{^1}(L,M)$ of extensions of $M$ by $L$. We claim that there is a locally free sheaf $V$ such that:
\beqa\label{Vasexact}
0 \to L \to V \to M \to 0
\eeqa
and hence automatically $V$ is semistable of degree $0$. To see the existence of such a $V$, we consider the local-global spectral sequence for $\text{Ext}$ (\cite[Section 4.2]{tohoku}) which gives (since $ \text{dim}(C) = 1$):
\beqa
\text{H}^{^1}(C, {\mathcal H}om(L,M)) \to \text{Ext}^{^1}(L,M) \to \text{H}^{^0}(C, {\mathcal E}xt^{^1}(L, M)) \to 0.
\eeqa
Note that $\text{H}^{^0}(C, {\mathcal E}xt^{^1}(L, M)) = \text{Ext}^{^1}_{_{A}}(L_{_c}, M_{_c})$, where $A = {\co}_{_{C,c}} \simeq \bc[x,y]/(xy)$. Locally we have ${\mathfrak m} = (x,y)$. Using these as generators, we have an embedding $\mathfrak m \hra \co_{_C} \oplus \co_{_C}$ and hence an extension:
\beqa
0 \to \mathfrak m \to \co_{_C} \oplus \co_{_C} \to \mathfrak m \to 0.
\eeqa
This gives an element in $\text{Ext}^{^1}_{_A}(L_{_c}, M_{_c})$ which lifts to give an element in $\text{Ext}^{^1}(L,M)$. Clearly this extension is locally free since it is so at the node and we get the required $V.$ This $V$ is semistable of degree $0$.

Giving a reduction of structure group of the principal $\text{GL}(2)$-bundle underlying $V$ is expressing it in an exact sequence of vector bundles \eqref{Vasexact}
and {\em the conjectural  definition of semistability} is equivalent to saying that for the
sub-bundle $L \otimes M^*\subset V \otimes V^*$, we have
\beqa
\text{deg}(\overline{L \otimes M^*}) \leq 0
\eeqa
where $\overline{L \otimes M^*}$ denotes the saturation in $V \otimes V^*$.

\noindent
{\it Claim}: 
\beqa
\text{deg}(\overline{L \otimes M^*}) = 1.
\eeqa
In particular $V \otimes V^*$  is not semistable.  Let $L'$
(resp. $M'$) denote $p^{^*}(L)/tors$ (resp. $p^{^*}(M)/tors$).  Then  the line sub-bundle of $p^{^*}(V)$ (resp. $p^{^*}(V^{^*})$) generated by $L'$
(resp. $M'$) is of the form $L'(y_{_1}+y_{_2})$ (resp. $M'(y_{_1}+y_{_2}))$.
We have $\text{deg}~L' = \text{deg}~M' = -1$, so that $\text{deg}~ L' (y_{_1}+y_{_2}) = \text{deg}~
M'(y_{_1}+y_{_2})=1$.  Then we see that the line bundle
\beqa
N = (L'(y_{_1}+y_{_2}) \otimes M' (y_{_1}+y_{_2}) (-y_{_1}- y_{_2})
\eeqa
descends to a torsion free subsheaf of $V \otimes V^*$, which is the saturation $\overline{L \otimes M^*}$.  Since $\deg N =0$, we see that $\text{deg}(\overline{L \otimes M^*}) = 1$. 
\brem The lesson is to avoid taking the saturation after taking tensor products. The degree exceeds the bound. Instead, one has to take some sort of a ``parabolic tensor product'' and then take a saturation, both of these operations need to be carried out on the normalization $Y$. This can be made precise. We proceed differently in \S 13 to achieve this.\erem 
\epigraph{\footnotesize{\it{AbhyavasthāH prajāyante pra vavrer vavriś ciketa, upasthe mātur vi cashte}}\\ States upon states are born, covering over covering awakens to knowledge, in the lap of the universal mother he wholly sees.}{\textit{Rig Veda, Mandala V, Hymn 19.1}}

\end{document}